\documentclass[twoside,a4paper,11pt]{article}
\usepackage[top=4cm, bottom=4cm, left=3cm, right=3cm]{geometry}

\usepackage{amsmath,amssymb,amsthm,enumerate}
\usepackage{xcolor}
\usepackage{hyperref}
\usepackage{graphicx}
\hypersetup{
    pdfborder={0 0 0.5 [3 3]},
    colorlinks=true,
    linkcolor=blue,
    filecolor=magenta,      
    urlcolor=cyan,
}

\usepackage{mathtools}
\newtheorem{theorem}{Theorem}[section]
\newtheorem{proposition}[theorem]{Proposition}
\newtheorem{lemma}[theorem]{Lemma}
\newtheorem{corollary}[theorem]{Corollary}

\makeatletter
\def\th@definition{%
  \thm@headfont{\bfseries}
  \thm@notefont{
  \fontseries{m}\fontshape{it}\selectfont}
  \itshape
}
\makeatother

\theoremstyle{definition}
\newtheorem{definition}[theorem]{Definition}

\theoremstyle{remark}
\newtheorem{remark}[theorem]{Remark}

\newcommand{\CAT}{\mathrm{CAT}}
\newcommand{\X}{\mathcal{X}}

\newcommand{\E}{\mathbb E}

\newcommand{\R}{\mathbb R}
\newcommand{\Sp}{\mathbb S}
 
\newcommand{\dist}{\mathrm{d}}
\newcommand{\Wone}{W_1}
\newcommand{\supp}{\operatorname{supp}}
\newcommand{\per}{\mathrm{per}}
\newcommand{\DS}{\displaystyle}
\newcommand{\dd}{\mathrm d}

\newcommand{\ybcomment}[1]{\par\noindent{\color{red}\textbf{YB:} #1}\par}
\newcommand{\bary}{\mathrm{bar}}
\newcommand{\conv}{\mathrm{conv}}

\newcommand{\dH}{\mathrm{d}_{\mathrm H}}
\newcommand{\diff}{\mathrm{d}}
\newcommand{\eps}{\epsilon}
\newcommand{\arccosh}{\mathrm{arccosh}}

\newcommand{\veb}[1]{\textcolor{magenta}{\textbf{VEB:}  #1}}

\renewcommand{\geqslant}{\geq}
\renewcommand{\leqslant}{\leq}

\begin{document}

\title{Geometry and barycenters in large $\CAT(\kappa)$ spaces}

\author{Yassine Boukhateb
\and Victor-Emmanuel Brunel
\and Austin J.\ Stromme
\thanks{\raggedright
Department of Statistics, CREST, ENSAE, IP Paris.
Emails: \texttt{firstname.lastname@ensae.fr}.}}

\date{}
\maketitle







\begin{abstract}
\emergencystretch=1em
We study the geometry and resulting properties of Fréchet functions and barycenters in metric spaces with curvature bounded above by a positive constant in the sense of Alexandrov. Existing results typically rely on the assumption that the support is small enough that the squared distance, and thus the Fréchet function, is strongly geodesically convex. We work in the sharper regime where existence and uniqueness of the barycenter are guaranteed, but where convexity typically fails. We develop new geometric properties of these spaces, including a quantitative four-point inequality and a geodesic stability result. Using these tools,
we establish new finite sample error bounds for two statistical estimators: the empirical barycenter and the iterated barycenter, as well as new results on sparse representations and approximations of barycenters.
\end{abstract}





\tableofcontents

 \section{Introduction}
 
 Given a metric space $(\X, \dist)$ and a Borel probability distribution $\mu$ on $\X$
 with finite second moment, a barycenter of $\mu$
 is defined to be any minimizer $b(\mu)$ of the Fréchet function
 $$
 F_\mu:b\in \X\mapsto \frac12 \int_\X \dist^2(b, y) \diff \mu(y).
 $$
 Barycenters generalize the statistical mean to non-Euclidean
 spaces, and are therefore a fundamental object connecting statistics
 to geometry~\cite{frechet1948,karcher1977}.
 
 From the statistical point of view, the most fundamental question one can
 ask is: how can we estimate barycenters given only samples from $\mu$?
 Concretely, suppose $X_1, \ldots, X_n \sim \mu$ are independent draws
 from $\mu$. 
 Perhaps the most natural approach to estimating a population barycenter $b(\mu)$ is
 to use a barycenter associated to the empirical measure $\mu_n := \frac1n \sum_{i = 1}^n \delta_{X_i}$. A barycenter $b(\mu_n)$ of $\mu_n$ is called an {\it empirical barycenter}.
 Another approach is to iteratively compute $b_{k + 1} = \gamma_k(t_k)$
 where $\gamma_k$ is a geodesic between $b_k$ and $X_{k + 1}$;
 the result $b_n$ is called the {\it iterative barycenter}~\cite{sturm2003}.
 Once these barycenters have been defined, it is natural to ask 
 how well they estimate $b(\mu)$.
The performance of these estimators
is intimately connected to the geometry of the underlying
space $(\X, \dist)$.
For example, one can equivalently characterize $\CAT(0)$
spaces through
the strong geodesic convexity of the Fréchet function $F_\mu$ over all
probability measures $\mu$ with finite second moment~\cite{sturm2003}.
Strong convexity of $F_\mu$ then yields, via simple proofs,
that there is a unique barycenter and that 
the empirical barycenter and the iterated barycenter
both estimate the true barycenter at a fast $1/n$ rate~\cite{sturm2003,legouic2022}.

 However, when one relaxes the $\CAT(0)$ assumption to 
 $\CAT(\kappa)$ for $\kappa > 0$,
 the behavior of barycenters and the Fréchet function becomes much more delicate.
 As can be immediately seen for the uniform measure on two antipodal points on a sphere,
 the Fréchet function can be non-convex, and barycenters may not be unique.
 Nonetheless, it is known that if one restricts the domain to what we here call a {\it $\kappa$-domain},
 namely a subset of a closed ball $B(o, \pi/(2\sqrt{\kappa}) - \eps)$
 for $\eps > 0$, then barycenters are guaranteed to be unique~\cite{afsari2011,yokota2016}, despite
 the fact that the Fréchet function may be non-convex.
 
 Because $\kappa$-domains
 are the largest sets in $\CAT(\kappa)$-spaces
 where barycenters are necessarily unique, they form
 a natural setting for the study of barycenters and Fréchet functions.
 However, the possible non-convexity of the Fréchet function
 obstructs a straightforward adaption of the techniques
 from the $\CAT(0)$ case.
 This obstruction has limited previous analyses of both empirical
 and iterated barycenters,
 either
 to significantly smaller domains~\cite{legouic2022,brunelserres2024,brunelserres2025},
or to ``slow" rates of the form $n^{-\beta}$, where
$\beta$ is a parameter that depends on the metric entropy of the space~\cite{schotz2019convergence,ahidar2020}.

The purpose of this work is therefore
to develop improved
understanding of the geometry of $\kappa$-domains, with
the goal of providing optimal rates
for barycenter estimation on the largest possible domains.
Our main results for empirical barycenters and iterated barycenters are given
in Section~\ref{sec:consequences}.
In particular, Theorem~\ref{th:empirical-expected-error}
and Corollary~\ref{cor:empirical-high-probability}
give $1/n$ rates in expectation and high probability, respectively,
for empirical barycenters of distributions supported on $\kappa$-domains. 
For iterated barycenters, we show $1/n$
rates in expectation and high probability in Theorem~\ref{th:iterated-expected-error}
and Corollary~\ref{cor:empirical-high-probability}.
As we discuss in Remark~\ref{rem:tightness},
all of these results are optimal in the sense that
$\kappa$-domains are the largest domains where it is possible
to get $1/n$ rates,
uniformly over probability distributions $\mu$. 
Finally, we obtain new
Carathéodry-type results
on $\CAT(\kappa)$ spaces in Theorem~\ref{thm:Carathéodory},
and we
study sparse approximations
of barycenters in Theorem~\ref{th:sparse_approximation}.

To better understand the geometry of
$\kappa$-domains, we develop
several new results on the geometry of $\kappa$-domains.
We highlight here a few of these new results which may be
of particular interest.
The first of these results, Proposition~\ref{prop:local-four-point} is a new quantitative four-point
inequality valid on $\kappa$-domains, generalizing
Reshetnyak's four-point comparison.
This four-point inequality is useful in our high probability
bounds for empirical barycenters in Corollary~\ref{cor:iterated-high-probability}.
For our bounds in expectation, the main tool is a new bound
on the optimality gap of the Fréchet function Theorem~\ref{thm:directional-derivative-transfer},
which gives a number of useful consequences, such as, for the first time, a Polyak-\L{}ojasiewicz inequality
for arbitrary distributions on $\kappa$-domain in Theorem~\ref{thm:pl-inequality},
as well as a more geometric version of Yokota's proof~\cite{yokota2016}
of existence and uniqueness of barycenters
in $\kappa$-domains.
Our results on the iterated barycenter are based on
different tools, in particular a new weighted distance
in Definition~\ref{def:weighted-distance},
and a resulting stability result in 
Proposition~\ref{prop:endpoint-inequality}.

\paragraph*{Organization of the paper.}
In Section~\ref{sec:background}, we recall background on $\CAT(\kappa)$ spaces.
In Section~\ref{sec:new_geometric_results}, we present our new geometric results on the quantitative four-point inequality and geodesic stability.
In Section~\ref{sec:geometric_properties_Frechet}, we give our main results on Fréchet functions, including a new version of Yokota's proof of barycenter existence
and uniqueness based on our new geometric insights from Section~\ref{sec:new_geometric_results}.
Finally, in Section~\ref{sec:consequences}, we give our results on empirical barycenters, iterated barycenters, and Carathéodory's theorem.

\paragraph*{Concurrent work of C. Schötz.}
In the final stages of preparing this work, two works of Schötz~\cite{schotz2026trapezoid,schotz2026weighted} appeared on the arXiv, each partially overlapping with our own. In particular, these works contain results related to our quantitative four-point inequality, variance inequality, and bounds for empirical barycenters.

For the quantitative four-point inequality, which appears here as Proposition~\ref{prop:local-four-point}, Schötz uses different techniques than us to obtain an inequality
with an optimal constant in terms of his diameter parameter in \cite[Theorem~1.4]{schotz2026trapezoid}.
Regarding the variance inequality, which appears here as Proposition~\ref{prop:quadratic-variance-cat-k},
Schötz uses a similar argument as ours to obtain sharper constants in \cite[Theorem~A.2]{schotz2026weighted}.
Lastly, bounds in expectation and with high probability for empirical barycenters on $\kappa$-domains appear here as Theorem~\ref{th:empirical-expected-error} and Corollary~\ref{cor:empirical-high-probability}, and in Schötz's work as \cite[Theorem~4.3 and Corollary~4.6]{schotz2026weighted}.

Our analysis of empirical barycenters in expected-error uses a tangent-cone argument, whereas Schötz uses a leave-one-out argument.
For high-probability bounds, both papers use the bounded differences inequality, as in~\cite{brunelserres2025}, and take the quantitative four-point inequality as input through barycenter stability. The refined bound in \cite[Remark~4.2 and equation~(31)]{schotz2026weighted}
gives a better dependence on the radius as it tends to $\pi/(2\sqrt{\kappa})$.
On the other hand, Schötz's expected-error bounds involve the metric variance $\frac12\E[\dist^2(X,X')]$, where $X,X'$ are independent with law $\mu$. Their conversion to the Fréchet variance $\sigma_\mu^2$ in \cite[equation~(35)]{schotz2026weighted} gives a coefficient tending to $2$ as $\kappa\downarrow0$.
Our bound instead recovers the optimal constant $1$ relative to $\sigma_\mu^2$.

\paragraph*{AI Usage Statement.} We occasionally used AI as a reference for 
standard trigonometric formulas, as well as for help with figures, translation and finding references.

 

\paragraph*{Notation.}

Let us set up the basic notation that we use throughout this article. Consider a metric space $(\X,\dist)$.  For all $x\in\X$ and $r\geq 0$, the closed ball with center $x$ and radius $r$ is denoted by $B(x,r)$. For $x\in \X$ and a non-empty subset $A\subseteq \X$, we let $\dist(x,A)=\inf_{y\in A}\dist(x,y)$ be the distance from $x$ to the set $A$.  The closure of a subset $A\subseteq \X$ is denoted by $\overline A$.
The collection of all Borel probability measures on $\X$ is denoted by $\mathcal P(\X)$.
For $k\geq 1$, we let $\mathcal P_k(\X)$ be the collection of all Borel probability measures on $\X$ such that $\int_\X \dist(x,x_0)^k\diff\mu(x)<\infty$ for some (and hence, all) $x_0\in\X$.
 Given $x \in \X$, we write the Dirac
    measure supported at $x$ as $\delta_x$.
 Given $\mu\in\mathcal P(\X)$, we denote its support by $\supp(\mu)$. Recall that it is defined as $$\supp(\mu):=\left\{x\in \X : \mu\bigl(B(x,\varepsilon)\bigr)>0 \text{ for every }\varepsilon>0\right\}.$$ This is a closed set. If $\X$ is separable, then $\mu(\supp(\mu))=1$ and $\supp(\mu)$ is the smallest closed subset of $\X$ with full $\mu$-measure (see \cite[Chapter II, Theorem 2.1]{parthasarathy1967probability}). Throughout, saying that a probability measure $\mu$ is supported on a Borel set $D$ means that $\mu(D)=1$.
 Finally, for $p\geqslant 1$,
 we define the $W_p$ Wasserstein
 distance between probability measures $\mu, \nu \in \mathcal{P}_p(\X)$ as
 $$
 W_p^p(\mu, \nu) := \inf_{\pi \in \Pi(\mu, \nu)}
 \int \dist^p(x, y) \dd \pi(x, y),
 $$ where $\Pi(\mu, \nu)$ denotes all of the couplings
 of $\mu$ and $\nu$.
 
\section{Background on \texorpdfstring{\boldmath$\CAT(\kappa)$}{CAT(kappa)} spaces}
\label{sec:background}
 In this section, we recall standard material
 on metric and Alexandrov geometry,
 following~\cite{burago2001course} and~\cite{bridsonhaefliger1999}.
 
\subsection{Geodesic metric
spaces and absolutely continuous curves}
 

Let $(\X,\dist)$ be a metric space. A {\it geodesic} is a map $\gamma \colon [0,1] \to \X$
such that $\dist(\gamma(s), \gamma(t)) = |s- t| \dist(\gamma(0), \gamma(1))$
for all $s,t\in[0,1]$. Given $x,y\in \X$, the collection of all geodesics $\gamma$ with $\gamma(0)=x$ and $\gamma(1)=y$ is denoted by $\Gamma(x,y)$. When $\Gamma(x,y)\neq\emptyset$ for all $x,y\in \X$, $(\X, \dist)$
is called a {\it geodesic space}.

A subset $C\subseteq \X$ is called geodesically convex if and only if $\gamma(t)\in C$ for all $t\in[0,1]$ and all geodesics $\gamma$ with $\gamma(0),\gamma(1)\in C$. 
Let $C\subseteq \X$ be geodesically convex and $f:C\to \R$. We say that $f$ is geodesically convex if $f(\gamma(t))\leq (1-t)f(\gamma(0))+tf(\gamma(1))$, for all $t\in[0,1]$ and all geodesics $\gamma$ with $\gamma(0),\gamma(1)\in C$. Given $\alpha>0$, we say that $f$ is geodesically $\alpha$-strongly convex if $f(\gamma(t))\leq (1-t)f(\gamma(0))+tf(\gamma(1))-\frac{\alpha}{2}t(1-t)\dist(\gamma(0),\gamma(1))^2$, for all $t\in[0,1]$ and all geodesics $\gamma$ with $\gamma(0),\gamma(1)\in C$.

Given $A\subseteq \X$, its geodesic convex hull, denoted by $\conv(A)$, is the smallest geodesically convex subset of $\X$ containing $A$. Equivalently, it is the intersection of all geodesically convex subsets of $\X$ containing $A$.



\begin{definition}[Absolute continuity and metric speed]
\label{def:metric-speed}
A curve $\gamma:[0,1]\to\X$ is absolutely continuous if there exists a nonnegative function $g\in L^1(0,1)$ such that
\[
\dist(\gamma(s),\gamma(t))
\leq\int_s^t g(u)\,\dd u,
\qquad 0\leq s\leq t\leq1.
\]
For such a curve, by \cite[Theorem~1.1.2]{ambrosio2005gradient}, the ratio $\dist(\gamma(s),\gamma(t))/|t-s|$, with $s\in[0,1]\setminus\{t\}$, has a limit as $s\to t$ for almost every $t\in(0,1)$ with respect to Lebesgue measure. We denote this limit by $|\Dot\gamma|(t)$ and call it the metric speed of $\gamma$ at $t$. Moreover, $|\Dot\gamma|\in L^1(0,1)$ and
\[
\dist(\gamma(s),\gamma(t))
\leq\int_s^t |\Dot\gamma|(u)\,\dd u,
\qquad 0\leq s\leq t\leq1.
\]
For $x,y\in\X$, we denote by $A(x,y)$ the set of all absolutely continuous curves $\gamma:[0,1]\to\X$ with $\gamma(0)=x$ and $\gamma(1)=y$. In particular, $\Gamma(x,y)\subseteq A(x,y)$.
\end{definition}

\subsection{Basics of upper-bounded Alexandrov curvature}
We now briefly recall the definition of metric spaces with global curvature upper bounds in Alexandrov's sense. For more details, we refer the reader to \cite{burago2022course,Alexandrovgeom} and \cite{bridson2013metric}.
First, we recall from \cite{brunelserres2025} the definition of model spaces of constant curvature. Fix $\kappa\in\R$.

If $\kappa = 0$, then $M_0 = \R^2$ is the Euclidean plane.
This is a geodesic space where geodesics are unique and given by line segments.

If $\kappa > 0$, then $M_\kappa=\frac{1}{\sqrt\kappa}\Sp^2$ is 
the $2$-sphere of radius $\frac{1}{\sqrt{\kappa}}$.
We have the closed
form of distances based on the arc length metric: $\dist_\kappa(x,y)=\frac{1}{\sqrt\kappa}\arccos(\kappa x^\top y)$, for all $x,y\in M_\kappa$. Geodesics are unique except for antipodal points, and given by arcs of great circles.

If $\kappa < 0$, then $M_\kappa=\frac{1}{\sqrt{-\kappa}}\mathbb H^2$ is
the hyperbolic space.
We can define $\mathbb H^2=\{(x_1,x_2,x_3)\in\R^3:x_3>0, x_1^2+x_2^2-x_3^2=-1\}$. 
Distances are given by $\dist_\kappa(x,y)=\frac{1}{\sqrt{-\kappa}}\textrm{arccosh}(\kappa\langle x,y\rangle)$, for all $x,y\in M_\kappa$, where $\langle x,y\rangle=x_1y_1+x_2y_2-x_3y_3$. This is a geodesic space where geodesics are always unique and are given by the intersections of $M_\kappa$ with planes going through the origin in $\R^3$.



The Alexandrov curvature of a geodesic metric space $(\X, \dist)$ is defined by comparing triangles
in $(\X, \dist)$ to corresponding triangles in the model spaces.
	A {\it geodesic triangle} in $\X$ is a set of three points in $\X$, the vertices, together with three geodesics connecting them, the sides.
	
For $\kappa\in\R$, denote by $D_\kappa$ the diameter of the model space $M_\kappa$. That is, 
$$D_\kappa=\begin{cases}
    \frac{\pi}{\sqrt\kappa} & \text{ if } \kappa>0 \\
    \infty & \text{ if } \kappa \leq 0. 
\end{cases}$$


Given three points $x,y,z\in \X$, we write $\Delta(x,y,z)$ to denote a triangle with vertices $x,y,z$, with no mention of which geodesics are chosen for the sides; these geodesics
may not necessarily be unique. The perimeter of a triangle $\Delta=\Delta(x,y,z)$ is defined as $\per(\Delta)=\dist(x,y)+\dist(y,z)+\dist(x,z)$.
Let $\kappa\in\R$ and $\Delta$ be a triangle in $\X$ with $\per(\Delta)<2D_\kappa$. A {\it comparison triangle} for $\Delta$ in the model space $M_\kappa$ is a triangle $\bar\Delta\subseteq M_\kappa$ with the same side lengths as $\Delta$.



Note that in $M_\kappa$, any side of a triangle with perimeter less than $2D_\kappa$ must be of length less than $D_\kappa$, and so the geodesics connecting the vertices are always unique. Moreover, given any positive numbers $a,b,c$ with $a\leq b+c$, $b\leq a+c$, $c\leq a+b$ and $a+b+c<2D_\kappa$, there exists a unique triangle in $M_\kappa$ with side lengths given by $a,b$ and $c$, up to rigid transformations. Therefore, comparison triangles are always unique up to isometries of the model space $M_\kappa$.

We may now define upper-bounded Alexandrov curvature. Intuitively, we say that $(\X, \dist)$ has global curvature bounded from above by $\kappa$ if all its triangles with perimeter smaller than $2D_\kappa$ are thinner than their comparison triangles in the model space $M_\kappa$.

\begin{definition}[Global curvature upper bounds and $\CAT(\kappa)$ spaces] \label{def:CAT}
	Let $(\X, \dist)$ be a metric space and $\kappa\in\R$.
    We say that $(\X, \dist)$ has {\it global curvature bounded from above by $\kappa$} if and only if for all triangles $\Delta\subseteq \X$ with $\per(\Delta)<2D_\kappa$ and for all $x,y\in \Delta$, $\dist(x,y)\leq \dist_\kappa(\bar x,\bar y)$, where $\bar x$ and $\bar y$ are the points on a comparison triangle $\bar\Delta$ in $M_\kappa$ that correspond to $x$ and $y$ respectively.
    We say that {\it $(\X, \dist)$ is a $\CAT(\kappa)$ space}
    if it is a geodesic space, complete, and has global curvature bounded from above by $\kappa$ (note that some authors do not require completeness in the definition of $\CAT(\kappa)$ spaces, but we do).
    We say that $(\X, \dist)$ is a $\CAT$ space if it is a $\CAT(\kappa)$ space for some $\kappa\in\R$.
\end{definition}
Importantly, by \cite[Theorem II.1.12]{bridsonhaefliger1999}, any $\CAT(\kappa)$ space is also a $\CAT(\kappa')$ space, for all $\kappa,\kappa'\in\R$ with $\kappa\leq \kappa'$. 
Moreover, it is clear that for all $\kappa\in\R$ and $\rho>0$, a metric space $(\X, \dist)$ is a $\CAT(\kappa)$ space if and only if the $\rho$-dilation $(\X,\rho\,d)$ is a $\CAT(\kappa/\rho^2)$ space. In particular, if $\kappa>0$ and $(\X,\dist)$ is a $\CAT(\kappa)$ space, then $(\X,\sqrt\kappa \dist)$ is a $\CAT(1)$ space. We will frequently
use this fact to simplify our proofs by assuming
$\kappa=1$ without loss of generality.

We can now formally define the subsets of the $\CAT(\kappa)$ spaces
that we study in this work.
\begin{definition}[$\kappa$-domain] 
A {\it $\kappa$-domain} is
any closed ball $B(o, r)$ with $r <\pi/(2\sqrt\kappa)$ if $\kappa>0$, and any ball or $\X$ itself if $\kappa\leq 0$. 
\end{definition}


An essential property of $\kappa$-domains is the following.

\begin{proposition} \cite[Proposition II.1.4]{bridsonhaefliger1999}
\label{prop:unique-geodesics-domain}
    Let $(\X,\dist)$ be a $\CAT(\kappa)$ space for some $\kappa\in\R$. Then, between any two points belonging to some $\kappa$-domain, there is a unique geodesic, and this geodesic does not leave the $\kappa$-domain. 
\end{proposition}


\subsection{Angles in \texorpdfstring{\boldmath$\CAT(\kappa)$}{CAT(kappa)} spaces}

Throughout this work, we use the
notion of comparison angles in $\CAT(\kappa)$ spaces.

\begin{definition}[$\kappa$-comparison angle] \label{def:comparison_angles} Let $(\X,\dist)$ be a $\CAT(\kappa)$ space for some $\kappa\in\R$. Let $x,y,z\in\X$ with $x\notin\{y,z\}$ satisfying $\dd(x,y)+\dd(y,z)+\dd(x,z)<2D_\kappa$.
The $\kappa$-comparison angle $\angle_x^{(\kappa)}(y, z)$ is the angle at $\bar x$ in the comparison triangle $(\bar x,\bar y,\bar z)$ in $M_\kappa$. That is, by the spherical/hyperbolic law of cosines \cite[I.2.2 and I.2.7]{bridsonhaefliger1999},
\begin{align*}
& \angle_x^{(\kappa)}(y,z) := \\
    & \qquad \begin{cases}
    \arccos\left(
    \dfrac{
        \cos\bigl(\sqrt{\kappa}\,\dist(y,z)\bigr)
        -\cos\bigl(\sqrt{\kappa}\,\dist(x,y)\bigr)
         \cos\bigl(\sqrt{\kappa}\,\dist(x,z)\bigr)
    }{
        \sin\bigl(\sqrt{\kappa}\,\dist(x,y)\bigr)
        \sin\bigl(\sqrt{\kappa}\,\dist(x,z)\bigr)
    }
    \right)
    & \text{if } \kappa>0,
    \\[1.2em]
    \arccos\left(
    \dfrac{
        \dist(x,y)^2+\dist(x,z)^2-\dist(y,z)^2
    }{
        2\dist(x,y)\dist(x,z)
    }
    \right)
    & \text{if } \kappa=0,
    \\[1.2em]
    \arccos\left(
    \dfrac{
        \cosh\bigl(\sqrt{-\kappa}\,\dist(x,y)\bigr)
        \cosh\bigl(\sqrt{-\kappa}\,\dist(x,z)\bigr)
        -\cosh\bigl(\sqrt{-\kappa}\,\dist(y,z)\bigr)
    }{
        \sinh\bigl(\sqrt{-\kappa}\,\dist(x,y)\bigr)
        \sinh\bigl(\sqrt{-\kappa}\,\dist(x,z)\bigr)
    }
    \right)
    & \text{if } \kappa<0.
\end{cases}
\end{align*}

In particular, $\angle_x^{(\kappa)}(y,y)=0$ whenever $y\neq x$ and the perimeter condition holds. One can easily check that the ratios inside the $\arccos$ are indeed always in $[-1,1]$, thanks to the triangle inequality and the perimeter assumption.
Note as well that,
for fixed distinct points $x,y,z$, the map $\kappa\mapsto\angle_x^{(\kappa)}(y,z)$ is continuous on the set of parameters satisfying the perimeter condition.

\end{definition}
\begin{remark}[Euclidean comparison angle]\label{rem:euclidean-angle}
 When $\kappa=0$, $\angle_x^{(0)}(y,z)$ is called the Euclidean comparison angle, and it is defined in any metric space provided $x\notin\{y,z\}$, without a
perimeter restriction.
\end{remark}
\begin{definition}[Alexandrov angle {\cite[Definition~I.1.12 and Proposition~II.3.1]
{bridsonhaefliger1999}}]
\label{def:alexandrov-angle}
Let $(\X,\dist)$ be a $\CAT(\kappa)$ space, with $\kappa\in\R$. For two non-constant geodesics $\gamma,\eta$ issuing from $x\in\X$, their Alexandrov angle at $x$ is defined by
\[
\angle_x(\gamma,\eta)
:=
\lim_{s,t\downarrow0}
\angle_x^{(0)}(\gamma(s),\eta(t)).
\]
The limit exists by the $\CAT(\kappa)$ property.

For three distinct points $x,y,z\in\X$ and, if $\kappa>0$, assuming $\dist(x,y),\dist(x,z)<D_\kappa$, let $\gamma\in\Gamma(x,y)$ and $\eta\in\Gamma(x,z)$ be the unique geodesics. We define
\[
\angle_x(y,z):=\angle_x(\gamma,\eta).
\]
We also set $\angle_x(y,y):=0$ whenever $y\neq x$ and the geodesic from $x$ to $y$ is unique.
\end{definition}

The Alexandrov angle depends only on the initial directions of the two geodesics. Thus, with the notation above, $\angle_x(y,z)=\angle_x(\gamma(s),\eta(t))$ for every $s,t\in(0,1]$. Moreover, the Euclidean comparison angle in its definition can be replaced by a comparison angle of any fixed curvature $\kappa'\in\R$:
\begin{equation}\label{eqn:def_AlexandrovAngle}
\angle_x(\gamma,\eta)=\lim_{s,t\downarrow0}\angle_x^{(\kappa')}(\gamma(s),\eta(t)).
\end{equation}
Indeed, for sufficiently small $s,t>0$, the perimeter condition is satisfied. If $\gamma(s)=\eta(t)$, we take the comparison angle to be zero.


In $\CAT(\kappa)$ spaces, Alexandrov and $\kappa$-comparison angles are comparable in the following sense.

\begin{proposition}[Angle comparison {\cite[Proposition~II.3.1]{bridsonhaefliger1999}}]
\label{prop:angle-comparison}
Let $(\X,\dd)$ be a $\CAT(\kappa)$ space, with $\kappa>0$. Let $x,y,z\in\X$. Assume that $x\notin\{y,z\}$ and that $\dd(x,y)+\dd(y,z)+\dd(x,z)<2D_\kappa$. Then,
\[
\angle_x^{(\kappa)}(y,z)\geq\angle_x(y,z).
\]
\end{proposition}

In most applications below, we compare the cosines of angles rather than the angles themselves. The resulting algebraic inequality extends to points in a $\kappa$-domain without any additional perimeter assumption, even when the comparison angle is not defined. 
\begin{proposition}[Cosine comparison in a $\kappa$-domain]
\label{prop:cosine-comparison-domain}
Let $(\X,\dist)$ be a $\CAT(\kappa)$ space, with $\kappa>0$, and let $D$ be a $\kappa$-domain. For $x,y,z\in D$, with $x\notin\{y,z\}$,
\begin{equation}
\label{eq:cosine-angle-comparison}
\frac{
\cos\bigl(\sqrt{\kappa}\,\dist(y,z)\bigr)
-
\cos\bigl(\sqrt{\kappa}\,\dist(x,y)\bigr)
\cos\bigl(\sqrt{\kappa}\,\dist(x,z)\bigr)
}{
\sin\bigl(\sqrt{\kappa}\,\dist(x,y)\bigr)
\sin\bigl(\sqrt{\kappa}\,\dist(x,z)\bigr)
}
\leq \cos\angle_x(y,z).
\end{equation}
\end{proposition}

For the reader's convenience, we
give a proof of Proposition~\ref{prop:cosine-comparison-domain} in Appendix~\ref{subsec:cosine_comparison_without_perimeter_restriction}.

Finally, we shall have use of the following inequality
involving cosines along geodesics. A proof is given
in Section~\ref{subsec:Proof_cos_ineq}. 

\begin{proposition}[Cosine inequality along geodesics]
\label{prop:geodesic-cosine-inequality}
Let $(\X,\dd)$ be a $\CAT(\kappa)$ space, with $\kappa>0$, and let $D$ be a $\kappa$-domain. Let $x,y,z\in D$ and $\gamma\in\Gamma(y,z)$.
Set
\[
\alpha_t:=
\frac{\sin((1-t)\sqrt{\kappa}\,\dd(y,z))}
     {\sin(\sqrt{\kappa}\,\dd(y,z))},
\qquad
\beta_t:=
\frac{\sin(t\sqrt{\kappa}\,\dd(y,z))}
     {\sin(\sqrt{\kappa}\,\dd(y,z))}.
\]
When $y=z$, these coefficients are understood as
$\alpha_t=1-t$ and $\beta_t=t$.
Then, for every $t\in[0,1]$,
\[
\cos\bigl(\sqrt{\kappa}\,\dd(x,\gamma(t))\bigr)
\geq
\alpha_t\cos\bigl(\sqrt{\kappa}\,\dd(x,y)\bigr)
+
\beta_t\cos\bigl(\sqrt{\kappa}\,\dd(x,z)\bigr).
\]
\end{proposition}



\subsection{Tangent cones in \texorpdfstring{\boldmath$\CAT(\kappa)$}{CAT(kappa)} spaces}
We shall have frequent use of tangent cones, which
are a generalization of tangent spaces to
non-smooth geometries, and here recall their basic
definition and properties.
\begin{definition}[Tangent cone and log function {\cite[Definition~II.3.18]{bridsonhaefliger1999}}]\label{def:tangent-cone}
Let $(\X,\dd)$ be a $\CAT(\kappa)$ space for some $\kappa\in\R$ and let $x\in \X$. Let $\Gamma_x=\bigcup_{y\in \X\setminus\{x\}}\Gamma(x,y)$ be the set of all non-constant geodesics issuing from $x$. We call any two pairs $(s,\gamma),(t,\eta)\in [0,\infty)\times \Gamma_x$ equivalent if and only if, ($s=t=0$) or ($s=t>0$ and $\angle_x(\gamma,\eta)=0$. We denote by $\tilde C_x\X$ the quotient space. For $a\geq0$ and $\gamma\in\Gamma_x$, we write $a[\gamma]$ for the equivalence class of $(a,\gamma)$ and if $a=0$, we simply write it as $0_x$, which is the vertex of the cone $\tilde C_x\X$. We equip $\tilde C_x\X$ with the distance 
\[
\dist_{\widetilde C_x\X}(s[\gamma],t[\eta])
:=
\left(
s^2+t^2-2st\cos\angle_x(\gamma,\eta)
\right)^{1/2},
\]
which does not depend on the choice of $\gamma$ and $\eta$ in the equivalence classes.

The tangent cone $C_x\X$ is the metric completion of $\widetilde C_x\X$, and its distance is denoted by $\dist_{C_x\X}$. If there are no non-constant geodesic issuing from $x$, we simply set $C_x\X:=\{0_x\}$.

For $v,w\in C_x\X$, we define
\[
|v|_x:=\dist_{C_x\X}(v,0_x),
\qquad
\langle v,w\rangle_x
:=
\frac12\left(
|v|_x^2+|w|_x^2-\dist_{C_x\X}(v,w)^2
\right).
\]
In particular, for all $s,t\geq0$ and $\gamma,\eta\in\Gamma_x$,
\[
|s[\gamma]|_x=s,
\qquad
\langle s[\gamma],t[\eta]\rangle_x
=st\cos\angle_x(\gamma,\eta).
\]
For every $\lambda\geq0$, the map $s[\gamma]\mapsto(\lambda s)[\gamma]$ extends uniquely and continuously to $C_x\X$. We denote its value at $v$
by $\lambda v$. It satisfies $\dist_{C_x\X}(\lambda v,\lambda w) =\lambda\dist_{C_x\X}(v,w)$.
\end{definition}

Even though $C_x\X$ is not a vector space in general, $|\cdot|_x$ and $\langle\cdot,\cdot\rangle_x$ play the roles of a norm and an inner product, respectively.


Fix $x\in\X$. For every $y\in\X\setminus\{x\}$ with $\dist(x,y)<D_\kappa$, let $\gamma\in\Gamma(x,y)$ be the unique geodesic and define
\[
\log_x y:=\dist(x,y)[\gamma],
\qquad \log_x x:=0_x.
\]
For $y\neq x$, the element $\log_x y\in C_x\X$ represents the direction of the geodesic from $x$ to $y$, with norm $\dd(x,y)$. For all $y,z\in\X\setminus\{x\}$ with $\dist(x,y),\dist(x,z)<D_\kappa$,
\[
\langle\log_x y,\log_x z\rangle_x
=
\dd(x,y)\dd(x,z)\cos\angle_x(y,z).
\]
If $y=x$ or $z=x$, the pseudo-inner product is zero by its definition.

\begin{proposition}[Geometry of tangent cones
{\cite[Theorem~II.3.19]{bridsonhaefliger1999}}]
\label{prop:tangent-cone-cat-zero}
Let $(\X,\dist)$ be a $\CAT(\kappa)$ space, with $\kappa\in\R$. Then, for every $x\in\X$, the tangent cone $C_x\X$ is a complete $\CAT(0)$ space.
\end{proposition}

\begin{proposition}[First variation formula {\cite[Corollary~II.3.6]{bridsonhaefliger1999}}]\label{prop:first-variation}
Let $(\X,\dd)$ be a $\CAT(\kappa)$ space, with $\kappa\in\R$, and let $D$ be a $\kappa$-domain. Let $x,y,z\in D$ with $x\neq y$ and $x\neq z$ and let $\gamma$ be the unique geodesic from $x$ to $z$. Then
\begin{equation}\label{eq:first-variation}
\lim_{t\downarrow0}
\frac{\dd(\gamma(t),y)-\dd(x,y)}{t}
=-\dd(x,z)\cos\angle_x(z,y)
=-\frac{\langle\log_x y,\log_x z\rangle_x}{\dd(x,y)}.
\end{equation}
\end{proposition}


\subsection{Projection onto geodesically convex sets}

We recall the projection property for geodesically convex sets, which will be used in Proposition~\ref{prop:bary_in_convex_hull}. A proof in the $\CAT(0)$ case is given in \cite[Proposition~II.2.4(1)--(3)]{bridsonhaefliger1999};
the positive-curvature case is addressed in \cite[Exercise~II.2.6(1)]{bridsonhaefliger1999}. In Appendix~\ref{subsec:proof_of_chebyshev}, we give an alternative proof for $\kappa>0$, using strong convexity of the squared distance on sufficiently small balls, as in \cite[Proposition~3.1(i)]{ohta2007}.
Our proof also applies when $\kappa\leq0$, using the strong convexity of the squared distance in $\CAT(0)$ spaces \cite[Proposition~2.3, inequality~(2.3)]{sturm2003}.

\begin{proposition}[Projection onto a geodesically convex set]
\label{lem:projection}
Let $(\X,\dist)$ be a $\CAT(\kappa)$ space, with $\kappa\in\R$. Let $D$ be a $\kappa$-domain of $\X$, and let $C\subseteq D$ be a nonempty, closed,
geodesically convex set. Let $x\in\X$ satisfy
\[
\dist(x,C):=\inf_{z\in C}\dist(x,z)<D_\kappa/2.
\]
Recall that $D_\kappa=+\infty$ when $\kappa\leq0$, so this restriction is then automatic. There exists a unique point $p_C(x)\in C$ such that
\[
\dist(x,p_C(x))=\dist(x,C).
\]
Furthermore, if $x\notin C$, then
\begin{equation}\label{eqn:obtuse}
\angle_{p_C(x)}(x,c)\geq\pi/2
\end{equation}
for every $c\in C\setminus\{p_C(x)\}$.
\end{proposition}

\section{Geometry of large \texorpdfstring{\boldmath$\CAT(\kappa)$}{CAT(kappa)} spaces}
\label{sec:new_geometric_results}

In this section, we state our
new four-point inequalities 
and geodesic stability results for $\CAT(\kappa)$ spaces, $\kappa>0$,
which will be essential to our proofs in
Sections~\ref{sec:geometric_properties_Frechet} 
and~\ref{sec:consequences}.

\subsection{A quantitative four-point inequality}

The $\CAT(\kappa)$ condition in Definition~\ref{def:CAT} above is expressed through triangle comparisons. For $\CAT(0)$ spaces, there
is an equivalent characterization of the $\CAT(0)$
condition via an inequality involving four points,
due to Reshetnyak~\cite{reshetnyak1968}; see also~\cite[Proposition~2.4]{sturm2003} and \cite[Theorem~1]{bergnikolaev2008}.

\begin{proposition}[Reshetnyak's four point inequality]
\label{prop:Reshetnyak}
A complete geodesic space $(\X,\dd)$ is a $\CAT(0)$ space if and only if,
for all $x_1,x_2,x_3,x_4\in\X$,
\[
\left|
\dd^2(x_1,x_3)-\dd^2(x_2,x_3)
-\dd^2(x_1,x_4)+\dd^2(x_2,x_4)
\right|
\leq
2\,\dd(x_1,x_2)\dd(x_3,x_4).
\]
\end{proposition}

In this section, we show how this inequality
can be extended to positively curved spaces
at the cost of a larger constant.
Previous work has considered analogs
of the four-point inequality in positively curved space~\cite{berg2015}, but under domain restrictions
that are too restrictive for the present work.
Our interest in four-point inequalities
is partial because they allow us to obtain,
in Proposition~\ref{prop:wasserstein-one-stability}, improved
contraction results for barycenters.

The main tool we will use is that of {\it Euclidean cones},
which we introduce here.

\begin{definition}[Euclidean cone over a metric space
{\cite[Definition~I.5.6]{bridsonhaefliger1999}}]\label{def:euclidean-cone}
Let $(\X,\dist)$ be a metric space. Consider the following equivalence relation on $[0,\infty)\times\X$: for all $s,t\geq0$ and $x,y\in\mathcal X$, $(t,x)\sim(s,y)$ if and only if
$t=s=0$, or $t=s>0$ and $x=y$. The quotient space is called the Euclidean cone over $(\X,\dist)$ and is denoted by
$C_0\mathcal X$.

For every $(t,x)\in[0,\infty)\times\X$, we write its equivalence class as $tx$. If $t=0$, we simply write it as $0$ and call it the vertex of $C_0\mathcal X$. Finally, for every $\kappa>0$, setting
\[
\dist_{C_0\X}^{(\kappa)}(tx,sy)
:=
\left(
t^2+s^2
-2st\cos\bigl(\min\{\pi,\sqrt\kappa\dist(x,y)\}\bigr)
\right)^{1/2}
\]
for all $x,y\in\X$ and $s,t\geq0$ defines a metric on $C_0\X$ by \cite[Proposition~I.5.9(1)]{bridsonhaefliger1999}.
For $\kappa=1$, we simply write $\dist_{C_0\X}^{(1)}=\dist_{C_0\X}$. Finally, for $A\subseteq\X$, we identify $C_0A$ with a subset of $C_0\X$.
\end{definition}

The benefit of considering Euclidean cones
is that, due to
Berestovskii's cone theorem \cite[Theorem~II.3.14]{bridsonhaefliger1999}, if $(\X,\dist)$ is a $\CAT(\kappa)$ space with $\kappa > 0$, then $(C_0\X,\dist_{C_0\X}^{(\kappa)})$ is a $\CAT(0)$ space.
We use this observation to transfer the 
four-point inequality from the Euclidean cone
to the underlying space, and obtain the following result.

\begin{proposition}[Quantitative four-point inequality]
\label{prop:local-four-point}
Let $(\X,\dist)$ be a $\CAT(\kappa)$ space, with $\kappa>0$.
Let $D$ be a $\kappa$-domain of $\X$.
Then, for every $x_1,x_2,x_3,x_4\in D$,
\begin{equation} \label{eq:squared-four-point}
\big|\dist^2(x_1,x_3)-\dist^2(x_2,x_3)-\dist^2(x_1,x_4)+\dist^2(x_2,x_4)\big|\leq c_{4\bullet}(r\sqrt\kappa)\,\dist(x_1,x_2)\dist(x_3,x_4),
\end{equation}
where $r\in(0,\pi/(2\sqrt{\kappa}))$ is the radius of $D$ and, for all $u\in [0,\pi/2)$,
\begin{equation}
\label{eq:Q-constant}
c_{4\bullet}(u):=
\frac{4u}{\sin(2u)}
+
2m(u)^2
\frac{
\sin(2u)
-2u\cos(2u)
}{
\sin^3(2u)
},
\end{equation}
with
\[
m(u)
:=
\max_{0\leq s\leq2u}\sin s
=
\begin{cases}
\sin(2u),
&0\leq u\leq\pi/4,\\
1,
&\pi/4<u<\pi/2.
\end{cases}
\]
The constant $c_{4\bullet}(u)$ is understood as $2$ when $u=0$. 
\end{proposition}

Let us make several remarks on this result. First,
one can obtain a simpler constant using
$c_{4\bullet}(u)\leq 64u^3/(3\sin^3(2u))$ for every $u\in(0,\pi/2)$.
Second, 
it is not possible to obtain a quantitative four-point inequality
on larger domains. This can be seen by taking $x_1, x_2, x_3 , x_4 \in \mathbb{S}^2$ on a great circle, such that $\dist(x_1, x_2), \dist(x_3, x_4)$
are small yet $x_1, x_3$ are antipodal and $x_2, x_4$ are antipodal.
Thirdly, for any $r>0$, $c_{4\bullet}(r\sqrt\kappa) \to 2$
as $\kappa \downarrow 0$, so we exactly recover
 Reshetnyak's quadruple comparison Proposition~\ref{prop:Reshetnyak}, in $\CAT(0)$ spaces.
 
Finally, we remark that, unlike Reshetnyak's quadruple comparison, Proposition~\ref{prop:local-four-point} unfortunately does not characterize
the $\CAT(\kappa)$ condition:
if geodesic space satisfies such an inequality, it does not necessarily
imply that the space has curvature bounded from above.
For example, as we verify in Appendix~\ref{subsec:l_4_example},
the unit ball in $(\R^2,\|\cdot\|_4)$ satisfies a four-point inequality with a finite constant, but has no finite upper curvature bound, even locally.

\begin{proof}[Proof of Proposition~\ref{prop:local-four-point}]
We work in the normalized case $\kappa=1$; the general case follows by rescaling the metric to $\sqrt{\kappa}\,\dd$.

Write $D=B(o,r)$ for some $o\in\X$ and $r\in (0,\pi/2)$.
Note that by the triangle inequality, for all $x,y\in D$, we have $\dist(x,y)\leq 2r<\pi$. 

Fix $x_1,x_2,x_3,x_4\in D$.

\medskip
\noindent\textbf{Step 1: A cosine four-point inequality.}
Consider the Euclidean cone $C_0\X$, equipped with its
canonical distance $d_{C_0\X}$.
By Berestovskii's cone theorem,
\cite[Theorem~II.3.14]{bridsonhaefliger1999},
this cone is $\CAT(0)$.
For $i=1,\ldots,4$, set $\bar x_i:=1x_i\in C_0\X$.
Since $\dist(x_i,x_j)<\pi$,
\[
d_{C_0\X}^2(\bar x_i,\bar x_j)
=2\bigl(1-\cos\dist(x_i,x_j)\bigr),
\qquad
d_{C_0\X}(\bar x_i,\bar x_j)
\leq\dist(x_i,x_j).
\]
Applying Proposition~\ref{prop:Reshetnyak} in $C_0\X$ yields
\[
\Big|d_{C_0\X}^2(\bar x_1,\bar x_3)-d_{C_0\X}^2(\bar x_2,\bar x_3)-d_{C_0\X}^2(\bar x_1,\bar x_4)+d_{C_0\X}^2(\bar x_2,\bar x_4)\Big| \leq 2d_{C_0\X}(\bar x_1,\bar x_2)d_{C_0\X}(\bar x_3,\bar x_4).
\]
Writing $a_{ij}:=1-\cos\dist(x_i,x_j)$ for $i,j=1,\ldots,4$ and dividing by $2$,
we obtain
\begin{equation}
\label{eq:appendix-cosine-four-point}
|a_{13}-a_{23}-a_{14}+a_{24}|
\leq\dist(x_1,x_2)\dist(x_3,x_4).
\end{equation}

\medskip
\noindent\textbf{Step 2: A scalar composition estimate.}
Let $H:I\to \R$ be a $C^2$ function on an non-empty interval $I$ and let
$L:=\max_I |H'|$ and $M:=\max_I|H''|$.
For $a,b\in I$, define
\[
J(a,b):=\int_0^1 H'((1-t)b+ta)\,\dd t
\]
so for all $a,b,c,d\in I$, $H(a)-H(b)=(a-b)J(a,b)$, $|J(a,b)|\leq L$, and
\[
|J(a,b)-J(c,d)|
\leq\frac M2\bigl(|a-c|+|b-d|\bigr).
\]
Consequently, using
\[
H(a)-H(b)-H(c)+H(d)
=(a-b-c+d)J(a,b)+(c-d)\bigl(J(a,b)-J(c,d)\bigr),
\]
we obtain
\begin{equation}
\label{eq:scalar-composition}
|H(a)-H(b)-H(c)+H(d)|\leq L|a-b-c+d|+\frac M2|c-d|
\bigl(|a-c|+|b-d|\bigr).
\end{equation}

\medskip
\noindent\textbf{Step 3: Transfer to squared distances.}
Take $I=[0,1-\cos(2r)]$ and 
\[
H(u):=\bigl(\arccos(1-u)\bigr)^2,
\qquad u\in[0,1-\cos(2r)].
\]
Thus $H(1-\cos s)=s^2$ for $s\in[0,2r]$.
For $s=\arccos(1-u)>0$,
\[
H'(u)=\frac{2s}{\sin s},
\qquad
H''(u)=2\frac{\sin s-s\cos s}{\sin^3s}.
\]
These derivatives extend continuously to $u=0$,
with $H'(0)=2$ and $H''(0)=2/3$.
Moreover,
\[
|H'|\leq L_r:=\frac{4r}{\sin(2r)},
\qquad
|H''|\leq M_r:=
2\frac{\sin(2r)-2r\cos(2r)}{\sin^3(2r)}.
\]
Indeed, $s/\sin s$ is increasing on $[0,\pi)$.
The function $f(s):=(\sin s-s\cos s)/\sin^3s$ is also
increasing, since
\[
f'(s)=\frac{g(s)}{\sin^4s},
\qquad
g(s):=s(1+2\cos^2s)-3\sin s\cos s,
\]
and
\[
g(0)=0,\qquad
g'(s)=4\sin s(\sin s-s\cos s)\geq0.
\]

Set $m:=m(r)=\max_{0\leq s\leq2r}\sin s$.
Since $s\mapsto1-\cos s$ is $m$-Lipschitz on $[0,2r]$,
the reverse triangle inequality gives
\[
\begin{aligned}
|a_{14}-a_{24}|&\leq m\dist(x_1,x_2),\\
|a_{13}-a_{14}|&\leq m\dist(x_3,x_4),\\
|a_{23}-a_{24}|&\leq m\dist(x_3,x_4).
\end{aligned}
\]
Applying \eqref{eq:scalar-composition} to
$(a,b,c,d)=(a_{13},a_{23},a_{14},a_{24})$,
and using \eqref{eq:appendix-cosine-four-point}, yields
\[
\begin{aligned}
\Big|
&\dist^2(x_1,x_3)-\dist^2(x_2,x_3)
-\dist^2(x_1,x_4)+\dist^2(x_2,x_4)
\Big|\\
&\leq
\bigl(L_r+M_rm^2\bigr)
\dist(x_1,x_2)\dist(x_3,x_4),
\end{aligned}
\]
proving the normalized case since $L_r+M_rm^2=c_{4\bullet}(r)$.
\end{proof}

\subsection{Stability of geodesics with common initial or terminal point}
Our proof of high probability bounds on iterated barycenters
is based on showing that iterated barycenters
are Lipschitz in the sequence of inputted points.
To show this Lipschitz statement, we require control
on how much geodesics in $\kappa$-domains vary when their initial and
terminal points are changed.

When the initial point varies, stability was studied by Ohta~\cite{ohta2007}.
\begin{proposition}[Stability with a common initial point
{\cite[Proposition~3.1(ii)]{ohta2007}}]
\label{eq:L-convexity}
Let $(\X,\dist)$ be a $\CAT(\kappa)$ space, with $\kappa>0$, and let $D$ be a $\kappa$-domain.
For $x,y,z\in D$, let $\gamma\in\Gamma(x,y)$ and $\eta\in\Gamma(x,z)$.
Then, for every $t\in[0,1]$,
\[
\dist(\gamma(t),\eta(t))
\leq
\frac{2r\sqrt\kappa}
     {\sin(2r\sqrt\kappa)}
\,t\,\dist(y,z)
\]
where $r\in[0,\pi/(2\sqrt{\kappa}))$ is the radius of $D$. The quotient is understood as $1$ when $r=0$.
\end{proposition}

We also require stability with respect
to the terminal point. While one could
use Proposition~\ref{eq:L-convexity} with
$1 - t$ in place of $t$, this would not give a strong
enough contraction for our purposes (namely, for Theorem~\ref{thm:stability of it bar} below), where we will rather need a bound with a factor of the form $(1-\alpha t)$ for some $\alpha>0$. To overcome this issue, we make make the following
crucial definition 
a weighted version of the distance $\dist$ on $\kappa$-domains. 

\begin{definition}[Weighted distance]
\label{def:weighted-distance}
Let $(\X,\dist)$ be a $\CAT(\kappa)$ space for some $\kappa>0$, and let $D=B(o,r)$ be a $\kappa$-domain, with
$0\leq r<\pi/(2\sqrt{\kappa})$.
For $x,y\in D$, define
\[
\dist_o(x,y):=
\inf_{\substack{\gamma\in A(x,y)\\\gamma([0,1])\subseteq D}}
\int_0^1
\frac{|\dot\gamma|(s)}
     {\cos\bigl(\sqrt{\kappa}\,\dist(o,\gamma(s))\bigr)}
\,\dd s,
\]
where $A(x,y)$ is defined in
Definition~\ref{def:metric-speed}, and
$|\dot\gamma|$ denotes the metric speed from
Definition~\ref{def:metric-speed}.
Then using \cite[Theorem~1.1.2]{ambrosio2005gradient}, $\dist_o$ defines a distance on $D$ satisfying
\begin{equation}
\label{eq:equivalence-inequality}
\dist(x,y)
\leq \dist_o(x,y)
\leq \frac{1}{\cos(r\sqrt\kappa)}\dist(x,y)
\end{equation}
for all $x,y\in D$.
\end{definition}

The fact that $\dd_o$ defines a distance on $D$ follows by reversing and concatenating curves, together with the bounds below. Moreover, one easily checks that for all $x,y\in D$,
\[
\dist(x,y)=\inf_{\substack{\gamma\in A(x,y)\\\gamma([0,1])\subseteq D}}
\int_0^1|\Dot\gamma|(s)\,\dd s.
\]
Indeed, the left-hand side is bounded from above by the right-hand side by Definition~\ref{def:metric-speed}, and equality is attained by choosing $\gamma\in\Gamma(x,y)\subseteq A(x,y)$, which stays in $D$ by Proposition~\ref{prop:unique-geodesics-domain} and has constant metric speed $\dist(x,y)$. The inequalities~\eqref{eq:equivalence-inequality} then follow from the bounds on the weight and ensure that $\dd_o$ is finite and vanishes only when $x=y$.


\begin{proposition}[Stability with common terminal point]
\label{prop:endpoint-inequality}
Let $(\X,\dist)$ be a $\CAT(\kappa)$ space with $\kappa>0$ and let $D=B(o,r)$ be a $\kappa$-domain, for some $o\in \X$ and $r\in [0,\pi/(2\sqrt\kappa))$. Let $x,y,z\in D$ and $\gamma\in \Gamma(x,z)$, $\eta\in\Gamma(y,z)$.
Then, for every $t\in [0,1]$,
\[
\dist_o(\gamma(t),\eta(t))\leq (1-\cos(r\sqrt\kappa)t) \dist_o(x,y)
\]
\end{proposition}

In order to prove this proposition, we use the following infinitesimal estimate,
whose proof appears in Appendix~\ref{subsec:endpoint_inequality}.

\begin{lemma} \label{lem:stability_same_endpoint}
With the notation of Proposition~\ref{prop:endpoint-inequality} above, we have, for all $t\in[0,1]$,
\[
\limsup_{\substack{y\to x \\ y\in D\setminus\{x\}}}
\frac{d(\gamma(t),\eta(t))}{d(x,y)}
\leq
\frac{\sin((1-t)d(x,z)\sqrt\kappa)}{\sin(d(x,z)\sqrt\kappa)}
\]
Here $\eta\in\Gamma(y,z)$ is the unique geodesic corresponding to each $y$, and the right-hand side is understood as $1-t$ if $x=z$.
\end{lemma}

\begin{figure}[htbp]
  \centering
  \includegraphics[width=0.50\textwidth]
    {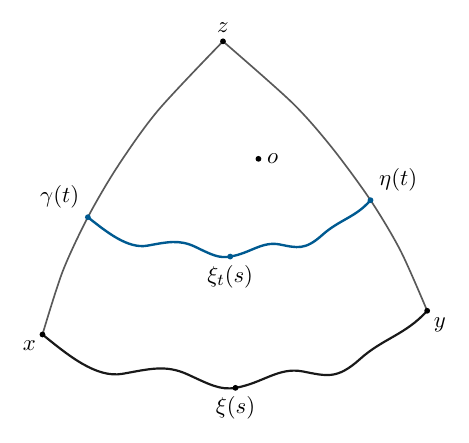}
  \caption{Setup in the proof of Proposition~\ref{prop:endpoint-inequality}.
    The curve $\xi$ joining $x$ to $y$ is mapped to
    $\xi_t(s)=\lambda_s(t)$, where $\lambda_s$ is the geodesic
    from $\xi(s)$ to $z$.}
  \label{fig:common-endpoint-stability}
\end{figure}

\begin{proof}[Proof of Proposition~\ref{prop:endpoint-inequality}]
We work in the normalized case $\kappa=1$; the general case follows by rescaling the metric to $\sqrt{\kappa}\,\dd$. The cases $t=0$ and $t=1$ are immediate, so fix $t\in(0,1)$.
Let $\xi:[0,1]\to B(o,r)$ be an absolutely continuous path joining $x$ and $y$. For each $s\in[0,1]$, let $\lambda_s\in\Gamma(\xi(s),z)$ and set $\xi_t(s):=\lambda_s(t)$: This defines a new path $\xi_t:[0,1]\to\X$. Since $B(o,r)$ is geodesically convex, $\xi_t(s)\in B(o,r)$ for all $s\in[0,1]$. For an illustration of this construction, see Figure~\ref{fig:common-endpoint-stability}.

If, for $p,q\in D$, we denote by $\gamma_{p,q}$ the unique geodesic from $p$ to $q$, then Lemma~\ref{lem:stability_same_endpoint} gives an infinitesimal bound for the radial map $p\in D\mapsto\gamma_{p,z}(t)$, uniformly bounded above by $2r(1-t)/\sin(2r)$. By the mean value inequality along geodesics in $D$, this map is Lipschitz. Its composition with $\xi$ is therefore absolutely continuous, showing that $\xi_t\in A(\gamma(t),\eta(t))$.

Therefore, by \cite[Theorem~1.1.2]{ambrosio2005gradient}, both ratios $\frac{\dd(\xi(s'),\xi(s))} {|s'-s|}$ and $\frac{\dd(\xi_t(s'),\xi_t(s))} {|s'-s|}$ have limits as $s'\to s$ for almost all $s\in [0,1]$, with respect to the Lebesgue measure -- these limits are $|\Dot\xi|(s)$ and $|\Dot\xi_t|(s)$ respectively, as per Definition~\ref{def:metric-speed}. Let $s\in [0,1]$ be such a number for which both limits exist. If $|\Dot\xi|(s)\neq 0$, then Lemma~\ref{lem:stability_same_endpoint} yields that 
\begin{equation*}
    \frac{|\Dot\xi_t|(s)}{|\Dot\xi|(s)} = \lim_{\substack{s'\to s \\ s'\ne s}} \frac{\dist(\xi_t(s'),\xi_t(s))}{\dist(\xi(s'),\xi(s))} \leq \frac{\sin((1-t)\dist(\xi(s),z))}{\sin\dist(\xi(s),z)}
\end{equation*}
where the right-hand side is understood as $1-t$ if $\xi(s)=z$. If $|\Dot\xi|(s)=0$, let $(s_n)_{n\geq 1}$ be any sequence of numbers in $[0,1]\setminus\{s\}$ with $s_n\to s$ as $n\to\infty$. If $\xi(s_n)=\xi(s)$ for all large enough $n$, then it also holds that $\xi_t(s_n)=\xi_t(s)$ for all large enough $n$, so
$$|\Dot\xi_t|(s)=\lim_{n\to\infty}\frac{\dist(\xi_t(s_n),\xi_t(s))}{|s_n-s|}=0=|\Dot\xi|(s).$$
Otherwise, consider a subsequence $(s_{k_n})_{n\geq 1}$, with $k_n\to\infty$ as $n\to\infty$, such that $\xi(s_{k_n})\neq\xi(s)$ for all $n\geq 1$. Then, again by Lemma~\ref{lem:stability_same_endpoint},
$$\infty>\limsup_{n\to\infty} \frac{\dist(\xi_t(s_{k_n}),\xi_t(s))}{\dist(\xi(s_{k_n}),\xi(s))}=\limsup_{n\to\infty} \frac{\dist(\xi_t(s_{k_n}),\xi_t(s))}{|s_{k_n}-s|} \left(\frac{\dist(\xi(s_{k_n}),\xi(s))}{|s_{k_n}-s|}\right)^{-1}$$
so necessarily, 
$$|\Dot\xi_t|(s)=\lim_{n\to\infty} \frac{\dist(\xi_t(s_{k_n}),\xi_t(s))}{|s_{k_n}-s|}=0.$$
In every case, we have established that $$|\Dot\xi_t|(s)\leq |\Dot\xi|(s)\frac{\sin((1-t)\dist(\xi(s),z))}{\sin\dist(\xi(s),z)}.$$

Dividing by $\cos\dd(o,\xi_t(s))>0$, we obtain
\begin{equation}\label{eqn:proof_31415}
\begin{aligned}
\frac{|\Dot\xi_t|(s)}{\cos\dd(o,\xi_t(s))}
&\leq
\frac{|\Dot\xi|(s)}{\cos\dd(o,\xi(s))}
\frac{\cos\dd(o,\xi(s))\sin((1-t)\dd(\xi(s),z))}
     {\cos\dd(o,\xi_t(s))\sin\dd(\xi(s),z)}.
\end{aligned}
\end{equation}

Our goal now is to show that the second ratio on the right-hand side is bounded by $1-t\cos r$,
independently of $s$.

Suppose first that $z\neq o$ and $\xi(s)\neq z$. By the triangle inequality,
\[
\dd(z,\xi(s))+\dd(z,o)+\dd(o,\xi(s))
\leq4r<2\pi.
\]
Hence, by
\cite[Proposition~II.1.7(3)]{bridsonhaefliger1999},
\[
\angle_z^{(1)}(o,\xi_t(s))
\leq\angle_z^{(1)}(o,\xi(s)).
\]
Since $\cos$ is decreasing on $[0,\pi]$,
\[
\cos\angle_z^{(1)}(o,\xi_t(s))
\geq\cos\angle_z^{(1)}(o,\xi(s)).
\]
Using Definition~\ref{def:comparison_angles}, we obtain
\[
\frac{
\cos\dd(o,\xi_t(s))
-\cos\dd(o,z)\cos\dd(\xi_t(s),z)}
{\sin\dd(o,z)\sin\dd(\xi_t(s),z)}\geq
\frac{
\cos\dd(o,\xi(s))
-\cos\dd(o,z)\cos\dd(\xi(s),z)}
{\sin\dd(o,z)\sin\dd(\xi(s),z)}.
\]

Multiplying by the positive denominator on the left gives
\[
\begin{aligned}
&\cos\dd(o,\xi_t(s))
-\cos\dd(o,z)\cos\dd(\xi_t(s),z)\\
&\qquad\geq
\left(
\cos\dd(o,\xi(s))
-\cos\dd(o,z)\cos\dd(\xi(s),z)
\right)
\frac{\sin\dd(\xi_t(s),z)}
     {\sin\dd(\xi(s),z)}.
\end{aligned}
\]
Note that if $z=o$ (but still, $\xi(s)\neq z$), this last inequality holds trivially, and it is in fact an equality. Rearranging and setting $a:=\dd(\xi(s),z)\in(0,\pi)$, we obtain
\begin{align*}
&\frac{\cos\dd(o,\xi(s))\sin((1-t)\dd(\xi(s),z))}
     {\cos\dd(o,\xi_t(s))\sin\dd(\xi(s),z)}\\
&\quad\leq
1-\frac{\cos\dd(o,z)}{\cos\dd(o,\xi_t(s))}
\left(
\cos\dd(\xi_t(s),z)
-\frac{\cos\dd(\xi(s),z)\sin\dd(\xi_t(s),z)}
      {\sin\dd(\xi(s),z)}
\right)\\
&\quad=
1-\frac{\cos\dd(o,z)}{\cos\dd(o,\xi_t(s))}
\frac{\cos((1-t)a)\sin a-\cos a\sin((1-t)a)}
     {\sin a}\\
&\quad=
1-\frac{\cos\dd(o,z)}{\cos\dd(o,\xi_t(s))}
\frac{\sin(ta)}{\sin a}.
\end{align*}

The desired bound follows since $\cos\dd(o,z)\geq\cos r$ and $\cos\dd(o,\xi_t(s))\leq1$. Moreover, concavity of $\sin$ on $[0,\pi]$ gives $\sin(ta)/\sin a\geq t$.

Now, when $\xi(s)=z$, the second ratio on the right-hand side of \eqref{eqn:proof_31415} is understood as $1-t$, which is also at most $1-t\cos r$.

Finally, we have established that for almost every $s\in[0,1]$, \begin{equation*}
\frac{|\Dot\xi_t|(s)}{\cos\dd(o,\xi_t(s))}
\leq
(1-t\cos r)
\frac{|\Dot\xi|(s)}{\cos\dd(o,\xi(s))}.
\end{equation*}
Integrating over $s\in[0,1]$ and recalling that $\xi_t\in A(\gamma(t),\eta(t))$, we obtain,
by Definition~\ref{def:weighted-distance},
\[
\begin{aligned}
\dd_o(\gamma(t),\eta(t))
&\leq
\int_0^1
\frac{|\Dot\xi_t|(s)}{\cos\dd(o,\xi_t(s))}\,\dd s\\
&\leq
(1-t\cos r)\int_0^1
\frac{|\Dot\xi|(s)}{\cos\dd(o,\xi(s))}\,\dd s.
\end{aligned}
\]
Taking the infimum over all $\xi\in A(x,y)$ with values in $B(o,r)$ yields the desired result.
\end{proof}

\section{Fréchet functions on large \texorpdfstring{\boldmath$\CAT(\kappa)$}{CAT(kappa)} spaces}
\label{sec:geometric_properties_Frechet}

Given $\mu\in\mathcal P_2(\X)$, we define its Fréchet function $$F_\mu:x\in \X \mapsto \frac{1}{2}\int_\X \dist(x,y)^2\diff\mu(y).$$
Any minimizer of $F_\mu$ is called a barycenter of $\mu$\footnote{One may also define barycenters for any $\mu\in \mathcal P_1(\X)$ by letting $F_\mu(x)=\frac{1}{2}\int_\X (\dist(x,y)^2-\dist(x_0,y)^2)\diff\mu(y)$ for some arbitrary $x_0\in \X$. However, for simplicity, we will only consider $\mu\in\mathcal P_2(\X)$ in this work.}.
Barycenters are the canonical generalization of average
to Euclidean spaces, since if $\X$ is Euclidean then $\mu\in\mathcal P_2(\X)$ has a unique barycenter given by the mean 
\(\int_{\X} y\,\dd\mu(y)\).

The study of barycenters goes back to Fréchet~\cite{frechet1948} and later work of Karcher~\cite{karcher1977} and barycenters have become standard tools in shape statistics \cite{kendall2009}, directional and spherical data analysis \cite{bussfillmore2001},
and optimal transport \cite{agueh2011,cuturi2014}.


In the rest of this section, and unless stated otherwise, we assume that $(\X,\dist)$ is a $\CAT(\kappa)$ space for a fixed $\kappa>0$, and that $\mu$ is a Borel probability measure supported on a $\kappa$-domain $D=B(o,r)$, with $o\in \X$ and $0<r<D_\kappa/2$.The case $r=0$ is immediate, since $\mu$ is then a Dirac measure. In particular, $\mu\in\mathcal P_2(\X)$.

\subsection{Bounds on differences
of a Fréchet function}



We shall require the following definition of
the directional derivative of the
Fréchet functional, whose proof is included in Appendix~\ref{subsec:proof_dir_derivative}.
\begin{proposition}[Directional derivative of the Fr\'echet function]
\label{prop:dir-derivative}
Let $(\X,\dist)$ be a $\CAT(\kappa)$ space with $\kappa\in\R$, let $D$ be a $\kappa$-domain of $\X$, and let $\mu\in\mathcal P_2(\X)$ be concentrated on $D$.
For $x,z\in D$, let $\gamma_{xz}\in\Gamma(x,z)$, which is unique by Proposition~\ref{prop:unique-geodesics-domain}.
We define
\[
DF_\mu(x)[\log_x z]
:=
\lim_{t\downarrow0}
\frac{F_\mu(\gamma_{xz}(t))-F_\mu(x)}{t}.
\]
Then this quantity exists and is given by
\begin{equation}\label{eq:dir-derivative}
DF_\mu(x)[\log_x z]
=
-\int_D
\big\langle\log_x y,\log_x z\big\rangle_x
\,\dd\mu(y).
\end{equation}
Moreover, defining the directional derivative along all
geodesics issuing from $x$, with arbitrary nonnegative speeds,
yields a function on $\widetilde C_x\X$ that extends uniquely
to a Lipschitz function on $C_x\X$, still denoted
\[
DF_\mu(x)[v]
=-\int_D\langle\log_x y,v\rangle_x\,\dd\mu(y),
\qquad v\in C_x\X.
\]
Its Lipschitz constant is at most
$\int_D\dist(x,y)\,\dd\mu(y)<\infty$.
In particular, if $D=B(o,r)$ with $r<\infty$,
this extension is $2r$-Lipschitz.
\end{proposition}

To study the behaviour of the Fréchet function,
we must get control on
the differences $F_\mu(x) - F_\mu(z)$
for points $x, z \in D$ in the $\kappa$-domain.
The following inequality is our first
step in this direction,
and relates differences of squared distances to differences of cosines, which are more
readily suited to 
comparison inequalities we get
from the 
 $\CAT(\kappa)$ condition.
 
\begin{lemma}[Fr\'echet-cosine inequality]
\label{lem:cosine-frechet-transfer}
For every $x,z\in D$,
\begin{equation}
\label{eq:integrated-transfer-upper}
F_\mu(x)-F_\mu(z)
\leq
-\frac1{\kappa}\int_D
\frac{\sqrt{\kappa}\,\dd(x,y)}
     {\sin\bigl(\sqrt{\kappa}\,\dd(x,y)\bigr)}
\Bigl[
\cos\bigl(\sqrt{\kappa}\,\dd(x,y)\bigr)
-
\cos\bigl(\sqrt{\kappa}\,\dd(z,y)\bigr)
\Bigr]
\,\dd\mu(y).
\end{equation}
The quotient $u/\sin u$ is understood as $1$ at $u=0$.
\end{lemma}

\begin{proof}
We work in the normalized case $\kappa=1$; the general case follows by rescaling the metric to $\sqrt{\kappa}\,\dd$. Consider
\[
    H\colon(-1,1]\longrightarrow[0,\pi^2),
    \qquad
    H(u):=\arccos^2(u).
\]
If \(u=\cos\theta\), with \(\theta\in(0,\pi)\), then
\[
    H'(u)=-\frac{2\theta}{\sin\theta},
    \qquad
    H''(u)
    =
    2\frac{\sin\theta-\theta\cos\theta}{\sin^3\theta}\geq0.
\]
Thus \(H\) is convex. 
Fix $y\in D$. Since $\dd(x,y),\dd(z,y)\leq2r<\pi$, the first-order convexity inequality for $H$ gives
\[
\begin{aligned}
H(\cos\dd(z,y))-H(\cos\dd(x,y))
&\geq H'(\cos\dd(x,y))\bigl(\cos\dd(z,y)-\cos\dd(x,y)\bigr).
\end{aligned}
\]
Using $H(\cos u)=u^2$ for $u\in[0,\pi)$, we obtain
\[
\begin{aligned}
\dd^2(x,y)-\dd^2(z,y)
&\leq
-2\frac{\dd(x,y)}{\sin\dd(x,y)}
\bigl(\cos\dd(x,y)-\cos\dd(z,y)\bigr).
\end{aligned}
\]
Integration with respect to \(\mu\), followed by division by \(2\), proves
\eqref{eq:integrated-transfer-upper}.
\end{proof}

Using this upper bound, we can obtain
the following useful bound on the
difference $F_\mu(x) - F_\mu(z)$ in terms
of the directional derivatives of $F_\mu$ at
$x$.

The following theorem provides a fundamental inequality that will be used throughout the rest of the paper. It is inspired by Afsari's analysis in the Riemannian
setting~\cite{afsari2011}.

\begin{theorem}[Directional Fr\'echet-cosine inequality]
\label{thm:directional-derivative-transfer}
For every $x,z\in D$,
\begin{equation}
\label{eq:directional-transfer-hemispherical}
\begin{aligned}
F_\mu(x)-F_\mu(z)
&\leq
-\frac{\sin\bigl(\sqrt{\kappa}\,\dd(x,z)\bigr)}
      {\sqrt{\kappa}\,\dd(x,z)}
    DF_\mu(x)[\log_x z]\\
&\quad
-\frac{1-\cos\bigl(\sqrt{\kappa}\,\dd(x,z)\bigr)}
      {\kappa\cos\bigl(\sqrt{\kappa}\,\dd(x,o)\bigr)}
    \left[
        \cos(r\sqrt\kappa)
        +
        \kappa
        \frac{\sin\bigl(\sqrt{\kappa}\,\dd(x,o)\bigr)}
             {\sqrt{\kappa}\,\dd(x,o)}
        DF_\mu(x)[\log_x o]
    \right].
\end{aligned}
\end{equation}
The quotient $\sin u/u$ is understood as $1$ at $u=0$.
\end{theorem}

\begin{proof}
We work in the normalized case $\kappa=1$; the general case follows by rescaling the metric to $\sqrt{\kappa}\,\dd$. The case $x=z$ is immediate, so assume $x\neq z$.
For every $y\in D\setminus\{x\}$, all distances between
$x,z,y$ are strictly smaller than $\pi$. By Proposition~\ref{prop:cosine-comparison-domain},
\begin{equation}
\label{eq:extended-angle-comparison-transfer}
\frac{
    \cos\dd(z,y)-\cos\dd(x,z)\cos\dd(x,y)
}{
    \sin\dd(x,z)\sin\dd(x,y)
}
\leq
\cos\angle_x(z,y).
\end{equation}
Multiplying \eqref{eq:extended-angle-comparison-transfer} by
$\dd(x,y)\sin\dd(x,z)$, we obtain
\[
\frac{\dd(x,y)}{\sin\dd(x,y)}
\Bigl(
    \cos\dd(z,y)
    -
    \cos\dd(x,z)\cos\dd(x,y)
\Bigr)
\leq
\frac{\sin\dd(x,z)}{\dd(x,z)}
\left\langle\log_x z,\log_x y\right\rangle_x.
\]
This inequality also holds when $y=x$, with both sides equal to zero, using the continuous extension of $u/\sin u$.

Together with
\[
\begin{aligned}
\cos\dd(z,y)-\cos\dd(x,y)
&=
\Bigl(\cos\dd(z,y)-\cos\dd(x,z)\cos\dd(x,y)\Bigr)\\
&\quad-
\bigl(1-\cos\dd(x,z)\bigr)\cos\dd(x,y),
\end{aligned}
\]
and recalling from Proposition~\ref{prop:dir-derivative} the directional derivative of $F_\mu$, Lemma~\ref{lem:cosine-frechet-transfer} gives
\[
    F_\mu(x)-F_\mu(z)
    \leq
-\frac{\sin d(x,z)}{d(x,z)}
    D F_\mu(x)[\log_xz]
    -
    \bigl(1-\cos d(x,z)\bigr)J_\mu(x),
\]
where
\begin{equation}\label{eq:def-Jmu}
J_\mu(x):=
\int_\X
\dist(x,y)\cot\bigl(\dist(x,y)\bigr)\,\diff\mu(y).
\end{equation}
Here $u\cot u$ is understood as $1$ at $u=0$.

We need to lower bound $J_\mu(x)$. Naively, one could use the boundedness of the support:
since $x,y\in D=B(o,r)$, we have $\dd(x,y)\leq2r<\pi$, and the monotonicity of $u\mapsto u\cot u$ gives $J_\mu(x)\geq2r\cot(2r)$. However, this lower bound is nonpositive when $r\geq\pi/4$. We therefore need a bound better suited to our setting. To address this issue, which will also arise later, the idea is to relate $x$ and $y$ in the integrand to
the centre $o$ of the ball. This is similar to the approach in \cite{afsari2011}.

If $x=o$, then $J_\mu(o)\geq\cos r$ and
$DF_\mu(o)[\log_o o]=0$, so the preceding bound proves
the claim. We may therefore assume $\dd(x,o)>0$.

By Proposition~\ref{prop:cosine-comparison-domain},
\[
\frac{\dd(x,y)}{\sin\dd(x,y)}
\Bigl[
\cos\dd(o,y)-\cos\dd(x,o)\cos\dd(x,y)
\Bigr]
\leq
\frac{\sin\dd(x,o)}{\dd(x,o)}
\langle\log_x o,\log_x y\rangle_x.
\]
The inequality also holds when $y=x$, by continuous extension. Integrating and using Proposition~\ref{prop:dir-derivative}, we obtain
\begin{equation}\label{eq:afsari-inequality}
\begin{aligned}
\cos\dd(x,o)\,J_\mu(x)
&\geq
\int_D\frac{\dd(x,y)}{\sin\dd(x,y)}
\cos\dd(o,y)\,\dd\mu(y)\\
&\quad+
\frac{\sin\dd(x,o)}{\dd(x,o)}
DF_\mu(x)[\log_x o].
\end{aligned}
\end{equation}
Since $\dd(x,y)/\sin\dd(x,y)\geq1$ and $\cos\dd(o,y)\geq\cos r>0$, substitution into the
preceding bound proves the desired inequality.
\end{proof}

\subsection{Existence and uniqueness of barycenters on \texorpdfstring{\boldmath$\kappa$}{kappa}-domains}

Theorem~\ref{thm:directional-derivative-transfer} will be our main tool for studying the Fréchet function in the following sections. We first use it to recover Yokota's existence and uniqueness result for probability measures supported on a $\kappa$-domain $D=B(o,r)$ of radius $r<D_\kappa/2$.
In \cite[Theorem~1]{bussfillmore2001}, for a discrete measure associated with $n$ points on the upper hemisphere of the unit sphere $\Sp^{d-1}$, with a positive mass in its interior, the authors proved existence and uniqueness. Compactness of the sphere guarantees existence and they show that every minimizer lies in the interior of this hemisphere. Their uniqueness argument uses the positive definiteness of the Hessian at critical points. In the Riemannian setting, related uniqueness results were established by Kendall~\cite[Theorem~7.3]{kendall1990} using auxiliary geodesically convex functions, later referred to as Kendall's functions by Yokota~\cite{yokota2016}. Afsari~\cite[Theorem~2.1]{afsari2011} proved the existence and uniqueness of global $L^p$ centers of mass, including the quadratic case $p=2$, under support-radius conditions involving the injectivity radius and an upper sectional curvature bound. The proof uses comparison and Hessian estimates, following an approach similar to that of \cite{bussfillmore2001}.

In our metric setting, we cannot rely on the same compactness argument, since closed bounded balls in $\X$ need not be compact. Indeed, $\X$ is not assumed to be locally compact, so the Hopf--Rinow--Cohn-Vossen theorem \cite[Theorem~2.5.28]{burago2001course} does not apply. Consider, for instance, the closed unit ball in an infinite-dimensional Hilbert space. Yokota~\cite[Theorem~B]{yokota2016} overcame this difficulty using Ekeland's variational principle~\cite{ekeland1979}. We recover his result using the same variational principle.

\begin{theorem}[Existence and uniqueness {(see \cite[Theorem~B]{yokota2016})}]
\label{thm:existence_uniqueness}
The Fréchet function $F_\mu$ admits a unique minimizer $b(\mu)$ over $\X$, and $b(\mu)\in D$.
\end{theorem}

The proof is given in Appendix~\ref{app:existence-uniqueness}; we outline its main steps here. Ekeland's variational principle provides a family $(x_\varepsilon)_{\varepsilon>0}$ whose Fréchet values converge to $\inf_\X F_\mu$ and which satisfies approximate first-order conditions.

Yokota's first localization lemma gives a uniform gap away from the support ball: for every $\alpha>0$, there exists $a_\alpha>0$ such that
\[
F_\mu(x)>\inf_{z\in B(o,r)}F_\mu(z)+a_\alpha
\qquad\text{whenever }x\notin B(o,r+\alpha).
\]
Consequently, no minimizer can lie outside $B(o,r)$: for any such point $x$, choose $0<\alpha<\dd(o,x)-r$ and apply the inequality. Since $\inf_\X F_\mu\leq\inf_{B(o,r)}F_\mu$, the same gap also forces the approximate minimizers eventually to belong to every ball $B(o,\rho)$ with $\rho>r$.

\begin{figure}[htbp]
    \centering
    \includegraphics[width=0.5\textwidth]{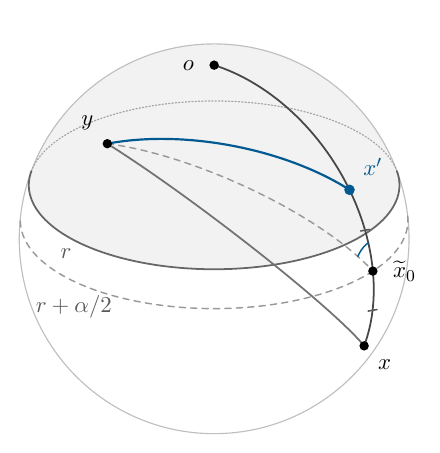}
    \caption{Radial reflection in the spherical comparison configuration. Matching marks indicate $\dd(x',\widetilde{x}_0)=\dd(\widetilde{x}_0,x)$. The blue arc marks the acute inward angle at $\widetilde{x}_0$, and the blue geodesic joining $y$ to $x'$ is shorter than the gray geodesic joining $y$ to $x$.}
    \label{fig:yokota-localization}
\end{figure}

The geometric idea behind this first localization step is illustrated in Figure~\ref{fig:yokota-localization}. In the normalized case $\kappa=1$, take $\alpha>0$ sufficiently small and consider $x$ such that $r+\alpha\leq\dd(x,o)<2r$. Reflect $x$ along the radial geodesic through $o$ about the point $\widetilde{x}_0$ at distance $r+\alpha/2$ from $o$. The resulting point $x'$ satisfies
\[
\dd(o,x')=2r+\alpha-\dd(o,x)\leq r,\qquad \dd(x',\widetilde{x}_0)=\dd(\widetilde{x}_0,x).
\]
In the spherical comparison configuration shown in the figure, the acute inward angle at $\widetilde{x}_0$ makes the blue geodesic from $y$ to $x'$ shorter than the gray geodesic from $y$ to $x$. Integrating such distance comparisons against $\mu$ yields a decrease of the Fréchet function. Together with estimates covering the remaining configurations, this gives the localization argument.

Yokota's second localization lemma provides an additional differential description: for every $\rho\in(r,D_\kappa/2)$, there exists $a_\rho>0$ such that
\[
DF_\mu(x)\left[\frac{\log_x o}{\dd(x,o)}\right]<-a_\rho
\qquad\text{for }x\in B(o,D_\kappa/2)\setminus B(o,\rho).
\]
Thus moving towards $o$ decreases $F_\mu$ at a uniformly positive rate in this region. Ekeland's approximate first-order condition bounds the same unit-direction derivative at $x_\varepsilon$ below by $-\varepsilon$. Once $x_\varepsilon$ lies in $B(o,D_\kappa/2)$, these estimates are incompatible for $\varepsilon<a_\rho$ unless $x_\varepsilon\in B(o,\rho)$.

Once the approximate minimizers are localized in $B(o,\rho)$ with $r<\rho<D_\kappa/2$, Theorem~\ref{thm:directional-derivative-transfer}, combined with their approximate first-order conditions, directly controls the distance between two such points and shows that the family is Cauchy. This replaces the use of Kendall's function in the proof of \cite[Theorem~B]{yokota2016}. Completeness and continuity then give a minimizer in $B(o,r)$, and the same comparison applied to two minimizers proves uniqueness.

By Theorem~\ref{thm:existence_uniqueness}, $b(\mu)$ is uniquely defined and belongs to $B(o,r)$. The following result provides a sharper localization.

\begin{proposition}[{\cite[Corollary~16]{yokota2016}}]
\label{prop:bary_in_convex_hull}
Let $C\subseteq D$ be a closed geodesically convex set such that $\mu(C)=1$. Denote by $b:=b(\mu)\in D$ the unique barycenter of $\mu$. Then $b\in C$.
\end{proposition}

If $\X$ is separable, the support of $\mu$ has full $\mu$-measure \cite[Chapter II, Theorem 2.1]{parthasarathy1967probability} so we may take $C=\overline{\conv(\supp(\mu))}$.

A proof of Proposition~\ref{prop:bary_in_convex_hull} is given in \cite[Proposition~3]{gietl2024}, under the stronger assumption that $r<\pi/(4\sqrt\kappa)$. Under our assumptions, \cite[Corollary~16]{yokota2016} gave a proof that is purely analytical. Here, we give a geometric approach, extending the proof of \cite[Proposition~3]{gietl2024} and relying on Proposition~\ref{lem:projection} above.

\begin{proof}
We work in the normalized case $\kappa=1$; the general case follows by rescaling the metric to $\sqrt{\kappa}\,\dd$. For brevity, write $b:=b(\mu)$.
Since $\mu(C)=1$,
\begin{equation*}
    2F_\mu(b) = \int_{\X}\dist^2(b,x)\diff\mu(x) = \int_{C}\dist^2(b,x)\diff\mu(x) \geq \dist^2(b,C)\mu(C) = \dist^2(b,C),
\end{equation*}
so 
$$\dd(b,C)^2\leq 2F_\mu(b)\leq 2F_\mu(o)\leq r^2,$$
hence $\dd(b,C)\leq r<\pi/2$. By contradiction, assume that $b\in B(o,r)\setminus C$.
By Proposition~\ref{lem:projection},
we can consider its projection onto $C$, denoted by
$p:=p_C(b)$.
Since $C\subseteq B(o,r)$, all distances between $b,p,y$,
for $y\in C$, are less than $\pi$.

Let $y\in C\setminus\{p\}$. Using Proposition~\ref{prop:cosine-comparison-domain} and the obtuse angle property \eqref{eqn:obtuse} from Proposition~\ref{lem:projection},
\begin{equation}\label{eq:projection-inequality}
\cos\dd(b,y)\leq\cos\dd(p,b)\cos\dd(p,y).
\end{equation}
This inequality also holds when $y=p$, with equality.

Unlike in \cite[Proposition~3]{gietl2024}, we cannot directly conclude that $\dd(p,y)\leq\dd(b,y)$ under our weaker support assumption.
Nevertheless, we can study the directional derivative of $F_\mu$ at $b$ towards $p$ (see~\eqref{eq:dir-derivative}). Using \eqref{eq:dir-derivative} and again
Proposition~\ref{prop:cosine-comparison-domain}, now applied at $b$ to $(b,p,y)$ (where $b\notin\{p,y\}$ by assumption), we obtain, since $\mu(C)=1$,
\[
\begin{aligned}
-DF_\mu(b)[\log_b p]
&=\int_C\langle\log_b p,\log_b y\rangle_b\,\dd\mu(y)\\
&\geq\frac{\dd(b,p)}{\sin\dd(b,p)}
  \int_C\frac{\dd(b,y)}{\sin\dd(b,y)}
  \left(\cos\dd(p,y)-\cos\dd(b,y)\cos\dd(b,p)\right)\dd\mu(y)\\
&\geq\frac{\dd(b,p)\sin\dd(b,p)}{\cos\dd(b,p)}J_\mu(b).
\end{aligned}
\]
Here $J_\mu$ is defined in \eqref{eq:def-Jmu}, and the last inequality follows from \eqref{eq:projection-inequality} and $\cos\dd(b,p)>0$. At the barycenter, \eqref{eq:afsari-inequality} and first-order optimality give $J_\mu(b)\geq\cos r>0$. Since $0<\dd(b,p)<\pi/2$, we obtain
$DF_\mu(b)[\log_b p]<0$, contradicting the first-order optimality condition at $b$.
\end{proof}

\subsection{Variance inequality}
In this section, we prove variance inequalities for the Fréchet
function on $\kappa$-domains.
When $\kappa>0$ and the support of $\mu$ is contained in a ball of radius $r<D_\kappa/4$, the Fréchet function $F_\mu$ is strongly convex on that ball. A quadratic variance inequality then follows directly from strong convexity and the minimality of the barycenter. This argument underlies \cite[Lemma 1 and Proposition 2]{brunelserres2025} In the larger regime $r<D_\kappa/2$, \cite[Proposition 19]{yokota2016} establishes 
\[
F_\mu(x)-F_\mu(b)
\geq C_{\kappa,r}\,d^\alpha(x,b),
\qquad \alpha>2,
\]
where $C_{\kappa,r}>0$ and $\alpha$ depend only on $\kappa$ and $r$. His proof uses Kendall's convex function from \cite[Theorem A]{yokota2016} and yields a weaker estimate than quadratic growth.
We mention also that quadratic growth inequalities for the 
Fréchet function are also known to be closely
related to extendibility of geodesics, as in~\cite{ahidar2020,legouic2022}.

The following theorem establishes a quadratic variance inequality throughout this larger radius regime, with an explicit geometric constant.

\begin{proposition}[Quadratic variance inequality]
\label{prop:quadratic-variance-cat-k}
Denote by $b:=b(\mu)\in D$ the unique barycenter of $\mu$.
For every $x\in D$,
\[
F_\mu(x)-F_\mu(b)
\geq
\frac{\cos(r\sqrt\kappa)}{\kappa}
\Bigl(1-\cos\bigl(\sqrt{\kappa}\,d(x,b)\bigr)\Bigr).
\]
Consequently, for every $x\in D$,
\[
F_\mu(x)-F_\mu(b)
\geq c_{\mathrm{var}}(r\sqrt\kappa)\,d^2(x,b),
\]
where for all $u\in[0,\pi/2)$
\[
c_{\mathrm{var}}(u):=
\cos(u)\,
\frac{1-\cos(2u)}{4u^2}
>0
\]
and $c_{\mathrm{var}}(0)$ is understood as $1/2$.
\end{proposition}

\begin{remark}[Zero-curvature limit]
For fixed $r>0$, $c_{\mathrm{var}}(r\sqrt{\kappa})\to\frac12$ as $\kappa\downarrow0$, consistently with the variance inequality in \cite[Proposition~4.4]{sturm2003}, which requires no bounded-support assumption.
\end{remark}

\begin{proof}

We work in the normalized case $\kappa=1$; the general case follows by rescaling the metric to $\sqrt{\kappa}\,\dd$. Apply Theorem~\ref{thm:directional-derivative-transfer} with first point \(b\) and second point \(x\). Since \(b\) minimizes \(F_\mu\), the
first-order condition gives
\[
    D F_\mu(b)[\log_{b}x]\geq0,
    \qquad
    D F_\mu(b)[\log_{b}o]\geq0.
\]
Consequently,
\[
F_\mu(x)-F_\mu(b)\geq
\frac{
    1-\cos d(x,b)
}{
    \cos d(b,o)
}
\cos r \geq
\cos r\bigl(1-\cos d(x,b)\bigr).
\]
Since \(d(x,b)\leq2r<\pi\) and the function
\[
    t\longmapsto\frac{1-\cos t}{t^2}
\]
is decreasing on \((0,\pi)\), it follows that
\[
    F_\mu(x)-F_\mu(b)
    \geq
    \cos r
    \frac{1-\cos(2r)}{4r^2}
    d^2(x,b).
\]
\end{proof}

\subsection{Wasserstein stability of the barycenter map}

In $\CAT(0)$ spaces, the fundamental contraction property of \cite[Theorem 6.3]{sturm2003} states that the barycenter map is $1$-Lipschitz with respect to the Wasserstein-$1$ distance:
\[
d(b(\mu),b(\nu))\leq\Wone(\mu,\nu).
\]
In the setting of positive curvature, \cite[Proposition 24]{yokota2016} establishes a related stability estimate via Kendall's convex function 
for measures concentrated in the same ball of radius
$r<\pi/(2\sqrt{\kappa})$:
\[
d(b(\mu),b(\nu))\leq C_{\kappa,r}\,W_p(\mu,\nu),
\]
where $C_{\kappa,r}<\infty$ and $p>1$ depend only on $\kappa$ and $r$. In the following
bound, we establish a stronger
stability result with respect to $W_1$,
throughout the same radius regime.

\begin{proposition}[Wasserstein-1 Lipschitz stability]
\label{prop:wasserstein-one-stability}
Let $\mu,\nu\in\mathcal P(\X)$ be supported on $D$, and denote by $b(\mu),b(\nu)\in D$ their respective barycenters. Then
\begin{equation}
\label{eq:wasserstein-one-stability}
\dd(b(\mu),b(\nu))
\leq c_{W_1}(r\sqrt\kappa)\,\Wone(\mu,\nu),
\end{equation}
where, for every $u\in[0,\pi/2)$, $c_{W_1}(u):=\frac{c_{4\bullet}(u)}{4c_{\mathrm{var}}(u)}$. Here $c_{4\bullet}(u)$ and $c_{\mathrm{var}}(u)$
are the constants in Propositions~\ref{prop:local-four-point} and~\ref{prop:quadratic-variance-cat-k}, respectively.
Thus, the barycenter map on probability measures supported on $D$ is $c_{W_1}(r\sqrt\kappa)$-Lipschitz with respect to $\Wone$ and $\dd$. 
\end{proposition}

\begin{remark}[Zero-curvature limit]
For fixed $r>0$, $c_{W_1}(r\sqrt\kappa)\to1$ as $\kappa\downarrow0$, consistently with \cite[Theorem~6.3]{sturm2003}, which requires no bounded-support assumption.
\end{remark}

\begin{proof}
We work in the normalized case $\kappa=1$; the general case follows by rescaling the metric to $\sqrt{\kappa}\,\dd$. Fix a coupling $\pi\in\Pi(\mu,\nu)$. Since both marginals assign
mass one to $B(o,r)$, the coupling assigns mass one to
$B(o,r)\times B(o,r)$. Set
\[
G(y):=d^2(b(\nu),y)-d^2(b(\mu),y).
\]
Applying the quantitative four-point inequality
Proposition~\ref{prop:local-four-point} in the domain
$D=B(o,r)$ of radius $r$ gives
\[
|G(y)-G(z)|
\leq
c_{4\bullet}(r)\,d(b(\mu),b(\nu))d(y,z)
\]
for $\pi$-almost every $(y,z)$. Consequently,
\begin{align}
&F_\mu(b(\nu))-F_\mu(b(\mu))
 +F_\nu(b(\mu))-F_\nu(b(\nu)) \notag\\
&\qquad=
\frac12\int_{\X\times\X}
\bigl(G(y)-G(z)\bigr)\,\dd\pi(y,z) \notag\\
&\qquad\leq
\frac{c_{4\bullet}(r)}2\,
d(b(\mu),b(\nu))
\int_{\X\times\X}d(y,z)\,\dd\pi(y,z).
\label{eq:coupling-stability-bound}
\end{align}
On the other hand, applying
Proposition~\ref{prop:quadratic-variance-cat-k} to $F_\mu$ at $b(\nu)$,
and to $F_\nu$ at $b(\mu)$, yields
\[
\begin{aligned}
F_\mu(b(\nu))-F_\mu(b(\mu))
+F_\nu(b(\mu))-F_\nu(b(\nu))
\geq
2c_{\mathrm{var}}(r)\,d^2(b(\mu),b(\nu)).
\end{aligned}
\]
Combining this bound with \eqref{eq:coupling-stability-bound} and taking the infimum
over $\pi\in\Pi(\mu,\nu)$ gives
$$22c_{\mathrm{var}}(r)\dist^2(b(\mu),b(\nu))\leq \frac{c_{4\bullet}(r)}{2}\dist(b(\mu),b(\nu))\Wone(\mu,\nu),$$
which proves \eqref{eq:wasserstein-one-stability}.
\end{proof}

\subsection{Polyak--\L{}ojasiewicz inequality for Fréchet functions
on \texorpdfstring{\boldmath$\kappa$}{kappa}-domains}
The {\it Polyak-\L{}ojasiewicz (PL) inequality}
is a natural
gradient domination condition which is significantly weaker
than strong convexity yet,
still implies many of the same optimization guarantees~\cite{polyak1964,lojasiewicz1963,karimi2016linear}.
In the context of barycenter estimation,
the PL inequality has been used to study algorithms for barycenter
estimation on Bures--Wasserstein space~\cite{chewi2020gradient,altschuler2021averaging}.
In this section, we apply Theorem~\ref{thm:directional-derivative-transfer}
to deduce a PL inequality for the Fréchet function.

\begin{definition}[Descending metric slope]
Following \cite[Definition 3.1]{ohtapalfia2015}, we denote by $|\nabla^-F_\mu|(x)$ the descending metric slope, also called the descending absolute gradient, defined by
\[
|\nabla^-F_\mu|(x)
:=
\max\left\{
0,\,
\limsup_{\substack{y\to x\\y\ne x}}
\frac{F_\mu(x)-F_\mu(y)}{d(x,y)}
\right\}.
\]
This is a nonnegative scalar quantity.
\end{definition}

The slope controls the directional derivative: for every $x,y\in B(o,r)$,
\begin{equation}
\label{eq:slope-directional-bound}
DF_\mu(x)[\log_x y]
\geq
-d(x,y)\,|\nabla^-F_\mu|(x).
\end{equation}
This is the observation following \cite[Definition 3.2]{ohtapalfia2015}.
Indeed, if $x\ne y$ and $\gamma$ is the constant-speed geodesic from $x$ to $y$, then $d(x,\gamma(t))=t\,d(x,y)$, so
\[
\begin{aligned}
-\frac{DF_\mu(x)[\log_x y]}{d(x,y)}
&=
\lim_{t\downarrow0}
\frac{F_\mu(x)-F_\mu(\gamma(t))}
     {d(x,\gamma(t))}
\\
&\leq |\nabla^-F_\mu|(x).
\end{aligned}
\]
The case $x=y$ is immediate.

\begin{theorem}[Polyak--\L{}ojasiewicz inequality]
\label{thm:pl-inequality}
Denote by $b:=b(\mu)\in D$ its unique barycenter.
Then, for every $x\in D$,
\[
F_\mu(x)-F_\mu(b(\mu))
\leq c_{\mathrm{PL}}(r\sqrt\kappa)\,
|\nabla^-F_\mu|(x)^2,
\]
where one may take for all $u\in[0,\pi/2)$
\[
c_{\mathrm{PL}}(u)
:=
\frac{1}{2\cos^2(u)}.
\]
\end{theorem}

\begin{remark}
For every fixed $r>0$, $c_{\mathrm{PL}}(r\sqrt\kappa)\to1/2$ as $\kappa\downarrow0$. This is the optimal constant in $\CAT(0)$ spaces, where the inequality follows from strong geodesic convexity of the squared distance \cite[Proposition~2.3, inequality~(2.3)]{sturm2003}, without bounded-support assumption, for every $\mu\in\mathcal P_2(\X)$.
\end{remark}

\begin{proof}
We work in the normalized case $\kappa=1$; the general case follows by rescaling the metric to $\sqrt{\kappa}\,\dd$.
Fix $x\in B(o,r)$ and set 
\[
h:=d(x,b),\qquad q:=d(x,o),\qquad
s:=|\nabla^-F_\mu|(x).
\]

From Theorem~\ref{thm:directional-derivative-transfer},
by directly applying~\eqref{eq:slope-directional-bound},
we obtain
\[
\begin{aligned}
F_\mu(x)-F_\mu(b)
&\leq
\bigl[\sin h+(1-\cos h)\tan q\bigr]s
-\frac{\cos r}{\cos q}(1-\cos h)\\
&=
\frac{2}{\cos q}
\left[
s\sin(h/2)\cos(h/2-q)
-\cos r\,\sin^2(h/2)
\right]\\
&\leq
\frac{2}{\cos q}
\left[
s\sin(h/2)-\cos r\,\sin^2(h/2)
\right]\\
&=:P\bigl(\sin(h/2)\bigr),
\end{aligned}
\]
where $P(u)=(2/\cos q)(su-u^2\cos(r))$ for all $u\in\R$.
Since $\cos r>0$ and $\cos q\geq\cos r$, we have $P(u)\leq\frac{s^2}{2\cos r\cos q}\leq\frac{s^2}{2\cos^2r}$ for all $u\in\R$, yielding the result.
\end{proof}

\subsection{Barycenters via tangent cones}
Our analysis of empirical barycenters in expectation
is inspired by the approach that is taken in \cite{bartl2025}, where
one relate the distance to the barycenter to averages of pseudo-inner products in the tangent cone.
However, the proof there
is specific to case where the tangent cone is a Hilbert space.
The purpose of this section is to extend this proof technique
to the case where the tangent cone may not be a vector space. In the rest of this section, we no longer assume that $\kappa>0$. In particular, when $\kappa\leq 0$, we allow the $\kappa$-domain to be the whole space $\X$.

\begin{lemma}[Dual norm identity in the tangent cone]
\label{lem:dual-norm-tangent-cone}
Let $(\X,\dist)$ be a $\CAT(\kappa)$ space, with $\kappa\in\R$, and let $b\in\X$. Let $\nu\in\mathcal P_2(C_b\X)$ and denote by $\beta$ its barycenter.
Then
\[
\sup_{\substack{u\in C_b\X\\ |u|_b\leq1}}
\int_{C_b\X}\langle u,v\rangle_b\,\diff\nu(v)
=|\beta|_b.
\]
\end{lemma}
We note that the existence and uniqueness of $\beta=b(\nu)$ is a consequence of \cite[Proposition~4.3]{sturm2003} because $C_b\X$ is a complete $\CAT(0)$ space, by Proposition~\ref{prop:tangent-cone-cat-zero}.
A proof of this result appears in 
Appendix~\ref{subsec:proof_dual_norm_identity}.

\begin{remark}
In a Hilbert space $H$, the barycenter is the mean $\beta=\int_H v\,\diff\nu(v)$, and the preceding identity reduces to the usual dual characterization of the norm
$\|\beta\| =\sup_{\|u\|\leq1}\langle u,\beta\rangle $.
\end{remark}

\begin{proposition}[Distance to the barycenter via tangent directions]
\label{prop:barycenter-tangent-process}
Let $(\X,\dist)$ be a $\CAT(\kappa)$ space with $\kappa\in\R$ and let $\mu\in\mathcal P_2(\X)$ be concentrated on a $\kappa$-domain $D$ of $\X$. Let $b=b(\mu)$ be the unique barycenter of $\mu$.

Then, for every $x\in D$,
\[
\dist(x,b)\leq c_{\rm tan}(r,\kappa)
\sup_{\substack{u\in C_x\X\\ |u|_x\leq1}}
\int_D\langle\log_x y,u\rangle_x\,\diff\mu(y),
\]
where
\[
c_{\rm tan}(r,\kappa):=
\begin{dcases}
1 & \mbox{if }\kappa\leq0,\\
\dfrac{r\sqrt{\kappa}}
{\sin(r\sqrt{\kappa})\cos^2(r\sqrt{\kappa})}
& \mbox{if }\kappa>0,
\end{dcases}
\]
and $r\in (0,\infty]$ is the radius of $D$.\footnote{Recall that $r$ can be infinite if $\kappa\leq 0$ and that $r<\pi/(2\sqrt\kappa)$ if $\kappa>0$.}
\end{proposition}

Existence and uniqueness of $b=b(\mu)$ follow from Theorem~\ref{thm:existence_uniqueness} if $\kappa>0$, or from \cite[Proposition~4.3]{sturm2003} if $\kappa\leq0$.

\begin{proof}
Suppose first that $\kappa>0$. We work in the normalized case $\kappa=1$; the general case follows by rescaling the metric to $\sqrt{\kappa}\,\dist$. In that case, note that the radius of $D$ satisfies $r<\pi/2$.
Set
\[
h:=\dist(x,b),\qquad
S:=\sup_{\substack{u\in C_x\X\\ |u|_x\leq1}}
\int_D\langle\log_x y,u\rangle_x\,\diff\mu(y).
\]
The case $h=0$ is immediate. As in
\eqref{eq:slope-directional-bound}, the directional derivative formula and the definition of $S$ give $DF_\mu(x)[\log_x z]\geq-\dist(x,z)S$ for every $z\in D$.
We then follow the proof of Theorem~\ref{thm:pl-inequality}, with $S$ in place of the descending slope. Theorem~\ref{thm:directional-derivative-transfer} gives
\[
F_\mu(x)-F_\mu(b)
\leq
\bigl[\sin h+(1-\cos h)\tan r\bigr]S
-\cos r(1-\cos h).
\]
By Proposition~\ref{prop:quadratic-variance-cat-k},
\[
F_\mu(x)-F_\mu(b)\geq\cos r(1-\cos h).
\]
Combining these inequalities and using the half-angle
identities, we obtain
\[
S\geq
\frac{2\cos r(1-\cos h)}
{\sin h+(1-\cos h)\tan r}
=
\frac{2\cos^2r\,\sin(h/2)}{\cos(r-h/2)}
\geq2\cos^2r\,\sin(h/2)
\geq\cos^2r\,\frac{\sin r}{r}\,h.
\]
where the last inequality follows from $h/2\leq r$ and the monotonicity of $\sin u/u$ on $(0,\pi)$. This proves the positive-curvature case.

Now, assume $\kappa\leq0$, and define $h$ and $S$ as above using the original metric $\dist$. The case $h=0$ is immediate. Otherwise, observe that $(\X,\dist)$ is a $\CAT(0)$ space,
so the Fr\'echet function is $1$-strongly geodesically convex by \cite[Proposition~2.3, inequality~(2.3)]{sturm2003}.

We give the details
here for the reader's convenience. Let $\gamma\in\Gamma(x,b)$ and note that strong geodesic convexity yields, for $t\in(0,1)$,
\[
F_\mu(\gamma(t)) \leq(1-t)F_\mu(x)+tF_\mu(b)
-\frac12t(1-t)h^2.
\]
Subtracting $F_\mu(x)$, dividing by $t$, and letting
$t\downarrow0$, we obtain
\[
DF_\mu(x)[\log_x b]
\leq F_\mu(b)-F_\mu(x)-\frac12h^2
\leq-h^2,
\]
where the last inequality is the variance inequality in $\CAT(0)$ spaces \cite[Proposition 4.4]{sturm2003}. Since $DF_\mu(x)[\log_x b]\geq-hS$, this yields $h\leq S$.
\end{proof}

\section{Consequences for barycenter estimation} \label{sec:consequences}
In this section, we give the consequences
of our previous results for barycenter estimation,
starting with empirical barycenters in Section~\ref{subsec:empirical_barycenters_consequences},
then moving to iterated barycenters in Section~\ref{subsec:iterated_barycenters_consequences}.
Throughout this section, let $(\X,\dd)$ be a separable $\CAT(\kappa)$ space, with $\kappa\in\R$.
Unless otherwise stated, we assume $\kappa>0$ and let $D=B(o,r)$ be a $\kappa$-domain, with $0<r<\pi/(2\sqrt\kappa)$.
    
\subsection{Empirical barycenters}
\label{subsec:empirical_barycenters_consequences}

Let $\mu\in\mathcal P(\X)$ be supported on $D$ and $X_1,\ldots,X_n$ be independent and identically distributed random variables with law $\mu$, where $n\geq 1$ is the sample size. In this section, we tackle the problem of estimating $b(\mu)$. We study two estimators: The empirical barycenter of $X_1,\ldots,X_n$ and an alternative version that is based on an iterative construction. 

Let us first focus on the empirical barycenter, defined as $\widehat b_n:=b(\mu_n)$, where $\mu_n=n^{-1}\sum_{i=1}^n \delta_{X_i}$ is the empirical measure associated with $X_1,\ldots,X_n$. That is, $\widehat b_n$ is the unique minimizer of the Fréchet function 
$$F_n(x):=F_{\mu_n}(x)=\frac{1}{2n}\sum_{i=1}^n \dist^2(x,X_i), \quad x\in \X.$$

Since $\mu$ is supported on $D$, so is $\mu_n$ almost surely, so $\widehat b_n$ is defined uniquely. Moreover, by Proposition~\ref{prop:wasserstein-one-stability}, $\widehat b_n$ is a Lipschitz function of the $n$-uple $(X_1,\ldots,X_n)$, hence it is a random variable.

The following theorem controls the expected statistical error for every $\kappa\in\R$. In particular, in $\CAT(0)$ spaces, it improves the constant $2$ in \cite[Theorem~3]{brunelserres2025} to the optimal constant $1$.

\begin{theorem}[Expected error for empirical barycenters]
\label{th:empirical-expected-error}
Let $(\X,\dist)$ be a $\CAT(\kappa)$ space, with $\kappa\in\R$, and let $X_1,\ldots,X_n$ be i.i.d.\ with law $\mu\in\mathcal P_2(\X)$, where $n\geq 1$ is an integer.
If $\kappa>0$, assume additionally that $\mu$ is concentrated
on a $\kappa$-domain of radius $r>0$.
Let $b:=b(\mu)$ and $\widehat b_n:=b(\mu_n)$, and set
$\sigma_\mu^2:=\E[\dist^2(X_1,b)]$. Then
\[
\E[\dist^2(\widehat b_n,b)]
\leq c_{\mathrm{emp}}(r,\kappa)\frac{\sigma_\mu^2}{n}
\]
where $c_{\mathrm{emp}}(r,\kappa):=c_{\mathrm{tan}}(r,\kappa)^2$, with $c_{\mathrm{tan}}$ defined in
Proposition~\ref{prop:barycenter-tangent-process}.
\end{theorem}
In the case where $\kappa > 0$,
this result is the first non-asymptotic rate for arbitrary
distributions on $\kappa$-domains,
which,
as we discuss in Remark~\ref{rem:tightness},
are the largest
possible domains where one can hope for $1/n$ convergence
bounds over arbitrary distributions.
The best previously known bounds were either restricted
to smaller domains with $r < \pi/(4\sqrt{\kappa})$,
where the Fréchet function
is strongly convex~\cite{brunelserres2025,legouic2022}, or achieved rates of the
form $n^{-\beta}$ where $\beta$
is a parameter depending on the metric entropy of the space~\cite{schotz2019convergence,ahidar2020}.

\begin{remark}[Zero-curvature
limit]
Recall that if $\kappa\leq 0$, then $c_{\textrm{tan}}(r,\kappa)=1$. Hence, in the $\CAT(0)$ case, the theorem improves the constant $2$ in \cite[Theorem~3]{brunelserres2025} to the optimal Hilbert-space constant $1$, addressing the nonpositive-curvature case of their Open Question~2.
Finally, we emphasize that
our tangent-cone argument, inspired by \cite{bartl2025}, requires no additional lower curvature bound, unlike \cite[Theorem~1.3]{legouic2022}.
\end{remark}

\begin{proof}
Suppose first that $\kappa>0$. We work in the normalized case $\kappa=1$; the general case follows by rescaling the metric to $\sqrt{\kappa}\,\dist$.
Set $V_i:=\log_bX_i$, and let $\beta_n$ be the barycenter in $C_b\X$ of $n^{-1}\sum_{i=1}^n\delta_{V_i}$. By Proposition~\ref{prop:barycenter-tangent-process} applied to $\mu_n$ at $b$, and Lemma~\ref{lem:dual-norm-tangent-cone},
\[
\dist(\widehat b_n,b)
\leq c_{\mathrm{tan}}(r,1)|\beta_n|_b.
\]
It therefore remains to show that $\E[|\beta_n|_b^2]\leq\sigma_\mu^2/n$.
Set $A_n(u):=n^{-1}\sum_{i=1}^n\langle u,V_i\rangle_b$.
The proof of Lemma~\ref{lem:dual-norm-tangent-cone} also gives, for every $u\in C_b\X$,
\[
A_n(u)\leq\langle u,\beta_n\rangle_b.
\]
Let $m_n$ be the point at time $1/(n+1)$ on the geodesic from $\beta_n$ to $V_{n+1}$. Then
\[
A_{n+1}(u)
\leq
\frac n{n+1}\langle u,\beta_n\rangle_b
+\frac1{n+1}\langle u,V_{n+1}\rangle_b.
\]
By Lemma~\ref{lem:dual-norm-tangent-cone}, applied to the empirical measure and to $\frac n{n+1}\delta_{\beta_n}+\frac1{n+1}\delta_{V_{n+1}}$,
\[
|\beta_{n+1}|_b
=\sup_{|u|_b\leq1}A_{n+1}(u)
\leq|m_n|_b.
\]
The $\CAT(0)$ inequality in $C_b\X$ therefore gives
\[
\begin{aligned}
|\beta_{n+1}|_b^2
&\leq
\frac n{n+1}|\beta_n|_b^2
+\frac1{n+1}|V_{n+1}|_b^2
-\frac n{(n+1)^2}
\dist_{C_b\X}^2(\beta_n,V_{n+1})\\
&=
\frac{n^2|\beta_n|_b^2+|V_{n+1}|_b^2
+2n\langle\beta_n,V_{n+1}\rangle_b}{(n+1)^2}.
\end{aligned}
\]
The dual norm identity and Cauchy--Schwarz give
\[
|\beta_n|_b^2\leq\frac1n\sum_{i=1}^n|V_i|_b^2,
\]
so $\E[|\beta_n|_b^2]\leq\sigma_\mu^2$. Consequently, $\langle\beta_n,V_{n+1}\rangle_b$
is integrable. First-order optimality at $b$ gives $\E[\langle u,V_1\rangle_b]\leq0$ for every $u\in C_b\X$. Since $\beta_n$ is measurable with respect to
$\mathcal F_n:=\sigma(X_1,\ldots,X_n)$ (by continuity of $\log_b$, see~\cite[Lemma~1]{brunelserres2025}, and the barycenter contraction property \cite[Theorem~6.3]{sturm2003}), independence yields
\[
\E[\langle\beta_n,V_{n+1}\rangle_b\mid\mathcal F_n]\leq0.
\]
Taking expectations, we obtain
\[
(n+1)^2\E[|\beta_{n+1}|_b^2]
\leq n^2\E[|\beta_n|_b^2]+\sigma_\mu^2.
\]
Since $\beta_1=V_1$ and $\E[|V_1|_b^2]=\sigma_\mu^2$, iteration gives
\[
\E[|\beta_n|_b^2]\leq\frac{\sigma_\mu^2}{n}.
\]
Combining this bound with
$\dist(\widehat b_n,b)\leq c_{\mathrm{tan}}(r,1)|\beta_n|_b$
proves the claim.

If $\kappa\leq 0$, the same argument applies without rescaling, with $c_{\mathrm{tan}}(r,\kappa)=1$, and gives \[
\E[\dist^2(\widehat b_n,b)]\leq\frac{\sigma_\mu^2}{n}.
\]
\end{proof}

\begin{corollary}[High probability bounds for empirical barycenters]
\label{cor:empirical-high-probability}
Let $(\X,\dist)$ be a $\CAT(\kappa)$ space, with $\kappa\in\R$, and let $X_1,\ldots,X_n$ be i.i.d. with law $\mu\in\mathcal P_2(\X)$, where $n\geq 1$ is an integer. Assume that $\mu$ is concentrated on a $\kappa$-domain of radius $r>0$. Let $b:=b(\mu)$ and $\widehat b_n:=b(\mu_n)$, and set $\sigma_\mu^2:=\E[\dist^2(X_1,b)]$. Then, for every $\delta\in(0,1)$, with probability at least
$1-\delta$,
\[
\dd(\widehat b_n,b)
\leq
\sqrt{c_{\mathrm{emp}}(r,\kappa)}
\sqrt{\frac{\sigma_\mu^2}{n}}
+
2r\,c_{\mathrm{emp,HP}}(r,\kappa)
\sqrt{\frac{\log(1/\delta)}{n}},
\]
where 
$$c_{\mathrm{emp,HP}}(r,\kappa)=\begin{cases}
c_{W_1}(r\sqrt{\kappa}) & \text{ if } \kappa>0 \\
1 & \text{ if } \kappa\leq 0
\end{cases}$$
and $c_{\mathrm{emp}}(r,\kappa)$ and $c_{W_1}(r\sqrt\kappa)$
are the constants appearing in Theorem~\ref{th:empirical-expected-error} and Proposition~\ref{prop:wasserstein-one-stability}, respectively. 
\end{corollary}

The best previous high-probability bounds
for general $\CAT(\kappa)$ spaces with $\kappa > 0$ were obtained in~\cite[Corollary~1]{brunelserres2025},
which required the domain to be small enough that the Fréchet function is strongly convex.
In Corollary~\ref{cor:empirical-high-probability}, we therefore relax the domain assumption to $\kappa$-domains, which
are the largest possible domains where
non-trivial convergence bounds can be expected, as we discuss in Remark~\ref{rem:tightness}.
When $\kappa\leq 0$, we choose to assume that $\mu$ has bounded support in this corollary for simplicity, even though more general assumptions on the tails of $\mu$ are studied in \cite{brunelserres2025}. Importantly, thanks to Theorem~\ref{th:empirical-expected-error}, Corollary~\ref{cor:empirical-high-probability} improves their $\CAT(0)$ bound by reducing the expected-error term to $\sigma_\mu/\sqrt n$.
We also mention that
high-probability bounds on empirical barycenters under lower bounds on the so-called
hugging function were obtained in both~\cite{legouic2022,bartl2025}.
In general, 
this latter condition restricts
the support to domains that are significantly
smaller than $\kappa$-domains.

\begin{proof}

When $\kappa\leq 0$, the same argument used in \cite[Corollary~1]{brunelserres2025} yields the desired result, with a better constant in front of the expected-error term thanks to Theorem~\ref{th:empirical-expected-error}. 
We also follow the sub-Gaussian concentration argument used in \cite[Corollary~1]{brunelserres2025} when $\kappa>0$.
By Proposition~\ref{prop:wasserstein-one-stability} and using the coupling that pairs corresponding observations,
\[
\Wone\left(
\frac1n\sum_{i=1}^n\delta_{x_i},
\frac1n\sum_{i=1}^n\delta_{y_i}
\right)
\leq \frac1n\sum_{i=1}^n\dd(x_i,y_i).
\]
Hence, the map
\[
(x_1,\ldots,x_n)\longmapsto\dd(\widehat b_n,b)
\]
is $c_{W_1}(r\sqrt\kappa)/n$-Lipschitz for the $\ell^1$-product distance. Since the observations are independent and supported on $B(o,r)$, the bounded
differences inequality implies that with probability at least $1-\delta$,
\[
\dd(\widehat b_n,b)
\leq \E[\dd(\widehat b_n,b)]
+2r\,c_{W_1}(r\sqrt\kappa)
\sqrt{\frac{\log(1/\delta)}{n}}.
\]
Finally, Cauchy--Schwarz inequality and Theorem~\ref{th:empirical-expected-error} give
\[
\E[\dd(\widehat b_n,b)]
\leq \sqrt{\E[\dd^2(\widehat b_n,b)]}
\leq
\sqrt{c_{\mathrm{emp}}(r,\kappa)}
\sqrt{\frac{\sigma_\mu^2}{n}},
\]
which proves the claim.
\end{proof}

\subsection{Iterated barycenters}
\label{subsec:iterated_barycenters_consequences}

Iterated barycenters are a natural estimator of
the population barycenter,
with the benefit over empirical barycenters
of being computationally
tractable as long as one can compute geodesics. They are defined as follows.

\begin{definition}[Iterated barycenter]
\label{def:iterated-barycenter}
Let $(\X,\dist)$ be a geodesic space and $n\geq 1$ be an integer. For $x_1,\ldots,x_n\in \X$ and weights $\lambda=(\lambda_2,\ldots,\lambda_n)\in(0,1)^{n-1}$, their iterated barycenter with weights $\lambda$ is $S_n$, defined recursively by
\[
S_1:=x_1,\qquad
S_k:=\gamma_k(\lambda_k),\quad k=2,\ldots,n,
\]
where $\gamma_k\in\Gamma(S_{k-1},x_k)$.
\end{definition}

While $\gamma_k$ is not necessarily unique, we will only work in setups where it is, at each step of the above construction. 

The iterated barycenter is an online estimator,
since each update uses only the previous iterate and the new observation. It can also be considered as a stochastic version of the proximal point method considered above, and in the case where $(\mathcal{X}, \dist)$ is a Riemannian manifold, it coincides exactly with 
stochastic gradient descent for the Fréchet function
$F_\mu$~\cite[Section~2.3]{brunelserres2025}.

For the canonical weights $\lambda_k=1/k$, this construction was studied by Sturm~\cite[Definition~4.6 and
Theorem~4.7]{sturm2003}, who proved a strong law of large numbers for bounded i.i.d.\ observations in $\CAT(0)$ spaces. In $\CAT(0)$ spaces, optimal bounds in expectation \cite[First claim in the proof of Theorem~4.7]{sturm2003} and with high probability \cite[Theorem 9]{brunelserres2025} are known, so from now on we only focus on the case of $\CAT(\kappa)$ spaces with $\kappa>0$.

We are interested in iterated barycenters with weights $\lambda_n:=A/(n-1+A)$, for $n\geq2$ and $A>0$, which are of the same form as those considered in \cite[Theorem~5]{brunelserres2025}. The following theorem establishes Lipschitz stability for iterated barycenters with a suitable choice of weights, thereby answering Open Question~1 in \cite{brunelserres2025}.
It is also basic to our high probability
bound for iterated barycenters in Corollary~\ref{cor:iterated-high-probability},
and is used in our 
proof of expected error of 
the iterated barycenter as well.

\begin{theorem}[Stability of iterated barycenters]
\label{thm:stability of it bar}
Let $(\X,\dist)$ be a $\CAT(\kappa)$ space with $\kappa>0$ and $D$ be a $\kappa$-domain. Let $n\geq 1$ be an integer, $x_1,\ldots,x_n,y_1,\ldots,y_n\in D$ and $\lambda=(\lambda_2,\ldots,\lambda_n)$ with 
$$\lambda_k:=\frac{A}{k-1+A}, \quad k=2,\ldots,n, \qquad A:=\frac{1}{\cos(r\sqrt\kappa)}.$$
Let $S_n^x$ be the iterated barycenter of $x_1,\ldots,x_n$ and $S_n^y$ that of $y_1,\ldots,y_n$, both with weights $\lambda$. Then, 
\[
\dd(S_n^x,S_n^y)
\leq
\frac{c_{\mathrm{it,Lip}}(r\sqrt\kappa)}n
\sum_{i=1}^n\dd(x_i,y_i),
\]
where, for every $u\in[0,\pi/2)$, one may take $c_{\mathrm{it,Lip}}(u):=2u/(\cos^3u\,\sin(2u))$, which is understood as $1$ if $u=0$.
\end{theorem}

\begin{remark}[Zero-curvature limit]
For fixed $r>0$, $c_{\mathrm{it,Lip}}(r\sqrt\kappa)\to1$ as $\kappa\downarrow0$, consistently with
\cite[Lemma~3.1]{Funano10} in $\CAT(0)$ spaces, which requires no bounded-domain assumption.
\end{remark}

\begin{proof}
We work in the normalized case $\kappa=1$; the general case follows by rescaling the metric to $\sqrt{\kappa}\,\dd$. First, the result is trivial for $n=1$, after noting that $c_{\mathrm{it,Lip}}(u)\geq 1$ for all $u\in [0,\pi/2)$. Let us now assume that $n\geq 2$. 

Set $L:=2r/(\cos r\sin(2r))$ in the following derivations.

For $p,q,u,v\in D$, $\gamma_{p,u}\in\Gamma(p,u)$, $\gamma_{p,v}\in\Gamma(p,v)$, $\gamma_{q,v}\in\Gamma(q,v)$ and $t\in[0,1]$, the triangle inequality gives
\begin{equation} \label{eqn:proof_lip12321}
\dd_o(\gamma_{p,u}(t),\gamma_{q,v}(t))\leq
\dd_o(\gamma_{p,v}(t),\gamma_{q,v}(t))+\dd_o(\gamma_{p,u}(t),\gamma_{p,v}(t)).
\end{equation}

The common-endpoint contraction~\eqref{prop:endpoint-inequality} bounds the first term by $(1-t\cos(r))\dd_o(p,q)$, while
Proposition~\ref{eq:L-convexity}, together with \eqref{eq:equivalence-inequality}, bounds the second term by $tL\dd_o(u,v)$. Hence, \eqref{eqn:proof_lip12321} becomes 
$$\dd_o(\gamma_{p,u}(t),\gamma_{q,v}(t)) \leq (1-t\,\cos(r))\dist_o(p,q)+ tL\dist_o(u,v).$$

Apply this bound to $p=S_{k-1}^x$, $q=S_{k-1}^y$, $u=x_k$, $v=y_k$ and $t=\lambda_k$ for all $k=2,\ldots,n$ yields
\begin{equation*}
    \dist_o(S_k^x,S_k^y) \leq \left(1-\frac{1}{k-1+A}\right)\dist_o(S_{k-1}^x,S_{k-1}^y)+\frac{AL}{k-1+A}\dist_o(x_k,y_k).
\end{equation*}
Now, multiplying both sides by $k-1+A$ yields
\begin{equation*}
    (k-1+A)\dist_o(S_k^x,S_k^y) \leq (k-2+A)\dist_o(S_{k-1}^x,S_{k-1}^y)+AL\dist_o(x_k,y_k).
\end{equation*}
Now, summing those inequalities from $k=2$ to $n$, we obtain
\begin{align*}
    (n-1+A)\dist_o(S_n^x,S_n^y) & \leq A\dist_o(S_1^x,S_1^y)+AL\sum_{k=2}^n \dist_o(x_k,y_k) \\
    & = A\dist_o(x_1,y_1)+AL\sum_{k=2}^n \dist_o(x_k,y_k) \\
    & \leq AL\sum_{i=1}^n \dist_o(x_i,y_i),
\end{align*}
after noting that $L\geq 1$. Therefore, we obtain that 
$$\dist_o(S_n^x,S_n^y) \leq \frac{AL}{n-1+A}\sum_{i=1}^n \dist_o(x_i,y_i)\leq \frac{AL}{n}\sum_{i=1}^n \dist_o(x_i,y_i)$$
since $n-1+A\geq n$, using that $A=1/\cos(r)\geq 1$. 
We finally obtain the desired result by recalling that $\dd\leq\dd_o\leq\dd/\cos r$ from \eqref{eq:equivalence-inequality} and that $c_{\mathrm{it,Lip}}(r)=A^2L$.
\end{proof}

We now turn to the expected error of the iterated barycenter. Yokota~\cite[Theorem~3]{yokota2018} established a strong law of large numbers for inductive barycenters in $\CAT(1)$ spaces under the same geometric assumptions as ours, using Kendall's function which is a function of two variables that is jointly geodesically convex. However, this function is comparable to $\dd^\alpha$ for an exponent $\alpha>2$, so the resulting estimates do not directly yield the desired rate for the expected squared error. Instead, we use a related function, whose form is attributed to J\"ager--Kaul \cite{jagerkaul1979} in \cite[Example~3]{kuwaesturm2006}, where its convexity is recalled
in the Riemannian setting. Given a $\CAT(\kappa)$ space $(\X,\dist)$ and a $\kappa$-domain $D=B(o,r)$, we define, for all $b\in D$, 
\[
\begin{aligned}
\Phi_b:B(o,r)&\longrightarrow[0,\infty),\\
x&\longmapsto
\frac{1-\cos\bigl(\sqrt{\kappa}\,\dd(x,b)\bigr)}
     {\kappa\cos\bigl(\sqrt{\kappa}\,\dd(x,o)\bigr)}.
\end{aligned}
\]

This function is closely related to the convex functions used in the study of harmonic maps and manifold-valued martingales. It is
geodesically convex (as a function of only one variable) and it comparable to $\dist^2(\cdot,b)$. As we show below, this choice yields almost sure convergence and an explicit bound of the form $C/n$ on the expected squared error for a suitable choice of step-sizes.

Related constructions appear in Kendall's work on martingales, harmonic maps, and convexity on hemispheres \cite{kendall1988,kendall1990,kendall1991}. 
Even a convexity estimate of the form $\E[\Phi_b(\gamma(t))]\leq(1-t)\Phi_b(x)+t\,\E[\Phi_b(Y)]$
would only give an additive term of order $t$, which is insufficient for the desired expected-error bound. The following lemma instead proves a instead exploits the specific form of $\Phi_b$ to obtain a one-step estimate with a remainder of order $t^2$.

\begin{lemma}[Descent in expectation]
\label{lem:iterated-exp}
Let $(\X,\dist)$ be a $\CAT(\kappa)$ space, with $\kappa>0$ and $D$ be $\kappa$-domain. Let $\mu\in\mathcal P(\X)$ with $\supp(\mu)\subseteq D$ and let $b=b(\mu)\in D$. Fix some $x\in D$ and let $Y$ be a random variable with law $\mu$. Denote by $\gamma\in\Gamma(x,Y)$ the almost surely unique, random geodesic from $x$ to $Y$. Then, for every $t\in(0,1)$,
\[
\E\bigl[\Phi_b(\gamma(t))\bigr]
\leq
\bigl(1-2\cos(r\sqrt\kappa)t\bigr)\Phi_b(x)
+
K(r\sqrt\kappa)\,\E[\dd^2(x,Y)]\,t^2.
\]
where
\[
K(u):=\frac{1}{2\cos u}
+\frac{\sin^2u}{\cos^3u}
+\frac{3u\sin u}{\cos^4u},
\qquad 0\leq u<\pi/2.
\]
\end{lemma}

\begin{proof}
We work in the normalized case $\kappa=1$. Indeed, rescaling the metric to $\sqrt{\kappa}\,\dd$ multiplies both $\Phi_b$ and the squared distance by $\kappa$; dividing the normalized inequality by $\kappa$ gives the general case.
Set
\[
c:=\cos r,\qquad \ell:=\dd(x,Y),\qquad
C:=\cos\dd(x,o),\qquad A:=1-\cos\dd(x,b).
\]
Define
\[
\alpha_t:=\frac{\sin((1-t)\ell)}{\sin\ell},
\qquad
\beta_t:=\frac{\sin(t\ell)}{\sin\ell},
\]
with continuous extensions when $\ell=0$.

The cosine inequality along geodesics, Proposition~\ref{prop:geodesic-cosine-inequality}, gives
\[
\begin{aligned}
1-\cos\dd(\gamma(t),b)
&\leq N(t)
:=1-\alpha_t\cos\dd(x,b)-\beta_t\cos\dd(Y,b),\\
\cos\dd(\gamma(t),o)
&\geq H(t)
:=\alpha_t\cos\dd(x,o)+\beta_t\cos\dd(Y,o).
\end{aligned}
\]
Since $H(t)\geq c>0$ (because $\dist(x,o),\dist(Y,o)\leq r<\pi/2$) and $N(t)\geq0$, we obtain
\[
\Phi_b(\gamma(t))\leq R(t):=\frac{N(t)}{H(t)}.
\]

For $z\in\{b,o\}$, denote by
\[
D_x(Y,z):=\frac{\ell}{\sin\ell}
\bigl[\cos\dd(Y,z)-\cos\ell\cos\dd(x,z)\bigr],
\]
with continuous extension when $\ell=0$.

Using $\alpha_t=\cos(t\ell)-\beta_t\cos\ell$, we obtain
\[
N(0)=A,\quad N'(0)=-D_x(Y,b),\qquad
H(0)=C,\quad H'(0)=D_x(Y,o).
\]
Consequently,
\[
R(0)=\Phi_b(x),\qquad R'(0)=-G(Y),
\]
where
\[
G(Y):=\frac{D_x(Y,b)}C+\frac{A\,D_x(Y,o)}{C^2}.
\]

We want to apply Taylor's formula to $R$. To do so, we need to control uniformly the remainder term, \emph{i.e.}, the second derivative of $R$. Set $U:=1-N$ and $W:=UH'-U'H$. Since $U''=-\ell^2U$ and $H''=-\ell^2H$, $W$ is constant. Differentiating $R=(1-U)/H$ gives
\[
R''=\frac{\ell^2}{H}
+\frac{2(H')^2}{H^3}
-\frac{2WH'}{H^3}.
\]

The definition of $H$ depends only on the side lengths of $(o,x,Y)$, not on the geodesic itself or the ambient space $\X$. We can therefore directly study $H$
and its derivative on the sphere. The triangle $(o,x,Y)$ has perimeter at most $4r<2\pi$.
Let $(\bar o,\bar x,\bar Y)$ be its spherical comparison triangle, and let $\bar\gamma\in\Gamma(\bar x,\bar Y)$.
Then
\[
H(t)=\cos\rho(t),
\qquad
\rho(t):=\dd_{\mathbb S^2}(\bar\gamma(t),\bar o).
\]

Since the spherical ball $B(\bar o,r)$ is geodesically convex, $\rho(t)\leq r$, and hence $H(t)\geq c$. For $\ell>0$ and $t<1$ such that $\rho(t)>0$, apply Proposition~\ref{prop:first-variation} to the geodesic $u\mapsto\bar\gamma(t+(1-t)u)$.
After dividing by $1-t$, this gives
\[
\rho'_+(t)
=-\ell\cos\angle_{\bar\gamma(t)}(\bar Y,\bar o).
\]
Consequently,
\[
H'(t)
=\ell\sin\rho(t)\cos\angle_{\bar\gamma(t)}(\bar Y,\bar o),
\qquad
|H'(t)|\leq\ell\sin r.
\]
If $\rho(t)=0$ or $\ell=0$, then $H'(t)=0$. The bound at $t=1$ follows by continuity.
Thus,
\[
H(t)\geq c,\qquad |H'(t)|\leq\ell\sin r.
\]

Now, $W$ is constant, so
\[
W=W(0)
=\frac{\ell}{\sin\ell}
\bigl[
\cos\dd(x,b)\cos\dd(Y,o)
-\cos\dd(Y,b)\cos\dd(x,o)
\bigr].
\]
Using the identity $ab-cd=(a-c)b+c(b-d)$, $|\cos(u)|\leq 1$, $\ell\leq2r$ and the monotonicity of $u/\sin u$ on $[0,\pi)$, we obtain
\[
\begin{aligned}
|W|
&\leq\frac{\ell}{\sin\ell}
\Bigl(
|\cos\dd(x,b)-\cos\dd(Y,b)|
+|\cos\dd(x,o)-\cos\dd(Y,o)|
\Bigr)\\
&\leq\frac{3\sin r\,\ell^2}{\sin\ell}\leq 3\ell\sin r\,\frac{2r}{\sin(2r)}=\frac{3r\ell}{c}.
\end{aligned}
\]
Consequently,
\[
R''\leq
\frac{\ell^2}{c}
+\frac{2\ell^2\sin^2r}{c^3}
+\frac{6r\ell^2\sin r}{c^4}
=2K(r)\ell^2.
\]
Taylor's formula with integral remainder therefore gives
\begin{equation}\label{eqn:proof_phib}
\Phi_b(\gamma(t))
\leq\Phi_b(x)-tG(Y)
+K(r)\,\dd^2(x,Y)t^2.
\end{equation}

The following estimate uses the integrated cosine-to-Fr\'echet transfer inequality~\eqref{eq:integrated-transfer-upper}, which also underlies Theorem~\ref{thm:directional-derivative-transfer} through angle comparison:
\[
\begin{aligned}
\E[G(Y)]
&\geq
\frac{F_\mu(x)-F_\mu(b)+\Phi_b(x)\cos r}
     {\cos\dd(x,o)}\\
&\geq
\frac{\cos r\bigl(1-\cos\dd(x,b)\bigr)}
     {\cos\dd(x,o)}
+\Phi_b(x)\cos r\\
&=2\cos r\,\Phi_b(x).
\end{aligned}
\]
Here the first inequality follows by expanding $G(Y)$ and using that $\dd(x,Y)/\sin\dd(x,Y)\geq1$ and that $\cos\dd(o,Y)\geq\cos r$ almost surely. The second follows from Proposition~\ref{prop:quadratic-variance-cat-k} and the fact that $0<\cos\dd(x,o)\leq1$.
Taking expectations in \eqref{eqn:proof_phib} then proves the normalized inequality
\[
\E[\Phi_b(\gamma(t))]
\leq(1-2\cos r\,t)\Phi_b(x)
+K(r)\E[\dd^2(x,Y)]\,t^2.
\]
\end{proof}

To get explicit convergence rate analysis for the expected error of iterated barycenter, we state a lemma similar to \cite[Lemma~5.8]{ohtapalfia2015}, which can be proved easily by induction. 

\begin{lemma}[Discrete recurrence estimates]
\label{lem:discrete-recurrence}
Let $(a_n)_{n\geq1}$ be a sequence of nonnegative numbers
satisfying
\[
a_{n+1}
\leq
\left(1-\frac{\alpha}{n+A}\right)a_n
+\frac{\beta}{(n+A)^2}, \qquad \forall n\geq 1,
\]
where $\alpha>1$, $\beta\geq0$, and $A\geq\alpha-1$.
Then, for every $n\geq1$,
\[
a_n\leq\frac{C}{n+A},
\qquad
C:=\max\left\{(1+A)a_1,\frac{\beta}{\alpha-1}\right\}.
\]
\end{lemma}

\begin{theorem}[Expected error for iterated barycenters]
\label{th:iterated-expected-error}
Let $(\X,\dist)$ be a $\CAT(\kappa)$ space, with $\kappa>0$ and $D$ be a $\kappa$-domain. Let $\mu\in\mathcal P(\X)$ with $\supp(\mu)\subseteq D$ and $X_1,\ldots,X_n$ be i.i.d. random variables with law $\mu$, for some integer $n\geq 1$. Let $b=b(\mu), \sigma_\mu^2=\E[\dist^2(X_1,b)]$ and $\widetilde b_n$ be the iterated barycenter of $(X_1,\ldots,X_n)$, with weights
\[
\lambda_k:=\frac{A}{k-1+A},\qquad 2\leq k\leq n,
\]
where we set $A:=1/\cos(r\sqrt\kappa)$. Then, for every $n\geq1$,
\[
\E[\dd^2(\widetilde b_n,b)]
\leq c_{\mathrm{E,it}}(r\sqrt{\kappa})
\frac{\sigma_\mu^2}{n+A},
\]
where, for every $u\in[0,\pi/2)$, one may take
\[
c_{\mathrm{E,it}}(u):=
\frac{4K(u)(c_{\mathrm{it,Lip}}(u)^2+1)}{\cos^4u},
\]
with $K$ and $c_{\mathrm{it,Lip}}$ defined in Lemma~\ref{lem:iterated-exp} and Theorem~\ref{thm:stability of it bar}, respectively.
\end{theorem}

In $\CAT(\kappa)$ spaces with
$\kappa > 0$, the only previously known bounds
for iterated barycenters
were in expectation and
for domains with radius $r < D_\kappa/4$, where
the squared distance is strongly convex~\cite{ohtapalfia2015,
brunelserres2025}, or in the particular case of 
the Bures--Wasserstein manifold~\cite{chewi2020gradient,altschuler2021averaging}.
Theorem~\ref{th:iterated-expected-error}
therefore relaxes the domain to 
$\kappa$-domains, which are the largest
possible as we discuss in Remark~\ref{rem:tightness}.

\begin{remark}[Zero-curvature limit]
For fixed $r>0$, $c_{\mathrm{E,it}}(r\sqrt{\kappa})\to4$ as $\kappa\downarrow0$. This is consistent with the variance-over-sample-size rate in \cite[proof of Theorem~4.7]{sturm2003}, which requires
no bounded-support assumption, but does not recover the optimal constant $1$. Obtaining this constant in the zero-curvature limit remains open for our argument, although we expect that a finer analysis could achieve this.
\end{remark}

\begin{proof}
We work in the normalized case $\kappa=1$; the general case follows by rescaling the metric to $\sqrt{\kappa}\,\dd$, which multiplies $\Phi_b$, the squared distances
and $\sigma_\mu^2$ by $\kappa$. Set $V_n:=\Phi_b(\widetilde b_n)$, $A:=1/\cos r$ and $L:=c_{\mathrm{it,Lip}}(r)$. By Theorem~\ref{thm:stability of it bar}, comparing with
the constant sample $(b,\ldots,b)$ and applying the Cauchy--Schwarz inequality, we obtain
\[
\E[\dd^2(\widetilde b_n,b)]\leq L^2\sigma_\mu^2.
\]
Consequently, by the triangle inequality,
\[
\E[\dd^2(\widetilde b_n,X_{n+1})]
\leq2(L^2+1)\sigma_\mu^2.
\]

To obtain the upper bound on the expected error, we apply Lemma~\ref{lem:iterated-exp} with $x=\widetilde b_n$,
$Y=X_{n+1}$, and $t=A/(n+A)$. Conditioning on $\mathcal F_n:=\sigma(X_1,\ldots,X_n)$
and then taking expectations gives
\[
\E[V_{n+1}]
\leq
\left(1-\frac{2}{n+A}\right)\E[V_n]
+\frac{2K(r)A^2(L^2+1)\sigma_\mu^2}{(n+A)^2}.
\]
Since $A\geq1$, Lemma~\ref{lem:discrete-recurrence}, with $\alpha=2$ and
$\beta=2K(r)A^2(L^2+1)\sigma_\mu^2$, yields
\[
\E[V_n]\leq
\frac{\max\{(1+A)\E[V_1],2K(r)A^2(L^2+1)\sigma_\mu^2\}}{n+A}.
\]
Since $\Phi_b(x)\leq\dd^2(x,b)/(2\cos r)$, $A\geq1$, $K(r)\geq1/2$ and $L\geq1$,
\[
(1+A)\E[V_1]
\leq\frac{A(1+A)}2\sigma_\mu^2
\leq A^2\sigma_\mu^2
\leq2K(r)A^2(L^2+1)\sigma_\mu^2.
\]
Thus,
\[
\E[V_n]\leq
\frac{2K(r)(L^2+1)}{\cos^2r} \frac{\sigma_\mu^2}{n+A}.
\]

Finally, the comparison of $\Phi_b$ with the squared distance and $\sin r\geq r\cos r$ give
\[
\E[\dd^2(\widetilde b_n,b)]\leq\frac{2}{\cos^2r}\E[V_n]\leq
\frac{4K(r)(L^2+1)}{\cos^4r}
\frac{\sigma_\mu^2}{n+A}=
c_{\mathrm{E,it}}(r)\frac{\sigma_\mu^2}{n+A}.
\]
This proves the claim.
\end{proof}

\begin{remark}
Lemma~\ref{lem:iterated-exp} gives the following conditional estimate:
$$\E[\Phi_b(\widetilde b_n)\mid\mathcal F_{n-1}]\leq (1-2\cos(r\sqrt\kappa)\lambda_n)\Phi_b(\widetilde b_{n-1})+4r^2K(r\sqrt\kappa)\lambda_n^2,$$
where we used $\dd(\widetilde b_{n-1},X_n)\leq2r$, and $\mathcal F_k:=\sigma(X_1,\ldots,X_k)$ is the $\sigma$-algebra generated by $X_1,\ldots,X_k$, for all $k\geq1$. Since $\sum_n\lambda_n=\infty$ and $\sum_n\lambda_n^2<\infty$, the Robbins--Siegmund theorem~\cite[Theorem~1]{robbinssiegmund1971} implies that $\widetilde b_n\xrightarrow[n\to\infty]{} b$ almost surely (for any choice of $A>0$ in the definition of $\widetilde b_n$). This recovers the strong law of large numbers for iterated barycenters proved by Yokota~\cite[Theorem~3]{yokota2018}.
\end{remark}

Combining Theorem~\ref{th:iterated-expected-error} with Theorem~\ref{thm:stability of it bar} yields the following high-probability bound for iterated barycenters.

\begin{corollary}[High-probability bound for iterated barycenters]
\label{cor:iterated-high-probability}
With the same notation and assumptions as in Theorem~\ref{th:iterated-expected-error}, for every $\delta\in(0,1)$, with probability at least $1-\delta$,
\[
\dd(\widetilde b_n,b)
\leq
\sqrt{\frac{c_{\mathrm{E,it}}(r\sqrt\kappa)\sigma_\mu^2}
           {n+A}}
+
2r\,c_{\mathrm{it,Lip}}(r\sqrt\kappa)
\sqrt{\frac{\log(1/\delta)}{n}}.
\]
Here $c_{\mathrm{E,it}}$ and $c_{\mathrm{it,Lip}}$ are the constants in Theorem~\ref{th:iterated-expected-error} and Theorem~\ref{thm:stability of it bar}, respectively.
\end{corollary}

\begin{proof}
Since the observations belong to $B(o,r)$, Theorem~\ref{thm:stability of it bar} shows that changing
one observation changes $\dist(\widetilde b_n,b)$ by at most $2r\,c_{\mathrm{it,Lip}}(r\sqrt\kappa)/n$.
Arguing as in the proof of
Corollary~\ref{cor:empirical-high-probability}, we obtain,
with probability at least $1-\delta$,
\[
\dd(\widetilde b_n,b)
\leq
\E[\dd(\widetilde b_n,b)]
+
2r\,c_{\mathrm{it,Lip}}(r\sqrt\kappa)
\sqrt{\frac{\log(1/\delta)}{n}}.
\]
Finally, the Cauchy--Schwarz inequality and
Theorem~\ref{th:iterated-expected-error} give
\[
\E[\dd(\widetilde b_n,b)]
\leq
\sqrt{\E[\dd^2(\widetilde b_n,b)]}
\leq
\sqrt{\frac{c_{\mathrm{E,it}}(r\sqrt\kappa)\sigma_\mu^2}
           {n+A}},
\]
which proves the claim.
\end{proof}

The only previously known high-probability
bounds for iterated barycenters
were in the $\CAT(0)$ case~\cite[Corollary~12 and Theorem~14]{brunelserres2024}. 
Corollary~\ref{cor:iterated-high-probability}
therefore shows the first known high-probabilty bounds 
for iterated barycenters on $\CAT(\kappa)$-spaces with
$\kappa > 0$, and in fact on the largest possible domains,
as discussed in the next remark.

\begin{remark}[Tightness of our results on barycenters]
\label{rem:tightness}
The rate $1/n$ in Theorem~\ref{th:empirical-expected-error}
and Theorem~\ref{th:iterated-expected-error}
and their respective corollaries
is optimal, as is immediately seen by specializing to the unit
interval
$(\X, \dist)= ([0,1], |\cdot|)$.
Beyond the rate, it is natural to ask if the restriction
$r < D_\kappa/2$ can be relaxed further.
However, as we now show,
any bound $\E[\dist^2(\widehat b_n, b)]
\leqslant \frac{c_\eps}{n}$,
valid uniformly over $\CAT(\kappa)$ spaces $(\X,\dist)$ and measures supported on $B(o, \pi/(2\sqrt\kappa) - \eps)$,
must have $c_\eps \gtrsim \frac{1}{\eps}$.
In particular, no $1/n$ bound with a finite constant can hold uniformly over this class when $r=D_\kappa/2$.

To see this, consider $(\mathcal{X}, \dist) = (\mathbb{S}^2, \dist_{\mathbb{S}^2})$, let $a \in [0,1]$, and let $\mu_a$ be the distribution supported on $\mathbb{S}^2$
defined by
$$
\mu_a := \frac12 (1 - a)(\delta_{x} + \delta_y) + a \delta_z.
$$ for $x = (\sqrt{1 - \eps^2}, 0, \eps), y = (-\sqrt{1 - \eps^2}, 0, \eps)$,
and $z = (0,\sqrt{1 - \eps^2}, \eps)$ with $\eps > 0$ small.
Then $x, y, z \in B(N, \arccos(\eps)) \subset B(N, \pi/2 - \eps)$ where $N = (0, 0, 1)$ denotes the north pole.
Let the barycenter of $\mu_a$ be $b_a$. A straightforward calculation
shows that $\dist_{\mathbb{S}^2}(N,
b_a) \geqslant \frac{ac}{\eps}$ so long
as $a \leqslant c' \eps$,
for constants $c,c' > 0$.
Let
$n = \lceil \frac1\eps \rceil$ and
$a =\frac{c'}{n}$,
and let $\widehat b_n$ denote the empirical barycenter
for a sample of size $n$ from $\mu_a$.
Then with
constant probability, the sample doesn't contain $z$ and $\dist_{\mathbb{S}^2}(N, \widehat b_n)
\leqslant \frac{2}{\sqrt{n}}$.
On this event,
$$
\dist_{\mathbb{S}^2}(\widehat b_n, b_a) \geqslant 
\dist_{\mathbb{S}^2}(N, b_a) - \dist_{\mathbb{S}^2}(N, \widehat b_n)
\geqslant \frac{ac}{\eps}-  \frac{2}{\sqrt{n}}.
$$ So long as $\eps$ is sufficiently small, this establishes
$\dist_{\mathbb{S}^2}(\widehat b_n, b_a)\gtrsim 1$ with constant probability.

Hence, any bound of the form $\mathbb{E}[\dist_{\mathbb{S}^2}(\widehat b_n, b_a)^2]
\leqslant \frac{c_\eps}{n}$ must satisfy
$$
\frac{c_\eps}{n} \gtrsim 1 \implies c_\eps \gtrsim n \gtrsim \frac1\eps.
$$ 
Therefore, the blow-up of the constant in Theorem~\ref{th:empirical-expected-error} is necessary.

An analogous argument applies to the iterated barycenter $\widetilde b_n$, showing that the restriction $r<D_\kappa/2$ is necessary for such uniform bounds and that the constant in Theorem~\ref{th:iterated-expected-error} must also blow up as $r\uparrow D_\kappa/2$.
\end{remark}

\subsection{Barycentric hulls}
\label{subsec:barycentric}


Let $(\X,\dist)$ be a geodesic space and let $A\subseteq\X$ be non-empty.
Recall that $\conv(A)$ denotes its geodesic convex hull.
If $(\X,\dist)$ is a $\CAT(\kappa)$ space and $A$ is contained in a $\kappa$-domain $D$, then $\conv(A)\subseteq D$, since $D$ is geodesically convex by Lemma~\ref{prop:unique-geodesics-domain}.

In Euclidean and Hilbert spaces, the convex hull of a non-empty set is the collection of all finite convex combinations of its elements or, equivalently, all
barycenters of probability measures finitely supported on that set. This characterization need not hold in non-Euclidean geometries.

\begin{definition}[Barycentric hull]
\label{def:barycentric-hull}
Let $(\X,\dist)$ be a metric space. The barycentric hull of a non-empty subset $A\subseteq\X$, denoted by $\bary(A)$, is the collection of all barycenters of probability measures finitely supported on $A$.
\end{definition}

The restriction to finitely supported measures ensures that $\bary(A)=\conv(A)$ in Hilbert spaces. Allowing arbitrary probability measures with finite second moment concentrated on $A$ may yield barycenters outside $\conv(A)$, though still in its closure.

Related constructions have been considered for triangulating Riemannian manifolds \cite{dyer2015riemannian,dyer2016barycentric} or extending principal component analysis to Riemannian contexts \cite{pennec2018barycentric}, and are well suited for applying Jensen's inequalities in non-linear spaces \cite{kuwae2014jensen}. Recently, algorithms were developed to provide membership oracles for barycentric hulls in spaces of non-positive curvature \cite{goodwin2025recognizing}.

In this section, we will only be interested in cases where $(\X,\dist)$ is a $\CAT(\kappa)$ space for some $\kappa\in\R$ and when $A$ is contained in a $\kappa$-domain, ensuring existence and uniqueness of barycenters of probability measures with support included in $A$. While the rest of the paper was dedicated to the case of $\kappa>0$, we also allow $\kappa\leq 0$ in this section, building on results of \cite{brunelserres2025}.

First, as a consequence of Proposition~\ref{prop:wasserstein-one-stability} on the $W_1$ stability of barycenters, we have the following result. In a metric space $(\X,\dd)$, we denote by $\dH$ the Hausdorff pseudo-distance, defined for any two nonempty
subsets $A,B\subseteq\X$ by
\[
\dH(A,B)
:=
\max\left\{
\sup_{x\in A}\dd(x,B),
\sup_{y\in B}\dd(y,A)
\right\}.
\]

\begin{theorem}
Let $(\X,\dist)$ be a $\CAT(\kappa)$ space for some $\kappa\in\R$, and let $A,B$ be two non-empty subsets of a bounded $\kappa$-domain of $\X$ of radius $r\geq0$.\footnote{If $\kappa>0$, all $\kappa$-domains are bounded by definition.}
Then,
\begin{align*}
\dH(\bary(A),\bary(B))
&\leq
\begin{cases}
\dH(A,B), & \text{if }\kappa\leq0,\\
c_{\mathrm{bh,Lip}}(r\sqrt{\kappa})\,\dH(A,B),
& \text{if }\kappa>0,
\end{cases}
\end{align*}
where $c_{\mathrm{bh,Lip}}(u):=c_{W_1}(u)$ for every $u\in[0,\pi/2)$. Here, $c_{W_1}(u)$ is the constant appearing in Proposition~\ref{prop:wasserstein-one-stability}.
\end{theorem}

\begin{proof}
    Let $\eta=\dH(A,B)$ and fix $\varepsilon>0$. Set $c:=1$ if $\kappa\leq0$, and
$c:=c_{\mathrm{bh,Lip}}(r\sqrt{\kappa})$ if $\kappa>0$. Let $x\in \bary(A)$. That is, there exists a measure $\mu$ of the form $\mu=\sum_{i=1}^n \omega_i\delta_{x_i}$, for some integer $n\geq 1$, some weights $\omega_1,\ldots,\omega_n\geq 0$ with $\omega_1+\ldots+\omega_n=1$ and some points $x_1,\ldots,x_n\in A$ such that $x=b(\mu)$. 
    For each $i=1,\ldots,n$, $\dist(x_i,B)\leq \eta$ so one can find $y_i\in B$ with $\dist(x_i,y_i)\leq\eta+\varepsilon$. Let $\nu=\sum_{i=1}^n \omega_i\delta_{y_i}$ and $y=b(\nu)$. Then, $y\in\bary(B)$. Using Proposition~\ref{prop:wasserstein-one-stability} if $\kappa>0$ or \cite[Theorem 6.3]{sturm2003} if $\kappa\leq 0$, we obtain that 
    $$\dist(x,y)=\dist(b(\mu),b(\nu))\leq c\,W_1(\mu,\nu).$$
    Now, choosing the coupling of $\mu$ and $\nu$ that simply pairs each $x_i$ with $y_i$, $i=1,\ldots,n$, we obtain 
    $$\dist(x,y)\leq c\,\sum_{i=1}^n \omega_i \dist(x_i,y_i)\leq c\,(\eta+\varepsilon).$$
    Since $\varepsilon>0$ was arbitrary, we thus have established
    that $\dist(x,\bary(B))\leq c\eta$. Similarly, we can also obtain that for all $y\in \bary(B)$, $\dist(y,\bary(A))\leq c\,\eta$, which then concludes the proof. 
\end{proof}

When $\kappa\leq 0$, the contraction constant is exactly $1$, and it is optimal in Euclidean spaces. In general, for $\kappa>0$, we do not know whether the constant is optimal.

Notably, a similar result can be proved for geodesic convex hulls, with a different constant.

\begin{theorem}\label{thm:Lip_conv_hulls}
Let $(\X,\dist)$ be a $\CAT(\kappa)$ space for some $\kappa\in\R$, and let $A,B$ be two non-empty subsets of a bounded $\kappa$-domain of $\X$ of radius $r\geq0$. Then,
\begin{align*}
\dH(\conv(A),\conv(B))
&\leq
\begin{cases}
\dH(A,B), & \text{if }\kappa\leq0,\\
c_{\mathrm{hull,Lip}}(r\sqrt{\kappa})\,\dH(A,B),
& \text{if }\kappa>0,
\end{cases}
\end{align*}
where, for every $u\in(0,\pi/2)$, $c_{\mathrm{hull,Lip}}(u):=u/(\sin u\cos u)$, understood as $1$ if $u=0$.
\end{theorem}

\begin{remark}
For every fixed $r>0$, $c_{\mathrm{hull,Lip}}(r\sqrt{\kappa})\to1$ as $\kappa\downarrow0$, recovering the bound in nonpositively curved spaces.
\end{remark}

The proof will heavily rely on the following facts, the first of which give a general, inductive construction of geodesic convex hulls.

\begin{lemma}[Convex hull representation
{\cite[Proposition 2.5.5]{papadopoulos2005metric}}] \label{lem:representation_conv_hull}
    Let $(\X,\dist)$ be a metric space and $A\subseteq \X$ be non-empty. Define the nested sequence of sets $(A_n)_{n\geq 0}$ inductively by setting $A_0=A$ and by letting, for all $n\geq 1$, $A_n$ be the collection of all points in $\X$ that are on a geodesic connecting two points of $A_{n-1}$. Then, 
    $$\conv(A)=\bigcup_{n\geq 0} A_n.$$
\end{lemma}

An important consequence of this is the following.

\begin{lemma} \label{lem:Jensen_conv}
Let $(\X,\dist)$ be a metric space and $A\subseteq \X$ be non-empty. Let $f:\X\to \R$ and assume that $f$ is geodesically convex on $\conv(A)$. Then, 
$$\sup_{x\in\conv(A)}f(x)\leq \sup_{x\in A}f(x).$$
\end{lemma}

\begin{proof}
    Using the notation of Lemma~\ref{lem:representation_conv_hull}, one can show, by induction, that for all $n\geq 1$ and all $x\in A_n$, $f(x)\leq\sup_{y\in A}f(y)$ and the conclusion follows immediately.
\end{proof}

The next two important results connect geodesic convexity in a metric space with geodesic convexity in its Euclidean cone. For simplicity, given two subsets $A,B$ of a metric space $\X$ with $A\subseteq B$, we identify the Euclidean cones $C_0A$ and $C_0B$ of $A$ and $B$ with subsets of $C_0\X$, satisfying the inclusion $C_0A\subseteq C_0B$. Moreover, given $A\subseteq \X$, we identify $\{1\}\times A$ with the set $\{1x:x\in A\}$, which is a subset of $C_0\X$.

\begin{lemma} \label{lem:conv_in_cone}
Let $(\X,\dist)$ be a $\CAT(\kappa)$ space, with $\kappa>0$, and let $D=B(o,r)$ be a $\kappa$-domain. Let $A\subseteq D$ be a non-empty set and let $C_0A=\{\lambda x:x\in A,\lambda\geq0\}$. We regard $C_0A$ as a subset of $C_0\X$, equipped with
the metric $\dist_{C_0\X}^{(\kappa)}$ from Definition~\ref{def:euclidean-cone}.
The following properties hold.
    \begin{itemize}
        \item $A$ is geodesically convex if and only if $C_0A$ is geodesically convex;
        \item $\conv(A)=\{x\in\X:\exists\lambda>0,\, \lambda x\in\conv(C_0A)\}$. Likewise, radial projection of $\conv(C_0A)\setminus\{0\}$ onto the unit section gives $\{1\}\times\conv(A)$;
        \item $\conv(\{1\}\times A)$ only contains elements of the form $tx$ with $x\in D$ and $t\geq\cos(\sqrt{\kappa}\,r)$.
    \end{itemize}
\end{lemma}

\begin{proof}[Proof of Lemma~\ref{lem:conv_in_cone}]
    We work in the normalized case $\kappa=1$; the general case follows by rescaling the metric to $\sqrt{\kappa}\,\dist$.
    The proof relies on the following fact, which can be found in the proof of \cite[Proposition I.5.10]{bridsonhaefliger1999}: Fix any two points $x,y\in D$ and any two positive numbers $s,u>0$. Let $\gamma:[0,1]\to \X$ be the unique geodesic from $x$ to $y$ in $\X$ and $\tilde\gamma:[0,1]\to C_0\X$ the unique geodesic from $sx$ to $uy$ in $C_0\X$. Then, $\gamma([0,1])=\{z\in D:\exists v>0,\ vz\in\tilde\gamma([0,1])\}$ (Proposition~\ref{prop:geodesics-cone} gives a more precise statement). That is, geodesics between nonzero points in $C_0D$ project onto geodesics in $D$, while geodesics from the vertex are radial segments. This proves the first assertion. For the second, the projected hull is convex by the same geodesic argument and contains $A$, so it contains $\conv(A)$. Conversely, $C_0\conv(A)$ is convex by the first assertion and contains $C_0A$, so the projected hull is contained in $\conv(A)$.
    For the last part of the lemma, since $D$ is geodesically convex, the first part shows that $C_0D$ is geodesically convex. Hence $\conv(\{1\}\times A)\subseteq C_0D$. Now, let $\gamma:t\geq 0\mapsto to\in C_0D$ (recall that $o$ is the center of the ball $D$). It is clear that $\gamma$ is a geodesic ray, since $\dist_{C_0\X}(\gamma(s),\gamma(t))=|s-t|$ for all $s,t\geq 0$. Since $C_0\X$ is a $\CAT(0)$ space, the Busemann function $b_\gamma:\lambda x\in C_0\X\mapsto \lim_{t\to\infty}\bigl(\dist_{C_0\X}(\lambda x,\gamma(t))-t\bigr)$ associated with $\gamma$ is geodesically convex on $C_0\X$ (see \cite[Definition II.8.17]{bridsonhaefliger1999} for the definition -- geodesic convexity follows from the geodesic convexity of $\dist$ in first argument). Now, one easily computes $b_\gamma(\lambda x)=-\lambda\cos(\dist(x,o))$ for all $x\in D$. Now, for all $x\in A$, $b_\gamma(1x)=-\cos(\dist(x,o))\leq -\cos(r)$, so by Lemma~\ref{lem:Jensen_conv}, all elements $\lambda x\in \conv(\{1\}\times A)$ must satisfy $b_\gamma(\lambda x)\leq -\cos(r)$, that is, $-\lambda\cos(\dist(x,o))\leq -\cos(r)$, implying that $\lambda\geq \cos(r)$.
\end{proof}

\begin{proof}[Proof of Theorem~\ref{thm:Lip_conv_hulls}]

    Let us first assume that $\kappa\leq 0$.
    Define two nested sequences of sets $(A_n)_{n\geq 0}$ and $(B_n)_{n\geq 0}$ inductively by setting $A_0=A$, $B_0=B$ and, for all $n\geq 1$, letting $A_n$ be the collection of all points in $\X$ that are on a geodesic connecting two points of $A_{n-1}$, and likewise for $B_n$. By Lemma~\ref{lem:representation_conv_hull}, $\conv(A)$ and $\conv(B)$ can be written as $\conv(A)=\bigcup_{n\geq 0} A_n$ and $\conv(B)=\bigcup_{n\geq 0} B_n$.
    Now, $\rho:=\dH(A,B)$ and let us check that $\dH(A_n,B_n)\leq \rho$, by showing that the sequence $(\dH(A_n,B_n))_{n\geq 0}$ is non-increasing. Let $n\geq 1$ and fix $x\in A_n$. Let $x_1,x_2\in A_{n-1}$ be such that $x=\gamma(t)$ for some $t\in [0,1]$, where $\gamma$ is the geodesic from $x_1$ to $x_2$. Fix $\varepsilon>0$ and let $y_1,y_2\in B_{n-1}$ such that $\dist(x_1,y_1)\leq \dH(A_{n-1},B_{n-1})+\varepsilon$ and $\dist(x_2,y_2)\leq \dH(A_{n-1},B_{n-1})+\varepsilon$. Let $\gamma'$ be the geodesic from $y_1$ and $y_2$. Then, $y:=\gamma'(t)\in B_n$ by definition, so $\dist(x,B_n)\leq \dist(x,y)$. Since $(\X,\dist)$ is a $\CAT(\kappa)$ space and $\kappa\leq 0$,  $\dist$ is jointly geodesically convex, yielding that $\dist(x,y)\leq (1-t)\dist(x_1,y_1)+t\dist(x_2,y_2)\leq \dH(A_{n-1},B_{n-1})+\varepsilon$. Therefore, $\dist(x,B_n)\leq \dH(A_{n-1},B_{n-1})+\varepsilon$. Similarly, for all $y\in B_n$, $\dist(y,A_n)\leq \dH(A_{n-1},B_{n-1})+\varepsilon$. Finally, we obtain the inequality $\dH(A_n,B_n)\leq \dH(A_{n-1},B_{n-1})+\varepsilon$, which must hold for all $\varepsilon>0$, yielding $\dH(A_n,B_n)\leq \dH(A_{n-1},B_{n-1})$. The inequality $\dH(\conv(A),\conv(B))\leq \dH(A,B)$ then follows readily, since the sequences $(A_n)_{n\geq 0}$ and $(B_n)_{n\geq 0}$ are nested.

    Let us now consider the case $\kappa>0$ and, by a scaling
    argument, we assume that $\kappa=1$ without loss of generality.
    Let $o\in\X$ such that $A$ and $B$ are included in the ball
    $B(o,r)$. This ball is a $\CAT(1)$ space itself, so its Euclidean cone $C$ is a $\CAT(0)$ space by Berestovskii's cone theorem \cite[Theorem II.3.14]{bridsonhaefliger1999}. We identify $C=C_0D$, where $D=B(o,r)$, with a convex subset of $C_0\X$, equipped with the restriction of $\dist_{C_0\X}$. Hausdorff distances between subsets of the cone are computed using this metric. Let $\tilde A=\{1x:x\in A\}$ and $\tilde B=\{1y:y\in B\}$. 
    
    Let $\pi:C_0\X\setminus \{0\}\to C_0\X$ be the map defined by setting $\pi(\lambda x)=1x$, for all $x\in \X$ and $\lambda>0$. 
    By the third part of Lemma~\ref{lem:conv_in_cone}, $\conv(\{1\}\times D)\subseteq C_0\X\setminus \{0\}$, so $\pi$ is well-defined on $\conv(\{1\}\times D)$. By the same geodesic projection argument as in the second part of Lemma~\ref{lem:conv_in_cone}, $\{1\}\times \conv(A)$ and $\{1\}\times \conv(B)$ can then be identified with $\pi(\conv(\tilde A))$ and $\pi(\conv(\tilde B))$ respectively. Moreover, one easily checks that $\dH(\tilde A,\tilde B)=2\sin(\dH(A,B)/2)$ and $\dH(\{1\}\times \conv(A),\{1\}\times \conv(B))=2\sin(\dH(\conv(A),\conv(B))/2)$, since $\dist_{C_0\X}(1x,1y)=2\sin(\dist(x,y)/2)$ for all $x,y\in B(o,r)$. 
    Denote by $\rho=\dH(A,B)$ and $\rho'=\dH(\conv(A),\conv(B))$. 

    By the first part of the proof (case $\kappa\leq 0$), we have that \begin{equation} \label{eqn:proof_ch1}
        \dH(\conv(\tilde A),\conv(\tilde B))\leq \dH(\tilde A,\tilde B)=2\sin(\rho/2).
    \end{equation}
    Moreover, 
    \begin{equation} \label{eqn:proof_ch2}
        2\sin(\dH(\conv(A),\conv(B))/2)=\dH(\pi(\conv(\tilde A)),\pi(\conv(\tilde B))).
    \end{equation}
    Let us now show that this restriction is $1/\cos(r)$-Lipschitz. Fix $x_1,x_2\in D$ and $\lambda_1,\lambda_2\geq \cos(r)$. Then,
    \begin{align*}
        \dist_{C_0\X}(\pi(\lambda_1x_1),\pi(\lambda_2x_2))^2 & = \dist_{C_0\X}(1x_1,1x_2)^2 \\
        & = 2-2\cos(\dist(x_1,x_2)) \\
        & = \frac{1}{\lambda_1\lambda_2} \left(2\lambda_1\lambda_2-2\lambda_1\lambda_2\cos(\dist(x_1,x_2))\right) \\
        & \leq \frac{1}{\lambda_1\lambda_2} \left(\lambda_1^2+\lambda_2^2-2\lambda_1\lambda_2\cos(\dist(x_1,x_2))\right) \\
        & = \frac{1}{\lambda_1\lambda_2} \dist_{C_0\X}(\lambda_1x_1,\lambda_2x_2)^2 \\
        & \leq \frac{1}{\cos(r)^2} \dist_{C_0\X}(\lambda_1x_1,\lambda_2x_2)^2
    \end{align*}
    where we have used the third part of Lemma~\ref{lem:conv_in_cone} in the last inequality.
    Therefore, the right-hand side of \eqref{eqn:proof_ch2} is bounded from above by $(1/\cos(r))\dH(\conv(\tilde A),\conv(\tilde B))$. Together with \eqref{eqn:proof_ch1}, we obtain that 
    \begin{equation} \label{eqn:proof_ch3}
        \sin(\rho'/2)\leq \frac{\sin(\rho/2)}{\cos(r)}.
    \end{equation}
    Since $A,B,\conv(A),\conv(B)\subseteq B(o,r)$, we have $\rho,\rho'\leq2r$. For $r>0$, the inequalities
    \[
    \frac{\sin r}{r}\,u\leq\sin u\leq u,
    \qquad 0\leq u\leq r,
    \]
    together with \eqref{eqn:proof_ch3}, give
    \[
    \frac{\sin r}{2r}\rho'
    \leq\sin(\rho'/2)
    \leq\frac{\sin(\rho/2)}{\cos r}
    \leq\frac{\rho}{2\cos r}.
    \]
    Hence
    \[
    \rho'\leq\frac{r}{\sin r\cos r}\rho.
    \]
    The case $r=0$ is immediate.
\end{proof}

By Proposition~\ref{prop:bary_in_convex_hull}, when $A$ is contained in a $\kappa$-domain of a $\CAT(\kappa)$ space, $\kappa\in\R$, $\bary(A)\subseteq \overline{\conv(A)}$. In general, however, $\bary(A)\neq \conv(A)$, except in specific cases such as Euclidean spaces, spheres and hyperbolic spaces, as seen in Theorem~\ref{thm:bary=conv} below. In fact, if $\X$ is a generic enough Riemannian manifold, then all geodesically convex subsets are either included in a geodesic line or full-dimensional (see \cite{lytchak2022every} for more details), whereas if $A$ is a finite subset of $\X$ with $p$ elements, then $\bary(A)$ is at most $(p-1)$-dimensional.

Yet, we show that barycentric hulls enjoy similar properties as convex hulls do in Euclidean spaces. One of these properties is the sparse representations of barycentric hulls. We give an extension of Carathéodory's theorem for exact, sparse representations of barycentric hulls in finite dimensional $\CAT$ spaces with curvature bounded from below. We also give a sparse, approximate representation of barycentric hulls in any $\CAT$ space. 

First, let us recall important properties of metric spaces with curvature bounds. In the sequel, we adopt \cite[Definition 8.2]{alexander2024alexandrov} to define curvature lower bounds -- In short, $(\X,\dist)$ is said to have curvature bounded from below if and only if all its triangles (with a perimeter restriction if needed) are at least as thick as in some fixed model space.

\begin{theorem}[{\cite[Corollary~7.11]{burago1992ad};
\cite[Corollary~15.16 and Theorem~8.14]{alexander2024alexandrov}}]
\label{thm:tangent-cone-dimension}
Let $(\X,\dist)$ be a $\CAT(\kappa)$ space with $\kappa\in\R$.
Assume that $(\X,\dist)$ has curvature bounded from below.
Then, for every $x\in\X$, $C_x\X$ has the same Hausdorff dimension as $\X$,
and this common dimension is either a nonnegative integer or infinite.
\end{theorem}

\begin{theorem}[{\cite[Lemma~4.6 and Theorem~4.11]{ohta2012barycenters};
\cite[Theorem~45]{yokota2012rigidity}}]
\label{thm:Yokota_Ohta}
Let $(\X,\dist)$ be a separable $\CAT(\kappa)$ space, with $\kappa\in\R$. Assume that $(\X,\dist)$ has curvature bounded from below. Let $\mu$ be a Borel probability measure supported on a bounded $\kappa$-domain $D$ of $\X$,
and let $b\in D$ be its unique barycenter. Then:
\begin{itemize}
    \item The tangent cone $C_b\X$ contains a subset $H$ isometric to a Hilbert space for the induced cone metric, with its origin identified with $0_b$ and $\log_b(\supp(\mu))\subseteq H$;
    \item For every $u\in C_b\X$,
    $\int_D\langle u,\log_b(x)\rangle_b\,\diff\mu(x)=0$.
    In particular, $\int_D\log_b(x)\,\diff\mu(x)=0_b$ in $H$;
    \item Any $b'\in D$ satisfying $\int_D\langle u,\log_{b'}(x)\rangle_{b'}\,\diff\mu(x)=0$ for every $u\in C_{b'}\X$ is equal to $b$.
\end{itemize}
\end{theorem}

The rest of this section is dedicated to the extension of Carathéodory's theorem to metric spaces with curvature bounds. Let us first recall Carathéodory's theorem.

\begin{theorem} [Carathéodory's theorem {\cite[Theorem~1.1.4]{schneider1993convex}}] \label{thm:classic_carat}
    Let $d\geq 1$ be an integer and let $A\subseteq \R^d$ be non-empty. Then, for all $x\in \conv(A)$, there exists $x_1,\ldots,x_{d+1}\in A$ with $x\in\conv(\{x_1,\ldots,x_{d+1}\})$.
\end{theorem}

This very classic theorem shows that in (finite dimensional) Euclidean spaces, convex hulls enjoy sparse representations. Here, we extend this theorem to barycentric hulls in finite dimensional metric spaces with curvature lower and upper bounds. Note that finite-dimensional spaces with curvature bounded from below are locally compact~\cite[Theorem~10.8.1]{burago2001course}. Since $\X$ is complete and geodesic, the Hopf--Rinow--Cohn-Vossen theorem \cite[Theorem~2.5.28]{burago2001course} implies that closed bounded balls are compact. Thus $\X$ is separable, being a countable union of compact balls.

\begin{theorem}[Metric version of Carathéodory's theorem] \label{thm:Carathéodory}

Let $(\X,\dist)$ be a $\CAT(\kappa)$ space for some $\kappa\in\R$ and with curvature bounded from below. Assume that $\X$ has finite Hausdorff dimension $d\geq 1$. Let $A\subseteq \X$ be a non-empty set included in a $\kappa$-domain. Then, for all $x\in\bary(A)$, there exists a measure $\mu$ with finite support of size at most $d+1$, included in $A$, whose barycenter is $x$. 
\end{theorem}

In other words, any element of $\bary(A)$ can be written as a weighted barycenter of at most $d+1$ points of $A$. That is, 
$$\bary(A)=\{b(\mu):\mu\in \mathcal P(\X), \supp(\mu)\subseteq A, \#\supp(\mu)\leq d+1\}.$$

\begin{proof}
    Let $x\in\bary(A)$. That is, there exist $x_1,\ldots,x_n\in A$ and weights $\lambda_1,\ldots,\lambda_n\geq 0$ with $\lambda_1+\ldots+\lambda_n=1$ such that $x$ is the unique minimizer of $y\in\X\mapsto\frac12\sum_{i=1}^n\lambda_i\dist(y,x_i)^2$. After discarding zero weights, we may assume that $\lambda_i>0$ for every $i$. If $\kappa\leq0$, the finite set $\{x,x_1,\ldots,x_n\}$ is contained in a bounded $\kappa$-domain, so Theorem~\ref{thm:Yokota_Ohta} applies. By Theorem~\ref{thm:Yokota_Ohta}, 
    $$0=\sum_{i=1}^n \lambda_i\log_x(x_i).$$
    That is, if $H$ is as in Theorem~\ref{thm:Yokota_Ohta}, $0$ is in the convex hull of $\log_x(x_1),\ldots,\log_x(x_n)$. Moreover, by Theorem~\ref{thm:Yokota_Ohta}, $H$ is isometrically embedded in $C_x\X$, whose Hausdorff dimension equals $d$ by Theorem~\ref{thm:tangent-cone-dimension}. Thus $\dim H\leq d$, and Carathéodory's theorem (Theorem~\ref{thm:classic_carat}) yields that one can extract at most $d+1$ vectors from $\{\log_x(x_1),\ldots,\log_x(x_n)\}$ whose convex hull in $H$ contains $0_x$, that is, 
    \begin{equation} \label{eqn:proof_metric_carat_181}
        0_x=\sum_{i=1}^{d+1}\omega_i \log_x(y_i)
    \end{equation}
    for some elements $y_1,\ldots,y_{d+1}\in\{x_1,\ldots,x_n\}$ and some weights $\omega_1,\ldots,\omega_{d+1}\geq0$ satisfying $\omega_1+\ldots+\omega_{d+1}=1$. 
    Denote by $\nu=\sum_{i=1}^{d+1}\omega_i \delta_{\log_x(y_i)}$. Then, $\nu \in \mathcal P_2(C_x\X)$ and it satisfies $\supp(\nu)\subseteq H$. Since $C_x\X$ is a $\CAT(0)$ space, $\nu$ has a unique barycenter $b(\nu)$. Moreover, since $H$ is isometric to a Hilbert space, it is a closed, geodesically convex subset of $C_x\X$. Therefore, since $C_x\X$ is a $\CAT(0)$ space,
    \cite[Proposition~6.1]{sturm2003} gives $b(\nu)\in H$. Now, \eqref{eqn:proof_metric_carat_181} yields that $0_x$ is the minimizer of $F_\nu$ on $H$ and, therefore, $b(\nu)=0_x$. Hence, letting $\mu_0:=\sum_{i=1}^{d+1}\omega_i\delta_{y_i}$, by Lemma~\ref{lem:dual-norm-tangent-cone}, 
    \begin{equation*}
    \sup_{u\in C_x\X, |u|_x\leq 1}
    \int_{\X}\langle u,\log_x(y)\rangle_x\,\diff\mu_0(y)
    =
    \sup_{u\in C_x\X, |u|_x\leq 1}
    \int_{C_x\X}\langle u,v\rangle_x\,\diff\nu(v)
    =0.
    \end{equation*}
    so Proposition~\ref{prop:barycenter-tangent-process} yields that $\dist(x,b(\mu_0))=0$, which concludes the proof. 

\end{proof}

The next corollary is an important consequence of Theorem~\ref{thm:Carathéodory}.

\begin{corollary}
    Let $(\X,\dist)$ be a finite-dimensional $\CAT(\kappa)$ space for some $\kappa\in\R$ with curvature bounded from below. Let $A\subseteq \X$ be a compact set that is included in some $\kappa$-domain. Then, $\bary(A)$ is compact.
\end{corollary}

\begin{proof}
    Set $d:=\dim_{\mathrm H}(\X)$, and let $D$ be a $\kappa$-domain containing $A$.
    Let $(x_n)_{n\geq1}$ be a sequence of elements of $\bary(A)$. By Theorem~\ref{thm:Carathéodory}, for each $n\geq1$, one can find $y_{1,n},\ldots,y_{d+1,n}\in A$ and weights $\lambda_{1,n},\ldots,\lambda_{d+1,n}\geq0$ with $\lambda_{1,n}+\ldots+\lambda_{d+1,n}=1$ such that $x_n$ is the unique minimizer of the map $F_n(y):=\frac12\sum_{i=1}^{d+1}
\lambda_{i,n}\dist(y,y_{i,n})^2$.
    The $d+2$ sequences $(y_{1,n})_{n\geq1},\ldots,(y_{d+1,n})_{n\geq1}$ and $((\lambda_{1,n},\ldots,\lambda_{d+1,n}))_{n\geq1}$ all take values in compact sets ($A$ for the first $d+1$ and the $d$-dimensional Euclidean simplex for the last), so one can extract a common converging subsequence.
    Up to renumbering, let us then simply assume that all $d+2$ sequences converge, and denote by $y_1,\ldots,y_{d+1}$ and $(\lambda_1,\ldots,\lambda_{d+1})$ their limits.
    Let $F(y):=\frac12\sum_{i=1}^{d+1}\lambda_i\dist(y,y_i)^2$.
    Since $y_1,\ldots,y_{d+1}$ are all in $A$, the measure $\sum_{i=1}^{d+1}\lambda_i\delta_{y_i}$ is supported on $D$, so its barycenter $b^*$ is the unique minimizer of $F$ on $\X$ and belongs to $D$.

    Moreover, $F_n$ converges uniformly to $F$ on all bounded subsets of $\X$. By Proposition~\ref{prop:bary_in_convex_hull}, all $x_n$ and $b^*$ belong to the bounded set $C:=\overline{\conv(A)}\subseteq D$.
    By the variance inequality, Proposition~\ref{prop:quadratic-variance-cat-k},\[
    c_{\mathrm{var}}(r\sqrt{\kappa})\,\dist^2(x_n,b^*)
    \leq F(x_n)-F(b^*)
    \leq 2\sup_{y\in C}|F_n(y)-F(y)|
    \longrightarrow0,
    \]
    where $r$ is the radius of $D$ if $\kappa>0$; for $\kappa\leq0$, the coefficient on the left-hand side is replaced by $1/2$, by \cite[Proposition~4.4]{sturm2003}. The second inequality uses $F_n(x_n)\leq F_n(b^*)$.
    Thus $x_n\to b^*\in\bary(A)$, which proves compactness of $\bary(A)$.
\end{proof}

We now extend a version of Carathéodory's theorem that gives sparse approximations of barycentric hulls, independently of the dimension of $\X$, in the same vein as \cite[Theorem~3.5]{de2019discrete} for normed vector spaces. In this next result, we no longer require $\X$ to have curvature lower bound or finite dimension. For a non-empty bounded set $A\subseteq\X$, we denote by c$\textrm{rad}(A):=\inf_{x\in A}\sup_{y\in A}\dist(x,y)$ its radius with center constrained to lie in $A$.

\begin{theorem}[Sparse approximations of barycentric hulls]
\label{th:sparse_approximation}
Let $(\X,\dist)$ be a $\CAT(\kappa)$ space for some $\kappa\in\R$ and let $A\subseteq\X$ be a non-empty set included in a bounded $\kappa$-domain of radius $r$. Fix $\varepsilon>0$.
Then, for every $x\in\bary(A)$, there exist $x_1,\ldots,x_N\in A$ with $N\leq\max\{1,\lceil
c_{\mathrm{Car}}(r,\kappa)\textrm{rad}(A)^2/\varepsilon^2 \rceil\}$ such that
$\dist(x,\bary(\{x_1,\ldots,x_N\}))\leq\varepsilon$, where $c_{\mathrm{Car}}(r,\kappa):=c_{\mathrm{emp}}(r,\kappa)$, with $c_{\mathrm{emp}}$ defined in Theorem~\ref{th:empirical-expected-error}.
\end{theorem}

\begin{proof}
    
    The proof relies on Theorem~\ref{th:empirical-expected-error} together with the probabilistic method. Let $x\in \bary(A)$. That is, there is some probability measure $\mu$ on $\X$ with (finite) support included in $A$, whose (unique) barycenter is $x$. Now, let $n\geq 1$ be any positive integer and let $X_1,\ldots,X_n$ be i.i.d random variables with distribution $\mu$. Set $\sigma_\mu^2:=\E[\dist^2(X_1,x)]$. By Theorem~\ref{th:empirical-expected-error}, the barycenter $\widehat b_n$ of $X_1,\ldots,X_n$ must satisfy $\E[\dist(\widehat b_n,x)^2]\leq c_{\mathrm{Car}}(r,\kappa)\sigma_\mu^2/n$. Therefore, with positive probability, $\dist(\widehat b_n,x)^2\leq c_{\mathrm{Car}}(r,\kappa)\sigma_\mu^2/n$. Moreover, all $X_i$'s, $i=1,\ldots,n$, belong to $A$ almost surely, since $\supp(\mu)\subseteq A$. Therefore, there exists a realization of $x_1,\ldots,x_n\in A$ whose barycenter is within distance $\sqrt{c_{\mathrm{Car}}(r,\kappa)\sigma_\mu^2/n}$ of $x$. Now, since $x$ minimizes the Fréchet function $F_\mu$ of $\mu$, we have that, for any choice $a\in A$, 
    $$\sigma_\mu^2 = 2F_\mu(x) \leq 2F_\mu(a)=\int_{\X} \dist^2(a,y)\diff\mu(y) \leq \sup_{y\in A}\dist^2(a,y)$$
    since $\dist(a,\cdot)\leq \sup_{y\in A}\dist(a,y)$ $\mu$-almost surely. Taking the infimum over $a\in A$, we obtain that $\sigma_\mu^2\leq\textrm{rad}(A)^2$.
    Choosing $n:=\max\{1,\lceil c_{\mathrm{Car}}(r,\kappa)\mathrm{rad}(A)^2/\varepsilon^2\rceil\}$, we obtain that $x$ is $\varepsilon$-close to the barycenter of at most $n$ points of $A$.
    \end{proof}

As is clear from the proof, in the theorem above, $x$ is in fact $\varepsilon$-close to the barycenter of $N$ points of $A$ with equal weights. 

In model spaces, the next theorem shows that barycentric and geodesic convex hulls coincide -- Our theorems above thus apply to geodesic convex hulls. 

\begin{theorem}[Equality of barycentric and convex hulls on model spaces] \label{thm:bary=conv}
    Let $(\X,\dist)$ be either a Euclidean space, a sphere or a hyperbolic space and $A\subseteq \X$ be non-empty. If $\X$ is a sphere, assume that $A$ is contained in an open hemisphere. Then, $\bary(A)=\conv(A)$. 
\end{theorem}

\begin{proof}
    When $\X$ is a Euclidean space, the result is classical. When $\X$ is a sphere, the result for finite sets has been shown in \cite[Theorem~7]{bussfillmore2001}, and the general case follows by taking unions over finite subsets of $A$. It stems from the fact that any geodesically convex subset of an open hemisphere can be written as the intersection of that hemisphere with a Euclidean convex cone (in the ambient Euclidean space). Let us focus on hyperbolic spaces. Without loss of generality, let us use their Beltrami-Klein representation. That is, let $\X$ be the open, Euclidean, unit ball, equipped with the metric
    $$\dist(x,y)=\arccosh\left(\frac{1-x^\top y}{\sqrt{(1-\|x\|^2)(1-\|y\|^2)}}\right), \quad x,y\in\X.$$
    where $x^\top y$ stands for the Euclidean inner product and $\|\cdot\|$ is the Euclidean norm. Then, geodesics in $\X$ are simply line segments, parameterized with non-constant Euclidean speed \cite[Proposition 5.28]{lee2018introduction}. Therefore, geodesically convex sets are exactly those subsets of $\X$ that are convex for the Euclidean geometry. Now, let $x_1,\ldots,x_n\in\X$, for some $n\geq1$. The case $n=1$ is immediate, so assume $n\geq2$. Let us show that $\conv(\{x_1,\ldots,x_n\})=\bary(\{x_1,\ldots,x_n\})$, which will yield the desired result. Let $y\in\bary(\{x_1,\ldots,x_n\})$. That is, by the first order condition, 
    \begin{equation} \label{eqn:proof_hyper1}
        \sum_{i=1}^n \omega_i \log_y(x_i)=0
    \end{equation}
    for some weights $\omega_1,\ldots,\omega_n\geq 0$ with $\omega_1+\ldots+\omega_n=1$. We identify $T_y\X$ with $\R^d$ through the Klein coordinates. Using again the characterization of geodesics in $\X$ as
    line segments, we have that $\log_y(x_i)=\lambda_i(x_i-y)$ for some $\lambda_i\geq 0$ that depends on $y$ and $x_i$. When $x_i\neq y$, $\lambda_i>0$ while when $x_i=y$, $\lambda_i$ can be chosen arbitrarily -- hence, we assume that all $\lambda_i$'s are positive. \eqref{eqn:proof_hyper1} then becomes 
    \begin{equation*} 
        \sum_{i=1}^n \omega_i\lambda_i (x_i-y)=0
    \end{equation*}
    which implies that $y$ belongs to the Euclidean convex hull
    of $x_1,\ldots,x_n$, and hence that $y\in\conv(\{x_1,\ldots,x_n\})$. This shows that $\bary(\{x_1,\ldots,x_n\})\subseteq \conv(\{x_1,\ldots,x_n\})$.
    
    Now, let $y\in \conv(\{x_1,\ldots,x_n\})$. Then, $y$ is in the Euclidean convex hull of $x_1,\ldots,x_n$, so $y$ must be a convex combination of those points, that is, 
    $$y=\sum_{i=1}^n \alpha_i x_i$$
    for some $\alpha_1,\ldots,\alpha_n\geq0$ with $\alpha_1+\ldots+\alpha_n=1$. That is, 
    $$\sum_{i=1}^n \alpha_i(x_i-y)=0.$$
    For $i=1,\ldots,n$, let $\tilde\omega_i=\alpha_i/\lambda_i$ and $\omega_i=\tilde\omega_i/(\sum_{j=1}^n\tilde\omega_j)$; where $\lambda_i$ is the same as defined above, so $\log_y(x_i)=\lambda_i(x_i-y)$. Then, 
    $$\sum_{i=1}^n \omega_i\lambda_i(x_i-y)=0$$
    that is,
    $$\sum_{i=1}^n \omega_i\log_{y}(x_i)=0$$
    which ensures that $y$ is the barycenter of the probability
    measure $\sum_{i=1}^n\omega_i\delta_{x_i}$, since every critical
    point of the Fr\'echet function is a minimizer in the hyperbolic space.
    That is, $y\in\bary(\{x_1,\ldots,x_n\})$. Since every point of either hull belongs to the corresponding
    hull of a finite subset of $A$, this proves the general case.
\end{proof}

\appendix

\section{Omitted proofs from Sections~\ref{sec:background}
and~\ref{sec:new_geometric_results}}

\subsection{Proof of cosine comparison without perimeter restriction}
\label{subsec:cosine_comparison_without_perimeter_restriction}
\begin{proof}[Proof of Proposition~\ref{prop:cosine-comparison-domain}]
We work in the normalized case $\kappa=1$; the general case follows by rescaling the metric to $\sqrt{\kappa}\,\dd$. Since $D$ has diameter strictly smaller than $\pi$,
all three distances are smaller than $\pi$,
and $\sin\dd(x,y)\sin\dd(x,z)>0$.

If $\dd(x,y)+\dd(x,z)+\dd(y,z)<2\pi$, the spherical law of cosines and
Proposition~\ref{prop:angle-comparison} give
\[
\frac{
\cos\dd(y,z)-\cos\dd(x,y)\cos\dd(x,z)
}{
\sin\dd(x,y)\sin\dd(x,z)
}
=
\cos\angle_x^{(1)}(y,z)
\leq\cos\angle_x(y,z).
\]

Otherwise, if $0<2\pi-\dd(x,y)-\dd(x,z)\leq\dd(y,z)<\pi$, then by the monotonicity of $\cos$ on $[0,\pi]$,
\[
\begin{aligned}
\cos\dd(y,z)
&\leq\cos\bigl(2\pi-\dd(x,y)-\dd(x,z)\bigr)\\
&=\cos\bigl(\dd(x,y)+\dd(x,z)\bigr)\\
&=\cos\dd(x,y)\cos\dd(x,z)
-\sin\dd(x,y)\sin\dd(x,z).
\end{aligned}
\]
Consequently,
\[
\frac{
\cos\dd(y,z)-\cos\dd(x,y)\cos\dd(x,z)
}{
\sin\dd(x,y)\sin\dd(x,z)
}
\leq-1
\leq\cos\angle_x(y,z).
\]
\end{proof}
\subsection{A space which is not
\texorpdfstring{\boldmath$\CAT(\kappa)$}{CAT(kappa)}
yet satisfies a four-point inequality}
\label{subsec:l_4_example}

Suppose $(\X, d) = (\R^2, \|\cdot\|_4)$.
Setting $N:=\|\cdot\|_4^2$, a direct computation gives $0\preceq \nabla^2N(u)\preceq 6I_2$ for all $u\neq0$, where $I_2$ is the $2\times 2$ identity matrix. Since $\nabla N$ extends continuously to the origin, it is globally $6$-Lipchitz for the Euclidean norm.
Integrating the difference of gradients along $[x_2,x_1]$ and using that $\|\cdot\|_2\leq 2^{1/4}\|\cdot\|_4$, we obtain 
\[
\begin{aligned}
&\bigl|N(x_1-x_3)-N(x_2-x_3)-N(x_1-x_4)+N(x_2-x_4)\bigr|\leq6\sqrt{2}\,\|x_1-x_2\|_4\|x_3-x_4\|_4
\end{aligned}
\]
for all $x_1,x_2,x_3,x_4\in\R^2$.
However, by \cite[Proposition~II.1.14]{bridsonhaefliger1999} and a rescaling argument, $(\R^2,\|\cdot\|_4)$ is not $\CAT(\kappa)$ for any $\kappa\in\R$, even locally.

\subsection{Proof of cosine inequality along geodesics} \label{subsec:Proof_cos_ineq}
In this section, we prove Proposition~\ref{prop:geodesic-cosine-inequality}.
To this end, we show the following result.

\begin{proposition}[Geodesics in the Euclidean cone]
\label{prop:geodesics-cone}
Let $(\X,\dd)$ be a $\CAT(\kappa)$ space, with $\kappa>0$, and let $C_0\X$ be its Euclidean cone. Let $y,z\in\X$ satisfy $\dd(y,z)<D_\kappa$, and let $\gamma\in\Gamma(y,z)$.
Let $\alpha_t, \beta_t$ be as defined in
Proposition~\ref{prop:geodesic-cosine-inequality}.
Then the unique geodesic $q$ in $(C_0\X,\dist_{C_0\X}^{(\kappa)})$ from $1y$ to $1z$ satisfies, for every $t\in[0,1]$,
\[
q\left(\frac{\beta_t}{\alpha_t+\beta_t}\right)
=
\frac{1}{\alpha_t+\beta_t}\,\gamma(t).
\]
\end{proposition}









\begin{proof}
We work in the normalized case $\kappa=1$; the general case follows by rescaling the metric to $\sqrt{\kappa}\,\dd$. The case $y=z$ is immediate, so assume $0<\dd(y,z)<\pi$. Choose $\bar y,\bar z\in\Sp^2$ such that
$\dd_{\Sp^2}(\bar y,\bar z)=\dd(y,z)$, and let $\bar\gamma\in\Gamma(\bar y,\bar z)$.
Then, for all $s,t\in[0,1]$,
\[
\dd(\gamma(s),\gamma(t))
=|s-t|\dd(y,z)
=\dd_{\Sp^2}(\bar\gamma(s),\bar\gamma(t)).
\]
The spherical interpolation formula \cite{shoemake1985} gives
\[
\bar\gamma(t)=\alpha_t\bar y+\beta_t\bar z.
\]
Since the Euclidean cone over $\Sp^2$ can be isometrically identified with $\R^3$, its geodesic joining $\bar y$ and $\bar z$ is simply the segment $\bar q(s)=(1-s)\bar y + s \bar z, s\in [0,1]$. Thus, setting $\sigma(t) :=\beta_t/(\alpha_t+\beta_t)$ and $\rho_t:=1/(\alpha_t+\beta_t)$,
we obtain 
\[
\bar q(\sigma(t))
=\rho_t\,\bar\gamma(t).
\]
This proves the formula in the spherical case and explains its origin. We now transfer it to the $\CAT(1)$ case.
Differentiating $\sigma$, one sees that it is strictly increasing and
therefore is a bijection from $[0,1]$ onto itself. We may therefore define a path $\widetilde\gamma:[0,1]\to C_0\X$ by setting 
\[
\widetilde\gamma(\sigma(t))
:=\rho_t\,\gamma(t), t\in [0,1].
\]
By the cone distance formula and the equality of distances along $\gamma$ and $\bar\gamma$, we have
\[
\begin{aligned}
&\dd_{C_0\X}\bigl(
\widetilde\gamma(\sigma(t)),
\widetilde\gamma(\sigma(s))\bigr)\\
&\qquad=
\|\bar q(\sigma(t))-\bar q(\sigma(s))\|_{\R^3}\\
&\qquad=
|\sigma(t)-\sigma(s)|\,\|\bar y-\bar z\|_{\R^3}\\
&\qquad=
|\sigma(t)-\sigma(s)|\,\dd_{C_0\X}(1y,1z).
\end{aligned}
\]
Hence $\widetilde\gamma$ is a geodesic joining $1y$ to $1z$. Uniqueness of geodesics in $\CAT(0)$ space therefore gives $\widetilde\gamma=q$, proving the claim.
\end{proof}

\begin{proof}[Proof of Proposition~\ref{prop:geodesic-cosine-inequality}]
We work in the normalized case $\kappa=1$; the general case follows by rescaling the metric to $\sqrt{\kappa}\,\dd$. The case $y=z$ is immediate, so assume $y\neq z$. When $\dd(x,y)+\dd(y,z)+\dd(x,z)<2\pi$, a proof using a comparison triangle is straightforward, by the $\CAT(1)$ inequality and the spherical law of cosines (see~\cite[Corollary~2.2]{KimuraSato2012}). However, as for the angle comparison in Proposition~\ref{prop:cosine-comparison-domain}, we want to remove this perimeter assumption. 

To do so, we will use the Euclidean cone $C_0\X$, which, equipped with $\dist_{C_0\X}$ (recall that we assume that $\kappa=1$ here), is a $\CAT(0)$ space, and, hence, whose squared distance to any given point is $2$-geodesically convex, by \cite[Proposition~2.3]{sturm2003}.

Let $q$ be the geodesic joining $1y$ to $1z$. By Proposition~\ref{prop:geodesics-cone}, we have $q(\sigma(t))=\rho_t\gamma(t)$ with $\sigma(t) =\beta_t/(\alpha_t+\beta_t)$ and $\rho_t=1/(\alpha_t+\beta_t)$, for all $t\in [0,1]$.
Using the strong geodesic convexity of $\dist_{C_0\X}^2(1x,\cdot)$,
\[
\dd_{C_0\X}^2(1x,\rho_t\gamma(t)) \leq (1-\sigma(t))\dd_{C_0\X}^2(1x,1y) + \sigma(t)\dd_{C_0\X}^2(1x,1z) - \sigma(t)(1-\sigma(t))\dd_{C_0\X}^2(1y,1z).
\]
All distances between points of $D$ are less than $\pi$. Expanding the cone distances therefore gives
\[
\begin{aligned}
1+\rho_t^2-2\rho_t\cos\dd(x,\gamma(t))
&\leq
2-2(1-\sigma(t))\cos\dd(x,y)-2\sigma(t)\cos\dd(x,z)\\
&\quad \quad \quad \quad \quad -2\sigma(t)(1-\sigma(t))
\bigl(1-\cos\dd(y,z)\bigr).
 \end{aligned}
\]
One can check that $\alpha_t^2+\beta_t^2 +2\alpha_t\beta_t\cos\dd(y,z)=1$, either using trigonometric formulas, or using the fact that, in the Euclidean space $\R^3$, $\alpha_t \bar y+\beta_t \bar z\in\mathbb S^2$, so $\alpha_t^2+\beta_t^2 +2\alpha_t\beta_t\cos\dd(y,z)=\|\alpha_t\bar y+\beta_t\bar z\|_{\R^3}^2=1$, where $\bar y,\bar z\in\mathbb S^2$ are comparison points on $\mathbb S^2$ for $y$ and $z$, as in the proof of Proposition~\ref{prop:geodesics-cone}. Dividing by $(\alpha_t+\beta_t)^2$, we obtain
\[
1+\rho_t^2
=
2-2\sigma(t)(1-\sigma(t))
\bigl(1-\cos\dd(y,z)\bigr).
\]
Substituting this equality into the preceding inequality and cancelling the common terms gives
\[
\rho_t\cos\dd(x,\gamma(t))
\geq
(1-\sigma(t))\cos\dd(x,y)
+\sigma(t)\cos\dd(x,z).
\]
Finally, dividing by $\rho_t$ yields the desired inequality.
\end{proof}

\subsection{Proof of endpoint geodesic stability in \texorpdfstring{$\kappa$}{kappa}-domains} \label{subsec:endpoint_inequality}

\begin{proof}[Proof of Lemma~\ref{lem:stability_same_endpoint}]
    
We work in the normalized case $\kappa=1$; the general case follows by rescaling the metric to $\sqrt{\kappa}\,\dd$. The case $t=1$ is immediate. If $x=z$, then $\dist(\gamma(t),\eta(t))=(1-t)\dist(x,y)$, so the result follows. We therefore assume that $x\neq z$. Finally, assume that $y\in D\setminus\{x\}$ is close enough to $x$ so $y\neq z$ and $\dist(x,y)+\dist(x,z)+\dist(y,z)<2\pi$ (this is possible, since the left hand side goes to $2\dist(x,z)$ as $y\to x$, and $2\dist(x,z)\leq 4r<2\pi$ by assumption). 
    
    For $t\in [0,1)$, let $\alpha(t)=\angle_z^{(1)}(\gamma(t),\eta(t))$. By \cite[Proposition II.3.1]{bridsonhaefliger1999}, $\alpha(t)$ is a non-increasing function of $t$. Moreover, as seen in Definition~\ref{def:comparison_angles}, we have, for $t\in [0,1)$, 
    \begin{equation} \label{eqn:proof_step254}
        \cos(\alpha(t))= \frac{\cos(\dist(\gamma(t),\eta(t)))-\cos((1-t)\dist(x,z))\cos((1-t)\dist(y,z))}{\sin((1-t)\dist(x,z))\sin((1-t)\dist(y,z))}.
    \end{equation}
    Rearranging \eqref{eqn:proof_step254}, we obtain that for all $t\in [0,1)$, 
\begin{align} 
    1 -\cos(\dist(\gamma(t),\eta(t))) & = 1-\cos((1-t)\dist(x,z))\cos((1-t)\dist(y,z)) \nonumber \\
    & \qquad -\sin((1-t)\dist(x,z))\sin((1-t)\dist(y,z))\cos(\alpha(t)) \nonumber \\
    & = 1-\cos((1-t)(\dist(x,z)-\dist(y,z))) \nonumber \\
    & \qquad +\sin((1-t)\dist(x,z))\sin((1-t)\dist(y,z))(1-\cos(\alpha(t))) \nonumber \\
    & = 2\sin^2\left(\frac{1-t}{2}(\dist(x,z)-\dist(y,z))\right) \label{eqn:proof_54231} \\
    & \qquad +2\sin((1-t)\dist(x,z))\sin((1-t)\dist(y,z))\sin^2(\alpha(t)/2) \nonumber \\
    & \leq 2\sin^2\left(\frac{1-t}{2}(\dist(x,z)-\dist(y,z))\right) \label{eqn:proof_5423} \\
    & \qquad +2\sin((1-t)\dist(x,z))\sin((1-t)\dist(y,z))\sin^2(\alpha(0)/2) \nonumber
\end{align}
where the last line makes use of the inequality $0\leq \alpha(t)/2\leq \alpha(0)/2\leq \pi/2$.

Applying the equality \eqref{eqn:proof_54231} with $t=0$ for the denominator and the inequality \eqref{eqn:proof_5423} for the numerator, and letting
$\beta_y:=\dist(x,z)-\dist(y,z)$, we obtain that 
\begin{align*}
    \, & \frac{1-\cos(\dist(\gamma(t),\eta(t)))}{1-\cos(\dist(x,y))} \\
    & \quad \quad \leq  \frac{\sin^2\left(\frac{1-t}{2}\beta_y\right)+\sin((1-t)\dist(x,z))\sin((1-t)\dist(y,z))\sin^2(\alpha(0)/2)}{\sin^2\left(\frac{1}{2}\beta_y\right)+\sin(\dist(x,z))\sin(\dist(y,z))\sin^2(\alpha(0)/2)}
\end{align*}

Since all summands are nonnegative, the right-hand side
is bounded by
\[
\max\left\{
\frac{\sin^2((1-t)\beta_y/2)}{\sin^2(\beta_y/2)},
\frac{\sin((1-t)\dist(x,z))\sin((1-t)\dist(y,z))}
     {\sin(\dist(x,z))\sin(\dist(y,z))}
\right\},
\]
where the first quotient is understood as $(1-t)^2$ when $\beta_y=0$. As $y\to x$, these quotients converge respectively to $(1-t)^2$ and $\sin^2((1-t)\dist(x,z))/\sin^2(\dist(x,z))$.
By concavity of $\sin$ on $[0,\pi]$,
\[
\sin((1-t)\dist(x,z))
\geq(1-t)\sin(\dist(x,z)),
\]
so the second limit is at least the first.
Taking the limsup therefore gives

\begin{equation} \label{eqn:proof98765}
    \limsup_{\substack{y\to x \\ y\in D\setminus\{x\}}} \frac{1-\cos(\dist(\gamma(t),\eta(t)))}{1-\cos(\dist(x,y))} \leq  \frac{\sin^2((1-t)\dist(x,z))}{\sin^2(\dist(x,z))}
\end{equation}

We conclude by noting that, as $y\to x$, both $\dist(x,y)$ and $\dist(\gamma(t),\eta(t))\to 0$ (geodesics vary continuously with their endpoints in $\kappa$-domains of $\CAT(\kappa)$ spaces, $\kappa\in\R$, see \cite[Proposition II.1.4]{bridsonhaefliger1999}), hence, the left-hand side of \eqref{eqn:proof98765} coincides with $\DS \limsup_{\substack{y\to x \\ y\in D\setminus\{x\}}}\frac{\dist(\gamma(t),\eta(t))^2}{\dist(x,y)^2}$, by $1-\cos u\sim u^2/2$ as $u\to0$.
Taking square roots proves the result.

\end{proof}

\subsection{Proof of the Chebyshev property in \texorpdfstring{\boldmath$\CAT(\kappa)$}{CAT(kappa)} spaces}
\label{subsec:proof_of_chebyshev}

\begin{figure}[htbp]
  \centering
  \includegraphics[width=0.60\textwidth]
    {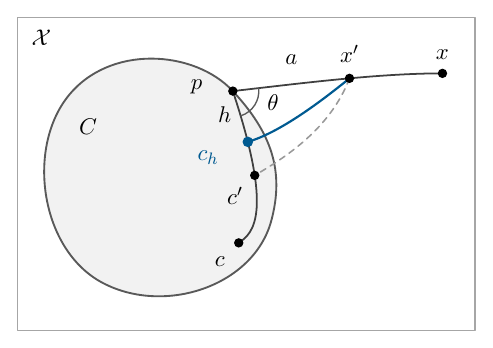}
  \caption{Configuration assumed for contradiction in the proof
    of the obtuse angle property.
    Here $p=p_C(x)$, $x'\in[p,x]$, and $c_h\in[p,c']\subset C$,
    with $a=\dd(p,x')$, $h=\dd(p,c_h)$, and
    $\theta=\angle_p(x,c)<\pi/2$.
    The blue geodesic joins $x'$ to $c_h$.}
  \label{fig:chebyshev-projection}
\end{figure}

\begin{proof}[Proof of Proposition~\ref{lem:projection} in the case $\kappa>0$]

We work in the normalized case $\kappa=1$; the general case follows by rescaling the metric to $\sqrt{\kappa}\,\dd$.

For the first part of the lemma, we provide a more straightforward proof than the one in \cite[Proposition II.2.4(1)--(3) and Exercise II.2.6(1)]{bridsonhaefliger1999}, \cite[Proposition 3.1(i)]{ohta2007}. Denote by $\rho:=\dist(x,C)$ and take an arbitrary $\varepsilon\in(0,\pi/2-\rho)$.
It is easy to see that finding a projection amounts to minimizing
\[
z\in C\cap B(x,\rho+\varepsilon)\longmapsto d^2(x,z).
\] 
The set $C \cap B(x,\rho+\varepsilon)$ is nonempty and convex, as an intersection of convex sets. It is also complete, since it is closed in the complete space $\X$. By spherical comparison, as in the proof of \cite[Proposition 3.1(i)]{ohta2007}, the fixed-center function $z \in B(x,\rho+\varepsilon) \mapsto d^2(x,z) $ is $\alpha(\varepsilon)-$strongly convex, with $\alpha(\varepsilon)>0$. Thus, since the function is continuous, we can apply \cite[Proposition 1.7]{sturm2003} to obtain a unique minimizer $p$ satisfying $d(x,p)=\rho$. Finally, any minimizer on $C$ belongs to
$C\cap B(x,\rho+\varepsilon)$, which proves the existence and uniqueness of the projection $p_C(x)=p$.

For completeness, let us now provide the proof of the obtuse angle property.
Assume $x\notin C$ and denote $p:=p_C(x)$. By contradiction, assume that there exists $c\in C\setminus\{p\}$ such that
\[
\angle_p(x,c)<\frac{\pi}{2}.
\]
Since the spherical comparison angles converge to the Alexandrov angle as both points approach $p$, we may choose $x'\in[p,x]\setminus\{p\}$ and $c'\in[p,c]\setminus\{p\}$ sufficiently close to $p$ such that the triangle $\Delta(p,x',c')$ has perimeter less than $2\pi$ and
\[
\alpha:=\angle_p^{(1)}(x',c')<\frac{\pi}{2}.
\]
Consider its spherical comparison triangle $\Delta(\bar p,\bar x',\bar c')$ and denote
$a:=\dd(p,x')\in(0,\pi/2)$. For $0<h<\dd(p,c')$, let $c_h\in[p,c']$ and $\bar c_h\in[\bar p,\bar c']$ be the corresponding points satisfying
\[
\dd(p,c_h)=\dd_1(\bar p,\bar c_h)=h.
\]
Since $C$ is geodesically convex, $c_h\in C$. For an illustration of this construction, see Figure~\ref{fig:chebyshev-projection}. Since $\bar c_h\in[\bar p,\bar c']$, the angle at $\bar p$ between $[\bar p,\bar x']$ and $[\bar p,\bar c_h]$ is $\alpha$.
By the spherical law of cosines, we get
\[
\cos\dd_1(\bar x',\bar c_h)
=\cos a\cos h+\sin a\sin h\cos\alpha.
\]
Hence,
\[
\begin{aligned}
\cos\dd_1(\bar x',\bar c_h)-\cos a
&=\cos a(\cos h-1)+\sin a\sin h\cos\alpha\\
&=\sin h\left[
\sin a\cos\alpha-\cos a\tan\left(\frac h2\right)
\right].
\end{aligned}
\]
Since $\sin a\cos\alpha>0$, the right-hand side is positive for sufficiently small $h>0$.
Therefore, $\dd_1(\bar x',\bar c_h)<a$. By the $\CAT(1)$ comparison and the triangle inequality,
\[
\begin{aligned}
\dd(x,c_h)
&\leq\dd(x,x')+\dd(x',c_h)\\
&\leq\dd(x,x')+\dd_1(\bar x',\bar c_h)\\
&<\dd(x,x')+\dd(x',p)=\dd(x,p).
\end{aligned}
\]
Since $c_h\in C$, this contradicts the definition of $p=p_C(x)$.
\end{proof}

\section{Omitted proofs from Section~\ref{sec:geometric_properties_Frechet}}

\subsection{Directional derivative of Fréchet functional}
\label{subsec:proof_dir_derivative}
\begin{proof}[Proof of Proposition~\ref{prop:dir-derivative}]
Suppose first that $\kappa>0$. The case $z=x$ is immediate, so assume $z\neq x$.For $y=x$, the squared-distance derivative is zero. For every $y\in D\setminus\{x\}$, we have $\dd(x,y),\dd(x,z)<D_\kappa$, so the first variation formula~\eqref{eq:first-variation} gives
\[
\lim_{t\downarrow0}
\frac{\frac12\dd^2(\gamma_{xz}(t),y)-\frac12\dd^2(x,y)}{t}
=
-\big\langle \log_x y,\log_x z\big\rangle_x.
\]
Since $\mu$ has bounded support, the difference quotients are uniformly bounded by an integrable constant. Hence, by the dominated convergence theorem,
\[
DF_\mu(x)[\log_x z]
=
-\int_{D}
\langle \log_x y,\log_x z\rangle_x
\,\dd\mu(y).
\]
To obtain the extension to the ambient tangent cone, fix $x\in D$ and consider any nonconstant geodesic $\eta$ issuing from $x$. Since $\dd(x,y)\leq2r<D_\kappa$ for every $y\in D$, the same first variation and dominated convergence argument applies along a sufficiently short initial segment of $\eta$. Allowing arbitrary nonnegative speeds therefore defines
\[
DF_\mu(x)[v]
=
-\int_D\langle\log_x y,v\rangle_x\,\dd\mu(y),
\qquad v\in\widetilde C_x\X.
\]
For every $\varepsilon>0$, the function $F_\mu$ is $(2r+\varepsilon)$-Lipschitz on $B(x,\varepsilon)$. Let $v,w\in\widetilde C_x\X$ be represented by geodesics
$\eta_v,\eta_w$ issuing from $x$, defined near $0$ with respective constant speeds $|v|_x$ and $|w|_x$. For sufficiently small $t>0$,
\[
\left|
\frac{F_\mu(\eta_v(t))-F_\mu(x)}{t}
-
\frac{F_\mu(\eta_w(t))-F_\mu(x)}{t}
\right|
\leq
(2r+\varepsilon)
\frac{\dd(\eta_v(t),\eta_w(t))}{t}.
\]
By the definition of the Alexandrov angle in Definition~\ref{def:alexandrov-angle} and the Euclidean comparison angle (see Remark~\ref{rem:euclidean-angle}),
\[
\lim_{t\downarrow0}
\frac{\dd(\eta_v(t),\eta_w(t))}{t}
=
\dd_{C_x\X}(v,w).
\]
Letting $t\downarrow0$ and then $\varepsilon\downarrow0$, we obtain
\[
|DF_\mu(x)[v]-DF_\mu(x)[w]|
\leq2r\,\dd_{C_x\X}(v,w).
\]
Thus, the extension theorem for uniformly continuous maps
\cite[Chapter~XI, Section~3, Theorem~1]{schwartz_analyse2}
gives a unique continuous extension of $DF_\mu(x)$ to $C_x\X$. Passing to the limit in the preceding inequality shows that this extension remains $2r$-Lipschitz.
The integral representation extends by continuity of the pseudo-inner product and dominated convergence.
 
If $\kappa\leq0$, the space is $\CAT(0)$ and the proof is similar, without a bounded-support assumption.
Indeed, for $\mu\in\mathcal P_2(\X)$, $\gamma\in\Gamma(x,z)$ and $h:=\dist(x,z)$,
\[
\left|
\frac{\dist^2(\gamma(t),y)-\dist^2(x,y)}{2t}
\right|
\leq h\,\dist(x,y)+\frac12h^2,
\qquad 0<t\leq1.
\]
The right-hand side is integrable, so dominated convergence gives the same directional derivative formula. The extension to $C_x\X$ follows as above, with $2r$
replaced by $M_x:=\int_\X\dist(x,y)\,\diff\mu(y)<\infty$: the function $F_\mu$ is $(M_x+\varepsilon)$-Lipschitz on $B(x,\varepsilon)$, and its directional derivative
therefore extends to an $M_x$-Lipschitz function on $C_x\X$.
\end{proof}

\subsection{Proof of the existence and uniqueness theorem}
\label{app:existence-uniqueness}

Let $(\X,\dd)$ be a $\CAT(\kappa)$ space, with $\kappa>0$, and let $\mu$ be a Borel probability measure supported on a $\kappa$-domain $D=B(o,r)$, where $r<D_\kappa/2$.
We first construct approximate minimizers using Ekeland's variational principle, then recall Yokota's localization lemmas and prove Theorem~\ref{thm:existence_uniqueness}.

\begin{lemma}[Ekeland's variational principle
{\cite[Theorem~1]{ekeland1979}}]
Let $(M,\dd)$ be a complete metric space, and let $f:M\to\R$ be lower semicontinuous and bounded below.
Let $\varepsilon>0$ and $x_0\in M$ satisfy
\[
f(x_0)<\inf_{x\in M}f(x)+\varepsilon.
\]
Then there exists $x_\varepsilon\in M$ such that
\[
f(x_\varepsilon)\leq f(x_0),
\qquad
\dd(x_0,x_\varepsilon)<1,
\]
and, for every $y\in M$,
\[
f(y)\geq f(x_\varepsilon)
-\varepsilon\,\dd(x_\varepsilon,y).
\]
\end{lemma}

We now apply Ekeland's variational principle to the Fréchet function associated with $\mu$.
Applying the principle on the whole space $\X$ yields the following corollary.

\begin{corollary}
There exists a family $(x_\varepsilon)_{\varepsilon>0}$ of points in $\X$ such that
\begin{equation}
\label{eq:def-xeps}
\lim_{\varepsilon\downarrow0}F_\mu(x_\varepsilon)
=
\inf_{x\in\X}F_\mu(x),
\end{equation}
and, for every $y\in\X$ and every $\varepsilon>0$,
\begin{equation}
\label{eq:ekeland-condition}
F_\mu(y)\geq F_\mu(x_\varepsilon)
-\varepsilon\,\dd(x_\varepsilon,y).
\end{equation}
\end{corollary}

We now study the behavior of this family as $\varepsilon\to0$.
Before proving its convergence, we record the two localization steps from \cite[Claims~12 and~13]{yokota2016}.
The first excludes points lying away from the support ball, while the second shows that moving towards its center strictly decreases the Fréchet function near its boundary.
Their role is to understand the behavior of the Fréchet function around the ball, in particular through its directional derivatives, which are also studied in \cite{bussfillmore2001,afsari2011}.

\begin{lemma}[Yokota's first localization
{\cite[Claim~12]{yokota2016}}]
\label{lem:first-localization}
For every $\alpha>0$, there exists $a_\alpha>0$ such that
\begin{equation}
\label{eq:first-localization}
F_\mu(x)>
\inf_{z\in B(o,r)}F_\mu(z)+a_\alpha
\qquad
\text{for every }x\notin B(o,r+\alpha).
\end{equation}
\end{lemma}

The directional derivative formula also applies at every $x\in B(o,D_\kappa/2)$, since
$\dist(x,y)\leq\dist(x,o)+r<D_\kappa$ for every $y\in D$.
This justifies the derivative used in the following lemma.

\begin{lemma}[Yokota's second localization
{\cite[Claim~13]{yokota2016}}]
\label{lem:second-localization}
For every $\rho\in(r,D_\kappa/2)$, there exists $a_\rho>0$ such that
\begin{equation}
\label{eq:negative-radial-derivative}
DF_\mu(x)\left[\frac{\log_x o}{\dd(x,o)}\right]
<-a_\rho
\end{equation}
for every
\[
x\in B(o,D_\kappa/2)\setminus B(o,\rho).
\]
\end{lemma}

The geometric idea of Lemma~\ref{lem:first-localization} is illustrated in Figure~\ref{fig:yokota-localization}.
In the normalized case $\kappa=1$, take $\alpha>0$ sufficiently small and consider $x$ such that $r+\alpha\leq\dd(x,o)<2r$.
Reflect $x$ along the radial geodesic through $o$ about the point $\widetilde{x}_0$ at distance $r+\alpha/2$ from $o$.
The resulting point $x'$ satisfies
\[
\dd(o,x')=2r+\alpha-\dd(o,x)\leq r,
\qquad
\dd(x',\widetilde{x}_0)=\dd(\widetilde{x}_0,x),
\]
where the latter equality is indicated by the matching marks in the figure.
The solid gray boundary represents the radius $r$, and the dashed boundary represents the radius $r+\alpha/2$.
For $y$ in the support ball, spherical comparison gives an acute angle between the inward radial direction and the segment towards $y$, indicated by the blue arc at $\widetilde{x}_0$.
The corresponding outward angle is obtuse.
The spherical law of cosines therefore shows that the blue geodesic from $y$ to $x'$ is shorter than the gray geodesic from $y$ to $x$.
The resulting uniform gap yields a decrease of $F_\mu$ after integration.
Configurations outside the comparison perimeter range, and points farther from $o$, are treated by side-length estimates, as in Yokota's proof.

For the second localization, choose $r<s<\rho$.
Since the measure is concentrated in the open ball of radius $s$, angle comparison and first variation give a uniformly negative inward derivative outside a ball of radius $s'<s<\rho$.
This yields the stated estimate.

\begin{corollary}
\label{cor:localization-approximate-minimizers}
For every $\rho\in(r,D_\kappa/2)$,
$x_\varepsilon\in B(o,\rho)$ for all sufficiently small $\varepsilon>0$.
\end{corollary}

\begin{proof}
We work in the normalized case $\kappa=1$; the general case follows by rescaling the metric to $\sqrt{\kappa}\,\dd$.
Fix $\rho\in(r,\pi/2)$ and choose $\alpha>0$ such that $\rho<r+\alpha<\pi/2$.
By \eqref{eq:def-xeps} and Yokota's first localization (Lemma~\ref{lem:first-localization}), we have $x_\varepsilon\in B(o,r+\alpha)$ for all sufficiently small $\varepsilon>0$.
Proposition~\ref{prop:dir-derivative}, applied in this enlarged ball, justifies the directional derivative.
If $x_\varepsilon=o$, the conclusion is immediate.
Otherwise, differentiating \eqref{eq:ekeland-condition} along the geodesic towards $o$ gives
\[
DF_\mu(x_\varepsilon)
\left[\frac{\log_{x_\varepsilon}o}{\dd(x_\varepsilon,o)}\right]
\geq-\varepsilon.
\]
For $\varepsilon<a_\rho$, this contradicts Yokota's second localization (Lemma~\ref{lem:second-localization}) if $x_\varepsilon\notin B(o,\rho)$.
Hence $x_\varepsilon\in B(o,\rho)$ for all sufficiently small $\varepsilon>0$.
\end{proof}

We now study the convergence of this family by showing that it is Cauchy.
In his proof, Yokota uses a so-called Kendall function (see \cite[Theorem~A]{yokota2016}), namely a geodesically convex function on $B(o,r)\times B(o,r)$, and controls $\Phi(x_\varepsilon,x_{\varepsilon'})$.
Since $\Phi$ is comparable to $\dd^\alpha$ for some $\alpha>2$, this yields the desired conclusion.

Our argument does not require Kendall's function: we can instead control directly the distance between the two points using Theorem~\ref{thm:directional-derivative-transfer}.
Together with the approximate first-order condition satisfied by $x_\varepsilon$, this inequality allows us to show that the family is Cauchy.

\begin{proof}[Proof of Theorem~\ref{thm:existence_uniqueness}]
We work in the normalized case $\kappa=1$; the general case follows by rescaling the metric to $\sqrt{\kappa}\,\dd$.
Fix $\rho\in(r,\pi/2)$.
By Corollary~\ref{cor:localization-approximate-minimizers}, $x_\varepsilon\in B(o,\rho)$ for all sufficiently small $\varepsilon>0$.
The directional derivatives are therefore well defined, and \eqref{eq:ekeland-condition} gives
\begin{equation}
\label{eq:diff-condition}
DF_\mu(x_\varepsilon)[\log_{x_\varepsilon}y]
\geq-\varepsilon\,\dd(x_\varepsilon,y),
\qquad y\in B(o,\rho).
\end{equation}

We apply Theorem~\ref{thm:directional-derivative-transfer} to $x_\varepsilon$ and $x_{\varepsilon'}$, using \eqref{eq:diff-condition} in the directions towards $x_{\varepsilon'}$ and $o$.
This yields
\begin{equation}
\label{eq:one-sided-eps}
\begin{aligned}
F_\mu(x_\varepsilon)-F_\mu(x_{\varepsilon'})
\leq{}&
\varepsilon\Bigl[
\sin\dd(x_\varepsilon,x_{\varepsilon'})
\\
&\quad+
\bigl(1-\cos\dd(x_\varepsilon,x_{\varepsilon'})\bigr)
\tan\dd(x_\varepsilon,o)
\Bigr]\\
&-
\bigl(1-\cos\dd(x_\varepsilon,x_{\varepsilon'})\bigr)
\frac{\cos\rho}{\cos\dd(x_\varepsilon,o)}.
\end{aligned}
\end{equation}

Interchanging $x_\varepsilon$ and $x_{\varepsilon'}$, we obtain the corresponding inequality with $\varepsilon'$ and $\dd(x_{\varepsilon'},o)$.
Adding the two inequalities and using $\tan\dd(x_\varepsilon,o)\leq\tan\rho$, $\tan\dd(x_{\varepsilon'},o)\leq\tan\rho$, and $1/\cos\dd(x_\varepsilon,o)\geq1$ (and similarly for $\varepsilon'$), we obtain
\[
\begin{aligned}
2\cos\rho\,
\bigl(1-\cos\dd(x_\varepsilon,x_{\varepsilon'})\bigr)
\leq{}&
(\varepsilon+\varepsilon')
\Bigl[\sin\dd(x_\varepsilon,x_{\varepsilon'})\\
&+
\bigl(1-\cos\dd(x_\varepsilon,x_{\varepsilon'})\bigr)
\tan\rho\Bigr].
\end{aligned}
\]
Since $\sin u\leq1$ and $1-\cos u\leq2$, we obtain
\[
0\leq
1-\cos\dd(x_\varepsilon,x_{\varepsilon'})
\leq
\frac{1+2\tan\rho}{2\cos\rho}
(\varepsilon+\varepsilon')
\longrightarrow0
\]
as $\varepsilon,\varepsilon'\downarrow0$.
Since $\dd(x_\varepsilon,x_{\varepsilon'})\leq2\rho<\pi$, it follows that $\dd(x_\varepsilon,x_{\varepsilon'})\to0$.
Thus, $(x_\varepsilon)_{\varepsilon>0}$ is a Cauchy family.
Completeness gives its limit $b(\mu)$.
Since $\rho>r$ was arbitrary, $b(\mu)\in B(o,r)$.
Continuity of $F_\mu$ and \eqref{eq:def-xeps} show that
$F_\mu(b(\mu))=\inf_\X F_\mu$.
Finally, any minimizer belongs to $B(o,r)$ by Lemma~\ref{lem:first-localization} and satisfies the first-order condition.
Applying the preceding summed inequality to two minimizers, with $\varepsilon=\varepsilon'=0$, therefore shows that their distance is zero.
\end{proof}

\subsection{Proof of dual norm identity on
the tangent cone}
\label{subsec:proof_dual_norm_identity}

\begin{proof}[Proof of Lemma~\ref{lem:dual-norm-tangent-cone}]
    Set $A(u):=\int_{C_b\X}\langle u,v\rangle_b\,\diff\nu(v)$ for $u\in C_b\X$.
    By the definition of the cone distance, for all $u\in C_b\X$,
    \[
    F_\nu(u)=F_\nu(0_b)+\frac12|u|_b^2-A(u),
    \]
    and the variance inequality
    \cite[Proposition~4.4]{sturm2003} at the barycenter $\beta$
    therefore gives
    \[
    A(u)\leq\langle u,\beta\rangle_b+K,
    \]
    where $K:=F_\nu(0_b)-F_\nu(\beta)-\frac12|\beta|_b^2$ does not depend on $u$.
    Scalar multiplication on $C_b\X$ is defined in Definition~\ref{def:tangent-cone}.
    Applying this inequality to $\lambda u$, dividing by $\lambda>0$, and letting $\lambda\to\infty$, homogeneity yields
    \[
    A(u)\leq\langle u,\beta\rangle_b\leq|u|_b|\beta|_b.
    \]
    Thus $\sup_{|u|_b\leq1}A(u)\leq|\beta|_b$. On the other hand, the same variance inequality, evaluated at $0_b$, gives
    \[
    A(\beta)
    =F_\nu(0_b)-F_\nu(\beta)+\frac12|\beta|_b^2
    \geq|\beta|_b^2.
    \]
    If $\beta\neq0_b$, taking $u=\beta/|\beta|_b$ and using homogeneity proves the reverse inequality. If $\beta=0_b$, equality follows by taking $u=0_b$.
\end{proof}

\section{Omitted proofs from Section~\ref{sec:consequences}}

\subsection{Proof of discrete recurrence result}
\label{subsec:discrete-recurrence}

\begin{proof}[Proof of Lemma~\ref{lem:discrete-recurrence}]
The bound holds for $n=1$ by the definition of $C$.
Suppose that $a_n\leq C/(n+A)$.
Since $n+A\geq\alpha$ and $(\alpha-1)C\geq\beta$,
\[
\begin{aligned}
a_{n+1}
&\leq
\left(1-\frac{\alpha}{n+A}\right)\frac{C}{n+A}
+\frac{\beta}{(n+A)^2}\\
&=
\frac{C}{n+A}
-\frac{\alpha C-\beta}{(n+A)^2}\\
&\leq
\frac{C}{n+A}-\frac{C}{(n+A)^2}\\
&\leq
\frac{C}{n+1+A}.
\end{aligned}
\]
This completes the induction.
\end{proof}

\bibliographystyle{alpha}
\bibliography{references}

\end{document}